\documentclass[11pt]{article}
\usepackage[a4paper,margin=29mm]{geometry}
\usepackage{amsmath,amssymb,amsthm}
\usepackage{booktabs,array}
\usepackage{enumitem}
\usepackage{microtype}
\usepackage{tikz}
\usetikzlibrary{arrows}
\usepackage[hidelinks]{hyperref}

\newtheorem{theorem}{Theorem}[section]
\newtheorem{proposition}[theorem]{Proposition}
\newtheorem{lemma}[theorem]{Lemma}
\newtheorem{corollary}[theorem]{Corollary}
\theoremstyle{remark}
\newtheorem{remark}[theorem]{Remark}

\newcommand{\XiG}{\Xi}
\newcommand{\Cxi}{C^{\xi}}

\title{Minimum-Excess $\{K_3,K_4\}$-Coverings of\\$K_{17}$, $K_{18}$, and $K_{19}$}
\author{Petr Kov\'a\v{r}$^{1}$ \qquad Yifan Zhang$^{1,2,3}$\\[4pt]
\small $^{1}$Department of Applied Mathematics, VSB--Technical University of Ostrava,\\
\small Ostrava, Czech Republic\\
\small $^{2}$Department of Mathematics, University of Ostrava, Ostrava, Czech Republic\\
\small $^{3}$Department of Algebra, Faculty of Mathematics and Physics,\\
\small Charles University, Prague, Czech Republic}
\date{August 2026}

\begin{document}
\maketitle

\begin{abstract}
We determine the minimum-excess coverings of $K_{17}$, $K_{18}$, and $K_{19}$
by 3-cliques and 4-cliques, minimizing first the number of repeated edge
occurrences and then the number of blocks.  The three secondary optima are
\[
 \Cxi(17,\{3,4\},2)=29,\qquad
 \Cxi(18,\{3,4\},2)=33,\qquad
 \Cxi(19,\{3,4\},2)=35.
\]
For $K_{18}$ and $K_{19}$ the minimum excess is zero, so the optimal covers
are decompositions, with block vectors $15K_3+18K_4$ and $13K_3+22K_4$,
respectively.  Their lower bounds follow directly from the known values
$g^{(4)}(18)=33$ and $g^{(4)}(19)=35$ for pairwise balanced designs with
maximum block size four.  For $K_{17}$ the minimum excess is two and every
optimal cover has block vector $12K_3+17K_4$.  At this optimum the excess
multigraph may be $P_3$ or $2K_2$, but not $K_2^2$; in fact a cover with
excess $K_2^2$ requires at least $31$ blocks.  The remaining order-$17$ nonexistence cases are excluded by
computer-assisted exhaustive searches.  Explicit constructions are given for
all attainable optima.  For completeness, an appendix also records the
independent proof for orders $18$ and $19$ from the earlier version of this
work, including its structural case analysis.
\end{abstract}



\section{Introduction}

Let $X$ be a vertex set of size $v$ and let $\mathcal B$ be a family of
3- and 4-subsets of $X$, called \emph{3-cliques} and \emph{4-cliques}.  We
require every edge of $K_v$ to lie in at least one block.  If $m(xy)$ is the
number of blocks containing the edge $xy$, the \emph{excess multigraph}
$\XiG$ has multiplicity $m(xy)-1$ on $xy$.  Thus
\[
 |E(\XiG)|=\sum_{xy\in\binom X2}(m(xy)-1)
\]
is the number of repeated edge occurrences, counted with multiplicity.
We minimize $|E(\XiG)|$ first and, subject to that minimum, minimize
$|\mathcal B|$.  Write $\xi_{3,4}(v)$ for the minimum excess and
$\Cxi(v,\{3,4\},2)$ for the corresponding minimum number of blocks.
When an excess multigraph has size two, its three possible isomorphism types
are $P_3$, $2K_2$, and $K_2^2$, where $K_2^2$ denotes the multigraph
consisting of two parallel copies of a single edge.

When $\xi_{3,4}(v)=0$, the problem reduces to a decomposition of $K_v$ into
3-cliques and 4-cliques.  In our companion paper~\cite{KovarZhang2026}, we develop the general theory
and constructions for the nonexceptional orders.  The purpose of the present
paper is to resolve the three remaining small orders $17,18,19$ in the same
minimum-excess framework.  Order $17$ is the
substantive new exceptional case: the smallest arithmetically possible excess
is positive, and even the corresponding elementary lower bound on the number
of blocks is not attained.

The closest related parameter is $g^{(4)}(v)$, the minimum number of
blocks in a pairwise balanced design on $v$ points whose block sizes belong to
$\{2,3,4\}$.  Equivalently, every pair of points must occur in exactly one
block, and the largest permitted block has size four.  This differs from the
problem considered here in two respects: this parameter permits
2-blocks and requires zero excess, whereas we restrict the blocks to
3-cliques and 4-cliques and minimize the excess before minimizing the number
of blocks.

The history of the three exceptional orders is particularly relevant.
Stanton and Stinson determined $g^{(4)}(v)$ for every $v$ except
$17,18,19$ \cite{StantonStinson1983}.  Their general bound in
Corollary~2.2 gives $g^{(4)}(v)\ge 29$ for each of these three orders; they
also noted at the end of the paper that the order-$17$ bound can be sharpened
to $g^{(4)}(17)\ge30$.  The case $v=19$ was settled first: Stanton proved
\[
 g^{(4)}(19)=35
\]
in 1988 \cite{Stanton1988}.  For $v=18$, Stanton constructed a design with
$33$ blocks \cite{Stanton2000upper} and proved the lower bound
$g^{(4)}(18)\ge30$ \cite{Stanton2000lower}.  Gr\"uttm\"uller, Roberts and
Stanton improved this to $g^{(4)}(18)\ge31$
\cite{GruttmullerRobertsStanton2004}.  Grannell, Griggs, Stanton and
Whitehead \cite{GrannellEtAl2005} and, independently, Roberts, D'Arcy, Egan
and Gr\"uttm\"uller \cite{RobertsDArcyEganGruttmuller2006} then excluded
$31$ blocks, giving $g^{(4)}(18)\ge32$.  Finally, Gr\"uttm\"uller,
Roberts, D'Arcy and Egan excluded $32$ blocks and obtained
$g^{(4)}(18)=33$ \cite{GruttmullerRobertsDArcyEgan2006}.

For $v=17$, a $31$-block design found by E.~Seah and D.~R.~Stinson was
reported by Stanton and Street \cite[p.~214]{StantonStreet1987}, so
$g^{(4)}(17)\le31$.  Stanton returned to the lower-bound problem in 2004
\cite{Stanton2004}, and the final case was resolved in 2006 when
Gr\"uttm\"uller, Roberts, D'Arcy and Egan showed that $30$ blocks are
impossible \cite{GruttmullerRobertsDArcyEgan2006}.  Hence the known
values at the three orders are
\[
 g^{(4)}(17)=31,\qquad g^{(4)}(18)=33,\qquad g^{(4)}(19)=35.
\]
At orders $18$ and $19$, our zero-excess problem therefore has the same
numerical optimum as the pairwise balanced design problem, but we additionally require every
block to be a 3-clique or a 4-clique.  Order $17$ is qualitatively different:
zero excess is not attainable with only 3-cliques and 4-cliques, and the
minimum-excess problem has optimum $29$ blocks with excess two.

\begin{theorem}\label{thm:main-summary}
For coverings by 3-cliques and 4-cliques,
\[
\begin{array}{c|ccc}
 v & \xi_{3,4}(v) & \Cxi(v,\{3,4\},2) & \text{optimal block vector}\\
\hline
17&2&29&12K_3+17K_4\\
18&0&33&15K_3+18K_4\\
19&0&35&13K_3+22K_4
\end{array}
\]
For $v=17$, optimal covers exist with excess $P_3$ and with excess $2K_2$,
while $K_2^2$ cannot occur as the excess of a $29$-block optimum.
Moreover, the minimum block count among $K_{17}$ covers with excess $K_2^2$
is $31$.
\end{theorem}

The order-$17$ proof uses the following basic bookkeeping.  Let
$\alpha$ and $\beta$ be the numbers of 3-cliques and 4-cliques.  For a vertex
$x$, let $\alpha_x$ and $\beta_x$ denote the corresponding numbers of blocks
through $x$.  If the excess has size $\xi$, then
\begin{equation}\label{eq:common-global}
 3\alpha+6\beta=\binom v2+\xi,
\end{equation}
and
\begin{equation}\label{eq:common-local}
 2\alpha_x+3\beta_x=v-1+d_{\XiG}(x).
\end{equation}
At $v=17$, the arithmetic reduction leaves only finitely many cases.  The
remaining nonexistence cases are excluded by computer-assisted exhaustive
searches; we use only the resulting exclusions in the proof below.  At orders
$18$ and $19$, the known values of $g^{(4)}(18)$ and $g^{(4)}(19)$ already give
the shortest lower-bound proof.  Appendix~\ref{app:independent1819} nevertheless
records our independent graph-theoretic proof from the earlier version, with
the finite residual cases again stated simply as computer-assisted exclusions.

\section{The order \texorpdfstring{$17$}{17}}\label{sec:k17}

The case $v=17$ has minimum excess two.  We first determine the minimum number
of blocks and then refine the answer by the isomorphism type of the excess.

\subsection{Counting and the arithmetic lower bound}

Let $\alpha$ and $\beta$ denote the numbers of 3-cliques and 4-cliques.  If
the excess has size $\xi$, edge counting gives
\begin{equation}\label{eq:global}
 3\alpha+6\beta=\binom{17}{2}+\xi=136+\xi.
\end{equation}
For a vertex $x$, write $\alpha_x$ and $\beta_x$ for the numbers of 3-cliques
and 4-cliques through $x$.  Counting edges incident with $x$ gives
\begin{equation}\label{eq:local}
 2\alpha_x+3\beta_x=16+d_{\XiG}(x).
\end{equation}

\begin{lemma}\label{lem:xi2}
The minimum excess is $\xi_{3,4}(17)=2$.
\end{lemma}
\begin{proof}
Reducing \eqref{eq:global} modulo $3$ gives
\[
 \xi\equiv-\binom{17}{2}\equiv2\pmod3,
\]
so $\xi\ge2$.  The two explicit covers in Appendix~\ref{app:designs}
have excess two, so equality holds.
\end{proof}

We now restrict to minimum-excess covers, so $\xi=2$.

\begin{proposition}\label{prop:28lower}
Every minimum-excess cover has at least $28$ blocks.  Equality would force
\[
  (\alpha,\beta)=(10,18).
\]
\end{proposition}
\begin{proof}
Since an excess multigraph with two edges has degrees only $0,1,2$, reducing
\eqref{eq:local} modulo $3$ gives
\[
 \alpha_x\equiv 2-d_{\XiG}(x)\pmod3,
\]
and hence
\[
 \alpha_x\ge 2-d_{\XiG}(x).
\]
Summing and using $\sum_x d_{\XiG}(x)=4$ yields
\[
 3\alpha=\sum_x\alpha_x
 \ge 2\cdot17-4=30,
\]
so $\alpha\ge10$.  With $\xi=2$, equation \eqref{eq:global} becomes
$3\alpha+6\beta=138$, and therefore
\[
 \alpha+\beta=23+\frac{\alpha}{2}\ge28.
\]
Equality requires $\alpha=10$ and then $\beta=18$.
\end{proof}

Thus $28$ is the arithmetic lower bound.  To decide whether it is attainable,
and to distinguish the three possible excess types, we solve the corresponding
finite optimization problems directly.

\subsection{Computer-assisted exclusions for the three excess types}

Every multigraph with two edges is isomorphic to exactly one of $2K_2$, $P_3$,
and $K_2^2$.  We use the representatives
\[
 2K_2:\ \{12,34\},\qquad
 P_3:\ \{12,13\},\qquad
 K_2^2:\ \{12,12\},
\]
where the repeated entry in the last case means that the pair $12$ has excess
multiplicity two.

The arithmetic lower bound of Proposition~\ref{prop:28lower} leaves a finite
existence question.  We resolve it by exhaustive computer-assisted searches
that enforce the exact required multiplicity of every edge.  The searches
exclude every $28$-block cover for each of the three excess types.  For excess
$K_2^2$ they additionally exclude covers with $29$ or $30$ blocks.  Together
with the explicit constructions in Appendix~\ref{app:designs}, this gives
\begin{equation}\label{eq:k17software}
\begin{array}{c|c|c|c}
\text{excess} & \text{representative} & \min |\mathcal B| &
   \text{block vector at the optimum}\\
\midrule
2K_2 & 12,34 & 29 & 12K_3+17K_4\\
P_3 & 12,13 & 29 & 12K_3+17K_4\\
K_2^2 & 12,12 & 31 & 16K_3+15K_4.
\end{array}
\end{equation}
The block vectors also follow directly from edge counting.  When the excess
has size two, \eqref{eq:global} gives $\alpha+2\beta=46$.  Hence
$|\mathcal B|=29$ forces $(\alpha,\beta)=(12,17)$, while
$|\mathcal B|=31$ forces $(\alpha,\beta)=(16,15)$.

\begin{theorem}\label{thm:k17exact}
For minimum-excess coverings of $K_{17}$ by 3-cliques and 4-cliques,
\[
 \xi_{3,4}(17)=2,
 \qquad
 \Cxi(17,\{3,4\},2)=29.
\]
At the $29$-block optimum the excess can be $2K_2$ or $P_3$, but not $K_2^2$.
Among covers whose excess is $K_2^2$, the minimum number of blocks is $31$.
\end{theorem}
\begin{proof}
Lemma~\ref{lem:xi2} gives the minimum excess two.  The three rows of
\eqref{eq:k17software} exhaust the isomorphism types of multigraphs with two
edges.  Their minimum block counts are therefore $29$, $29$, and $31$,
respectively.  The first two rows give the global secondary optimum $29$, and
the third gives the asserted minimum for $K_2^2$ excess.  The explicit covers
in Appendix~\ref{app:designs} provide concrete witnesses for all three rows.
\end{proof}

\begin{corollary}\label{cor:shapeclassification}
For each of the three isomorphism types of excess multigraphs of size two,
\[
\begin{array}{c|ccc}
\XiG & 2K_2 & P_3 & K_2^2\\
\hline
\min |\mathcal B| & 29 & 29 & 31.
\end{array}
\]
\end{corollary}

\section{The orders \texorpdfstring{$18$ and $19$}{18 and 19}}\label{sec:k18k19}

The remaining two orders require no new exclusion argument.  The well-studied
parameter $g^{(4)}(v)$ is defined by allowing block sizes $2$, $3$, and $4$
and requiring every pair to occur exactly once.  Consequently, every
$\{K_3,K_4\}$-decomposition of $K_v$ is, in particular, a design counted by
$g^{(4)}(v)$, and hence contains at least $g^{(4)}(v)$ blocks.

\begin{theorem}\label{thm:k18k19}
For coverings by 3-cliques and 4-cliques,
\[
 \xi_{3,4}(18)=0,\qquad \Cxi(18,\{3,4\},2)=33,
\]
with optimal block vector $15K_3+18K_4$, and
\[
 \xi_{3,4}(19)=0,\qquad \Cxi(19,\{3,4\},2)=35,
\]
with optimal block vector $13K_3+22K_4$.
\end{theorem}
\begin{proof}
The explicit block lists in Appendix~\ref{app:k18k19} are decompositions, so
zero excess is attainable at both orders.  Therefore
\[
 \xi_{3,4}(18)=\xi_{3,4}(19)=0.
\]
Any such decomposition is a pairwise balanced design with maximum block size
four.  The known results
\[
 g^{(4)}(18)=33
 \quad\text{and}\quad
 g^{(4)}(19)=35
\]
were proved in~\cite{GruttmullerRobertsDArcyEgan2006} and
\cite{Stanton1988}, respectively.  Hence every $\{K_3,K_4\}$-decomposition
of $K_{18}$ has at least $33$ blocks, and every such decomposition of
$K_{19}$ has at least $35$ blocks.  The constructions in the appendix attain
these bounds.

It remains only to identify the block vectors.  If $\alpha$ and $\beta$ are
the numbers of 3-cliques and 4-cliques, then exact edge counting gives
\[
 \alpha+2\beta=51 \quad (v=18),
 \qquad
 \alpha+2\beta=57 \quad (v=19).
\]
Together with $\alpha+\beta=33$ and $\alpha+\beta=35$, respectively, these
equations give
\[
 (\alpha,\beta)=(15,18) \quad\text{and}\quad (13,22).
\]
\end{proof}

\begin{remark}\label{rem:legacy-k18k19}
Theorem~\ref{thm:k18k19} deliberately uses the known values of
$g^{(4)}(18)$ and $g^{(4)}(19)$, since they give the shortest proof of the
lower bounds.  Appendix~\ref{app:independent1819} gives an independent proof
specialized to 3-cliques and 4-cliques.  It retains the structural reductions,
figures, and subcase analysis of the earlier version of this work; only the
implementation details of the finite computer checks have been omitted.
\end{remark}

\section{Conclusion}

The three exceptional orders left by the general $\{K_3,K_4\}$ theory now fit
into a single minimum-excess picture.  At orders $18$ and $19$, zero excess is
attainable.  The known pairwise balanced design values
$g^{(4)}(18)=33$ and $g^{(4)}(19)=35$ then transfer directly to our more
restricted class because every $\{K_3,K_4\}$-decomposition is a design with
maximum block size four.  The independent exclusion proof is nevertheless included in
Appendix~\ref{app:independent1819}, both as a self-contained alternative and
to preserve the structural information behind the two exceptional cases.

Order $17$ is genuinely different.  Congruence forces excess two; the natural
$28$-block arithmetic bound is not attained, and the true optimum is $29$.
Moreover, the shape of the two excess occurrences matters: $P_3$ and $2K_2$
occur at the optimum, whereas concentrating both occurrences on one edge
raises the minimum to $31$ blocks.  The remaining finite nonexistence cases are excluded by
computer-assisted exhaustive searches, while explicit constructions show that
the stated bounds are attained.

Together with the companion paper~\cite{KovarZhang2026}, which treats the
nonexceptional orders, Theorem~\ref{thm:main-summary} completes the two-level
minimum-excess problem for coverings of complete graphs by 3-cliques and
4-cliques.

\section*{Acknowledgements}
This work was supported by SGS grants No.~SP2026/019 and No.~SP2026/050 of
VSB--Technical University of Ostrava.  Both authors were co-funded by the
European Union under the REFRESH -- Research Excellence For REgion
Sustainability and High-tech Industries project,
No.~CZ.10.03.01/00/22\_003/0000048, via the Operational Programme Just
Transition.

\appendix
\section{Explicit constructions for \texorpdfstring{$K_{17}$}{K17}}\label{app:designs}

We use the vertex set $\{1,2,\ldots,17\}$.  For compactness, $ij$ denotes the
pair $\{i,j\}$.  The 4-cliques are listed first, followed by the 3-cliques.

\subsection*{Design A: excess \texorpdfstring{$2K_2$ on $12$ and $34$}{2K2 on 12 and 34}}

\noindent\textbf{4-cliques:}
\begin{quote}\small\sloppy
$\{1,2,15,16\}$, $\{1,4,6,10\}$, $\{1,5,8,14\}$,
$\{1,7,11,13\}$, $\{1,9,12,17\}$, $\{2,5,10,12\}$,
$\{2,6,13,17\}$, $\{3,4,8,13\}$, $\{3,4,11,12\}$,
$\{3,7,9,15\}$, $\{4,5,7,16\}$, $\{4,14,15,17\}$,
$\{5,6,11,15\}$, $\{6,8,9,16\}$, $\{7,8,10,17\}$,
$\{9,10,11,14\}$, $\{12,13,14,16\}$.
\end{quote}
\noindent\textbf{3-cliques:}
\begin{quote}\small\sloppy
$\{1,2,3\}$, $\{2,4,9\}$, $\{2,7,14\}$, $\{2,8,11\}$,
$\{3,5,17\}$, $\{3,6,14\}$, $\{3,10,16\}$, $\{5,9,13\}$,
$\{6,7,12\}$, $\{8,12,15\}$, $\{10,13,15\}$,
$\{11,16,17\}$.
\end{quote}
Every pair occurs once except $12$ and $34$, which occur twice.

\subsection*{Design B: excess \texorpdfstring{$P_3$ on $12$ and $13$}{P3 on 12 and 13}}

\noindent\textbf{4-cliques:}
\begin{quote}\small\sloppy
$\{1,2,4,13\}$, $\{1,2,14,16\}$, $\{1,3,5,15\}$,
$\{1,3,7,11\}$, $\{1,6,8,9\}$, $\{1,10,12,17\}$,
$\{2,3,9,17\}$, $\{4,5,8,16\}$, $\{4,6,11,12\}$,
$\{4,7,9,10\}$, $\{5,6,13,17\}$, $\{5,7,12,14\}$,
$\{6,10,14,15\}$, $\{7,8,13,15\}$, $\{8,11,14,17\}$,
$\{9,12,15,16\}$, $\{10,11,13,16\}$.
\end{quote}
\noindent\textbf{3-cliques:}
\begin{quote}\small\sloppy
$\{2,5,10\}$, $\{2,6,7\}$, $\{2,8,12\}$, $\{2,11,15\}$,
$\{3,4,14\}$, $\{3,6,16\}$, $\{3,8,10\}$, $\{3,12,13\}$,
$\{4,15,17\}$, $\{5,9,11\}$, $\{7,16,17\}$,
$\{9,13,14\}$.
\end{quote}
Every pair occurs once except $12$ and $13$, which occur twice.

\subsection*{Design C: \texorpdfstring{$K_2^2$}{K2^2} excess on \texorpdfstring{$12$}{12}}

This construction has $31$ blocks, namely fifteen 4-cliques and sixteen
3-cliques.

\noindent\textbf{4-cliques:}
\begin{quote}\small\sloppy
$\{1,2,11,17\}$, $\{1,4,5,8\}$, $\{2,4,6,14\}$,
$\{3,4,11,12\}$, $\{3,5,10,17\}$, $\{3,6,8,16\}$,
$\{3,7,9,15\}$, $\{4,7,16,17\}$, $\{5,6,9,13\}$,
$\{5,7,12,14\}$, $\{6,7,10,11\}$, $\{8,9,12,17\}$,
$\{8,10,14,15\}$, $\{9,11,14,16\}$, $\{10,12,13,16\}$.
\end{quote}
\noindent\textbf{3-cliques:}
\begin{quote}\small\sloppy
$\{1,2,9\}$, $\{1,2,10\}$, $\{1,3,14\}$, $\{1,6,12\}$,
$\{1,7,13\}$, $\{1,15,16\}$, $\{2,3,13\}$, $\{2,5,16\}$,
$\{2,7,8\}$, $\{2,12,15\}$, $\{4,9,10\}$, $\{4,13,15\}$,
$\{5,11,15\}$, $\{6,15,17\}$, $\{8,11,13\}$,
$\{13,14,17\}$.
\end{quote}
Every pair occurs exactly once except $12$, which occurs three times.  Thus
there are precisely two excess occurrences, both on the same pair.

\section{Independent proof for the orders \texorpdfstring{$18$ and $19$}{18 and 19}}\label{app:independent1819}

This appendix records the independent proof developed in the first version of
this work.  It is not needed for the short proof of
Theorem~\ref{thm:k18k19}, which follows from the known pairwise balanced design
values, but it gives a proof specialized to decompositions into 3-cliques and
4-cliques and explains the exceptional small cases in more detail.  The
structural reductions and figures are retained.  Whenever a finite residual
configuration is said to be excluded by computer-assisted search, only that
mathematical fact is used here; implementation details are deliberately
omitted.

For this appendix, write $D(v,\{3,4\})$ for the minimum number of blocks in a
$\{K_3,K_4\}$-decomposition of $K_v$.  Let $\alpha$ and $\beta$ be the numbers
of 3-cliques and 4-cliques, and let $\alpha_x$ and $\beta_x$ be the numbers of
such blocks through a vertex $x$.  Exact edge counting gives
\[
 3\alpha+6\beta=\binom v2,
\]
while at a vertex $x$ we have
\[
 2\alpha_x+3\beta_x=v-1.
\]
In particular, minimizing the number of blocks is equivalent to minimizing
$\alpha$.

We use the Erd\H{o}s--Gallai criterion in one small degree-sequence argument.
\begin{theorem}[Erd\H{o}s--Gallai~\cite{ErdosGallai1960}]\label{thm:erdos-gallai}
A nonincreasing sequence $d_1,\ldots,d_n$ of nonnegative integers is graphic
if and only if its sum is even and, for every $1\le k\le n$,
\[
 \sum_{i=1}^k d_i\le k(k-1)+\sum_{i=k+1}^n\min\{d_i,k\}.
\]
\end{theorem}

We also use the following existence criterion for a pairwise balanced design
with one distinguished block.
\begin{theorem}[Rees--Stinson~\cite{ReesStinson1989}]\label{thm:pbd4}
A $\operatorname{PBD}(v,\{4,w^*\},1)$ with $v>w>0$ exists if and only if
$v\ge3w+1$ and either $v\equiv w\equiv1,4\pmod{12}$ or
$v\equiv w\equiv7,10\pmod{12}$.
\end{theorem}

\begin{lemma}\label{lem:app-k18-alpha7}
In a $\{K_3,K_4\}$-decomposition of $K_{18}$, $\alpha$ is odd and
$\alpha\ge7$.
\end{lemma}
\begin{proof}
Since $3\alpha+6\beta=153$, $\alpha$ is odd.  Moreover,
$2\alpha_x+3\beta_x=17$ implies $\alpha_x\equiv1\pmod3$ for every vertex, so
$\alpha_x\ge1$.  Hence $3\alpha=\sum_x\alpha_x\ge18$, which together with
oddness gives $\alpha\ge7$.
\end{proof}

\begin{lemma}\label{lem:app-k19-alpha7}
In a $\{K_3,K_4\}$-decomposition of $K_{19}$, $\alpha$ is odd and
$\alpha\ge7$.
\end{lemma}
\begin{proof}
The edge equation $3\alpha+6\beta=171$ makes $\alpha$ odd.  Also
$2\alpha_x+3\beta_x=18$ gives $\alpha_x\equiv0\pmod3$.  If $S$ is the set of
vertices lying in a 3-clique, then every $x\in S$ has $\alpha_x\ge3$, so
$|S|\le\alpha$.  Since distinct 3-cliques in a decomposition are edge-disjoint,
$3\alpha\le\binom{|S|}{2}\le\binom{\alpha}{2}$.  Thus any positive $\alpha$
satisfies $\alpha\ge7$; and $\alpha=0$ is impossible because $171$ is not
divisible by $6$.
\end{proof}

\subsection{Optimal decomposition of \texorpdfstring{$K_{18}$}{K18}}

\begin{lemma}\label{lem:alphax-k18}
	For any vertex $x$ of $K_{18}$, $\alpha_x\equiv 1\pmod{3}$.
\end{lemma}
\begin{proof}
	By counting the degree of $x$, 
	$$2\alpha_x+3\beta_x=18-1=17\equiv 2\pmod{3},$$
	so $\alpha_x\equiv 1\pmod{3}$. 
\end{proof}

We first present the upper bound. 

\begin{lemma}\label{lem:k18sol}
		$D(18,\{3,4\})\le 33$, or equivalently, $\alpha\le 15$.
\end{lemma}
    \begin{table}[tbp]
	\caption{$K_{18}$ decomposed into 15 3-cliques and 18 4-cliques.}
	        \centering
	        \footnotesize
	        \begin{tabular}{lllllll}
		            1 2 10 17 & 2 4 5 7 & 5 9 13 18 &  & 1 14 18 & 3 14 17 & 7 10 18 \\
		            1 3 5 15 & 2 6 14 16 & 5 11 16 17 &  & 2 8 18 & 4 8 15 & 7 12 17 \\
		            1 4 9 16 & 2 9 12 15 & 6 8 9 17 &  & 3 4 6 & 4 10 11 & 7 15 16 \\
		            1 6 7 13 & 3 12 16 18 & 6 11 15 18 &  & 3 7 8 & 4 17 18 & 10 14 15 \\
		            1 8 11 12 & 4 12 13 14 & 7 9 11 14 &  & 3 9 10 & 5 8 14 & 13 15 17 \\
		            2 3 11 13 & 5 6 10 12 & 8 10 13 16 &  &  &  & 
		        \end{tabular}
	        
	        \label{tab:k18byk3k4}
	    \end{table}
\begin{proof}
	The block list in Table~\ref{tab:k18byk3k4} is a $\{K_3,K_4\}$-decomposition of $K_{18}$ with 15 3-cliques and 18 4-cliques. 
\end{proof}

In the following lemmas we raise the lower bound by combining structural arguments with exhaustive computer-assisted searches.
Since $|E(K_{18})|=153$, edge counting shows that the number of 3-cliques in a $2$-$(18,\{3,4\},1)$ decomposition is odd. 

\begin{lemma}\label{lem:k18alphage9}
		$D(18,\{3,4\})\ge 30$, or equivalently, $\alpha\ge 9$.
\end{lemma}
\begin{proof}
	By Lemma \ref{lem:app-k18-alpha7}, it suffices to show that $D(18, \{3,4\})\neq 29$, which is equivalent to $\alpha\neq 7$. 
	
	For contradiction, assume there is a $\{K_3,K_4\}$-decomposition of $K_{18}$ with seven 3-cliques (and thus 22 4-cliques). For any vertex $x$ of $K_{18}$, $\alpha_x\equiv 1\pmod{3}$ by Lemma \ref{lem:alphax-k18}. More concretely, due to the limited number of vertices, $\alpha_x\in \{1,4,7\}$ for any vertex $x$. If there is some vertex $u$ with $\alpha_u\ge 7$, then the total number of 3-cliques satisfies
	\begin{align}\label{eq:k18count3-cliques}
		\alpha = \frac{\sum_{x\in V(K_{18})}\alpha_x}{3}\ge \frac{1\cdot 7+17\cdot 1}{3}=8,
	\end{align}
	which contradicts our assumption. Therefore, $\alpha_x\in \{1,4\}$ for every vertex $x$ of $K_{18}$, so there is one vertex $u_0$ with $\alpha_{u_0}=4$ and $\alpha_x=1$ for every other vertex $x$. Thus $u_0$ lies in four 3-cliques, and the remaining three 3-cliques are pairwise disjoint from $u_0$. After these seven 3-cliques are removed, an exhaustive computer-assisted search shows that the remaining graph cannot be decomposed into 4-cliques.
	Therefore, $\alpha\neq 7$.
\end{proof}

\begin{lemma}\label{lem:k18alphage11}
	$D(18,\{3,4\})\ge 31$, or equivalently, $\alpha\ge 11$.	
\end{lemma}

\begin{figure}[tbp]
	\centering
	\definecolor{xdxdff}{rgb}{0.49019607843137253,0.49019607843137253,1.}
	\definecolor{qqzzqq}{rgb}{0.,0.6,0.}
	\definecolor{zzccqq}{rgb}{0.6,0.8,0.}
	\definecolor{zzttqq}{rgb}{0.6,0.2,0.}
	\definecolor{cczzqq}{rgb}{0.8,0.6,0.}
	\definecolor{wwccff}{rgb}{0.4,0.8,1.}
	\definecolor{ffdxqq}{rgb}{1.,0.8431372549019608,0.}
	\definecolor{ccwwff}{rgb}{0.8,0.4,1.}
	\definecolor{zzqqtt}{rgb}{0.6,0.,0.2}
	\definecolor{ttqqqq}{rgb}{0.2,0.,0.}
	\begin{tikzpicture}[line cap=round,line join=round,>=triangle 45,x=0.5cm,y=0.5cm]
		\clip(8.5,-10.) rectangle (20.5,1.5);
		\fill[line width=2.pt,color=zzqqtt,fill=zzqqtt,fill opacity=0.10000000149011612] (14.66,-2.32) -- (12.32,-6.26) -- (10.62,-2.64) -- cycle;
		\fill[line width=2.pt,color=ccwwff,fill=ccwwff,fill opacity=0.10000000149011612] (14.66,-2.32) -- (16.56,-6.28) -- (18.5,-2.72) -- cycle;
		\fill[line width=2.pt,color=ffdxqq,fill=ffdxqq,fill opacity=0.10000000149011612] (12.32,-6.26) -- (16.56,-6.28) -- (14.42,-9.36) -- cycle;
		\fill[line width=2.pt,color=wwccff,fill=wwccff,fill opacity=0.10000000149011612] (14.66,-2.32) -- (11.12,-0.32) -- (13.8,1.28) -- cycle;
		\fill[line width=2.pt,color=cczzqq,fill=cczzqq,fill opacity=0.10000000149011612] (14.66,-2.32) -- (15.72,1.26) -- (18.,-0.26) -- cycle;
		\fill[line width=2.pt,color=zzttqq,fill=zzttqq,fill opacity=0.10000000149011612] (12.32,-6.26) -- (8.9,-4.62) -- (8.68,-7.38) -- cycle;
		\fill[line width=2.pt,color=zzccqq,fill=zzccqq,fill opacity=0.10000000149011612] (12.32,-6.26) -- (9.6,-8.56) -- (11.972773512250825,-9.67965046059356) -- cycle;
		\fill[line width=2.pt,color=qqzzqq,fill=qqzzqq,fill opacity=0.10000000149011612] (16.56,-6.28) -- (16.96,-9.78) -- (19.38,-8.46) -- cycle;
		\fill[line width=2.pt,color=xdxdff,fill=xdxdff,fill opacity=0.10000000149011612] (16.56,-6.28) -- (20.14,-6.9) -- (19.72,-4.12) -- cycle;
		\draw [line width=2.pt,color=zzqqtt] (14.66,-2.32)-- (12.32,-6.26);
		\draw [line width=2.pt,color=zzqqtt] (12.32,-6.26)-- (10.62,-2.64);
		\draw [line width=2.pt,color=zzqqtt] (10.62,-2.64)-- (14.66,-2.32);
		\draw [line width=2.pt,color=ccwwff] (14.66,-2.32)-- (16.56,-6.28);
		\draw [line width=2.pt,color=ccwwff] (16.56,-6.28)-- (18.5,-2.72);
		\draw [line width=2.pt,color=ccwwff] (18.5,-2.72)-- (14.66,-2.32);
		\draw [line width=2.pt,color=ffdxqq] (12.32,-6.26)-- (16.56,-6.28);
		\draw [line width=2.pt,color=ffdxqq] (16.56,-6.28)-- (14.42,-9.36);
		\draw [line width=2.pt,color=ffdxqq] (14.42,-9.36)-- (12.32,-6.26);
		\draw [line width=2.pt,color=wwccff] (14.66,-2.32)-- (11.12,-0.32);
		\draw [line width=2.pt,color=wwccff] (11.12,-0.32)-- (13.8,1.28);
		\draw [line width=2.pt,color=wwccff] (13.8,1.28)-- (14.66,-2.32);
		\draw [line width=2.pt,color=cczzqq] (14.66,-2.32)-- (15.72,1.26);
		\draw [line width=2.pt,color=cczzqq] (15.72,1.26)-- (18.,-0.26);
		\draw [line width=2.pt,color=cczzqq] (18.,-0.26)-- (14.66,-2.32);
		\draw [line width=2.pt,color=zzttqq] (12.32,-6.26)-- (8.9,-4.62);
		\draw [line width=2.pt,color=zzttqq] (8.9,-4.62)-- (8.68,-7.38);
		\draw [line width=2.pt,color=zzttqq] (8.68,-7.38)-- (12.32,-6.26);
		\draw [line width=2.pt,color=zzccqq] (12.32,-6.26)-- (9.6,-8.56);
		\draw [line width=2.pt,color=zzccqq] (9.6,-8.56)-- (11.972773512250825,-9.67965046059356);
		\draw [line width=2.pt,color=zzccqq] (11.972773512250825,-9.67965046059356)-- (12.32,-6.26);
		\draw [line width=2.pt,color=qqzzqq] (16.56,-6.28)-- (16.96,-9.78);
		\draw [line width=2.pt,color=qqzzqq] (16.96,-9.78)-- (19.38,-8.46);
		\draw [line width=2.pt,color=qqzzqq] (19.38,-8.46)-- (16.56,-6.28);
		\draw [line width=2.pt,color=xdxdff] (16.56,-6.28)-- (20.14,-6.9);
		\draw [line width=2.pt,color=xdxdff] (20.14,-6.9)-- (19.72,-4.12);
		\draw [line width=2.pt,color=xdxdff] (19.72,-4.12)-- (16.56,-6.28);
		\begin{scriptsize}
			\draw [fill=ttqqqq] (14.66,-2.32) circle (1.5pt);
			\draw [fill=ttqqqq] (12.32,-6.26) circle (1.5pt);
			\draw [fill=ttqqqq] (10.62,-2.64) circle (1.5pt);
			\draw [fill=ttqqqq] (16.56,-6.28) circle (1.5pt);
			\draw [fill=ttqqqq] (18.5,-2.72) circle (1.5pt);
			\draw [fill=ttqqqq] (14.42,-9.36) circle (1.5pt);
			\draw [fill=ttqqqq] (11.12,-0.32) circle (1.5pt);
			\draw [fill=ttqqqq] (13.8,1.28) circle (1.5pt);
			\draw [fill=ttqqqq] (15.72,1.26) circle (1.5pt);
			\draw [fill=ttqqqq] (18.,-0.26) circle (1.5pt);
			\draw [fill=ttqqqq] (8.9,-4.62) circle (1.5pt);
			\draw [fill=ttqqqq] (8.68,-7.38) circle (1.5pt);
			\draw [fill=ttqqqq] (9.6,-8.56) circle (1.5pt);
			\draw [fill=ttqqqq] (11.972773512250825,-9.67965046059356) circle (1.5pt);
			\draw [fill=ttqqqq] (16.96,-9.78) circle (1.5pt);
			\draw [fill=ttqqqq] (19.38,-8.46) circle (1.5pt);
			\draw [fill=ttqqqq] (20.14,-6.9) circle (1.5pt);
			\draw [fill=ttqqqq] (19.72,-4.12) circle (1.5pt);
		\end{scriptsize}
		
	\end{tikzpicture}
	\caption{Nine 3-cliques in $K_{18}$.}
	\label{fig:k18alpha9case}
\end{figure}

\begin{proof}
	According to Lemma \ref{lem:k18alphage9}, it suffices to show that $D(18, \{3,4\})\neq 30$, which is equivalent to $\alpha\neq 9$. 
	
	For contradiction, assume there is a $\{K_3,K_4\}$-decomposition of $K_{18}$ with nine 3-cliques (and thus 21 4-cliques). For any vertex $x$ of $K_{18}$, we know $\alpha_x\in \{1,4,7\}$ from Lemma \ref{lem:alphax-k18}, just as in the proof of Lemma \ref{lem:k18alphage9}. 
	
	If there is some vertex $u$ with $\alpha_u=7$, by counting the total number of 3-cliques as in \eqref{eq:k18count3-cliques}, we deduce that there must be some vertex $u'\neq u$ with $\alpha_{u'}=4$, and then any of the remaining vertices $x$ has $\alpha_x=1$. 
	\begin{itemize}
		\item If $u$ and $u'$ belong to a common 3-clique, then the total number of vertices satisfies $v\ge 1 + 7\cdot 2 + (4-1)\cdot 2 = 21$, a contradiction.
		\item If $u$ and $u'$ do not share a common 3-clique, then the total number of vertices satisfies $v\ge 1+7\cdot 2 + 1+4\cdot 2 = 24$, a contradiction.
	\end{itemize}
	Therefore, for any vertex $x$ in $K_{18}$, we have $\alpha_x\in \{1,4\}$. And from the total number of 3-cliques it is easy to deduce that there must be 3 vertices $u_1$, $u_2$, $u_3$ with $\alpha_{u_1}=\alpha_{u_2}=\alpha_{u_3}=4$, and any of the other vertices $x$ has $\alpha_x=1$. 
	\begin{itemize}
		\item If $\{u_1,u_2,u_3\}$ is a 3-clique in the decomposition, then the total number of vertices satisfies $v\ge 3 + 3\cdot (4-1)\cdot 2 = 21$, a contradiction.
		\item If $\{u_1,u_2,u_3\}$ is not a 3-clique in the decomposition, then the total number of vertices satisfies $v\ge 3 + 3 + 6\cdot 2 = 18$, and the equality holds if and only if any two of the three vertices ${u_1,u_2,u_3}$ share a 3-clique, as illustrated in Figure \ref{fig:k18alpha9case}. After removing the nine 3-cliques shown in Figure \ref{fig:k18alpha9case}, an exhaustive computer-assisted search shows that the remaining graph cannot be decomposed into 4-cliques. 
	\end{itemize}
	Therefore, $\alpha\neq 9$. 
\end{proof}

\begin{lemma}\label{lem:k18alphax1or4}
	In a $\{K_3,K_4\}$-decomposition of $K_{18}$, if $\alpha=11$, then for any element $x$, $\alpha_x\in \{1,4\}$.
\end{lemma}
\begin{proof}
	Suppose that a $\{K_3,K_4\}$-decomposition of $K_{18}$ exists with $\alpha=11$. For any vertex $x$ of $K_{18}$, we know $\alpha_x\in \{1,4,7\}$ from Lemma \ref{lem:alphax-k18}.
	
	If there is some vertex $u$ with $\alpha_u=7$, by counting the total number of 3-cliques as in \eqref{eq:k18count3-cliques}, we deduce that there are two cases.
	\begin{itemize}
		\item Two vertices $u,u'$ with $\alpha_u=\alpha_{u'}=7$, another vertex $w$ with $\alpha_w=4$ and the remaining 15 vertices are in 1 3-clique each. Then both $u$ and $u'$ have at least 14 neighbours that are covered by 3-cliques with $u$ and $u'$, respectively. So at least $2\cdot 14-18=10$ vertices among them overlap. However, there are only three vertices $u$, $u'$ and $w$ that are contained in more than one 3-clique. This is a contradiction. 
		\item One vertex $u$ with $\alpha_u=7$, another three vertices $w_1,w_2,w_3$ with $\alpha_{w_i}=4$ for $i\in\{1,2,3\}$, and the remaining 14 vertices are in 1 3-clique each. Then $u$ has 14 neighbours that are covered by 3-cliques, and we denote by $x,y,z$ the three vertices in $K_{18}$ that are not contained in a common 3-clique with $u$. 
		Suppose $\{w_1,w_2,w_3\}=\{x,y,z\}$, then the remaining $11-7=4$ 3-cliques need to be formed in the $K_3$-subgraph on the vertices $w_1,w_2,w_3$, which is impossible due to insufficient edges. Therefore, at least one of $w_1,w_2,w_3$ must be among the 14 neighbours of $u$ that share a 3-clique with $u$. Without loss of generality, we assume $w_1$ and $u$ are contained commonly in a 3-clique. Then there are three more 3-cliques containing $w_1$, which means $w_1$ has 6 more neighbours that share a 3-clique with $w_1$. However, there are at most 5 candidate vertices to form a new 3-clique, namely $w_2$, $w_3$, $x$, $y$ and $z$. This is a contradiction.
	\end{itemize}
	Therefore, for any vertex $x$, $\alpha_x\in \{1,4\}$. 
\end{proof}

\begin{lemma}\label{lem:k18w3}
	In a $\{K_3,K_4\}$-decomposition of $K_{18}$, if $\alpha=11$, then there must be a 3-clique $\{w_1,w_2,w_3\}$ in the decomposition satisfying $\alpha_{w_1}=\alpha_{w_2}=\alpha_{w_3}=4$. 
\end{lemma}
\begin{proof}
	Suppose that a $\{K_3,K_4\}$-decomposition of $K_{18}$ exists with $\alpha=11$. For any vertex $x$ of $K_{18}$, we know from Lemma \ref{lem:k18alphax1or4} that $\alpha_x\in\{1,4\}$. By counting the total number of 3-cliques as in \eqref{eq:k18count3-cliques}, we deduce that there must be 5 vertices, denoted by $w_1$, $w_2$, $w_3$, $w_4$ and $w_5$, satisfying $\alpha_{w_i}=4$ for $i\in\{1,2,3,4,5\}$, and the remaining 13 vertices are only in 1 3-clique each. 
	
	Assume that no 3 of the 5 vertices $w_1,w_2,w_3,w_4,w_5$ share a 3-clique, and suppose that $t$ pairs between these 5 vertices share a 3-clique. We clearly have $0\le t\le \binom{5}{2}=10$. Since each 3-clique that contains 2 vertices among $w_1,w_2,w_3,w_4,w_5$ requires an additional vertex, and each 3-clique containing only 1 vertex among $w_1,w_2,w_3,w_4,w_5$ requires 2 additional vertices, we have
	\begin{align*}
		5+t+2\cdot(10-t)\cdot 2\le 18,
	\end{align*}
	whence $t\ge 9$. Thus, $t\in \{9,10\}$. So there are two cases:
	\begin{itemize}
		\item each pair between the vertices $w_1,w_2,w_3,w_4,w_5$ shares a 3-clique, and there is an independent 3-clique $\{x,y,z\}$ with $\alpha_x=\alpha_y=\alpha_z=1$;
		\item each pair between the vertices $w_1,w_2,w_3,w_4,w_5$, except one, shares a 3-clique, and without loss of generality we assume $w_1$ and $w_2$ do not share a 3-clique, so there are two more 3-cliques $\{w_1,x_1,y_1\}$ and $\{w_2,x_2,y_2\}$, satisfying $\alpha_{x_i}=\alpha_{y_i}=1$ for $i\in\{1,2\}$. 
	\end{itemize}
	Both cases above are excluded by computer-assisted exhaustive search. Therefore, there must be 3 of the 5 vertices $w_1,w_2,w_3,w_4,w_5$ that share a 3-clique, which we assume, without loss of generality, to be $\{w_1,w_2,w_3\}$. 
\end{proof}

\begin{lemma}\label{lem:k18alphage13}
	$D(18,\{3,4\})\ge 32$, or equivalently, $\alpha\ge 13$.
\end{lemma}
\begin{proof}
	By Lemma \ref{lem:k18alphage11}, it suffices to show that $D(18,\{3,4\})\neq 31$, which is equivalent to $\alpha\neq 11$. 
	
	For contradiction, assume there is a $\{K_3,K_4\}$-decomposition of $K_{18}$ with eleven 3-cliques (and thus 20 4-cliques). By the proof of Lemma \ref{lem:k18w3}, there must be 5 vertices $w_1,w_2,w_3,w_4,w_5$ satisfying $\alpha_{w_i}=4$ for $i\in\{1,2,3,4,5\}$ and the remaining 13 vertices are only in 1 3-clique each, and moreover there must be a 3-clique among $w_1,w_2,w_3,w_4,w_5$. Without loss of generality, we assume that there is 3-clique $T=\{w_1,w_2,w_3\}$. Denote by $N_i$ the subset of $V(K_{18})\setminus T$ consisting of the vertices that share a 3-clique with $w_i$, for $i\in \{1,2,3\}$. Since each of $w_4$ and $w_5$ is in at least one 3-clique that does not contain any vertex of $T$, and at least one neighbour of $w_4$ and $w_5$ is not in $N_1\cup N_2\cup N_3$, we have
	\begin{align*}
		|N_1\cup N_2\cup N_3|\le 18-3-1=14.
	\end{align*}
	Meanwhile, the inclusion-exclusion principle gives
	\begin{align*}
		|N_1\cup N_2\cup N_3| = 3\cdot 6 - \sum_{1\le i<j\le 3}|N_i\cap N_j| + |N_1\cap N_2\cap N_3|,
	\end{align*}
	since $|N_1|=|N_2|=|N_3|=6$. We thus deduce
	\begin{align}\label{eq:k18n3}
		|N_1\cap N_2|+|N_1\cap N_3|+|N_2\cap N_3|-|N_1\cap N_2\cap N_3|\ge 4.
	\end{align}
	Notice that for any pair $i,j$ satisfying $1\le i<j\le 3$, $N_i\cap N_j\subseteq \{w_4,w_5\}$ holds. We claim that 
	\begin{align*}
		N_1\cap N_2\cap N_3 = \{w_4,w_5\}.
	\end{align*}

	We prove the claim in four steps.
	\begin{enumerate}
		\item First suppose there is some $N_i$ that does not contain either $w_4$ or $w_5$, where $i\in\{1,2,3\}$. Without loss of generality, we may assume $i=1$, and then $N_1\cap N_2=N_1\cap N_3=N_1\cap N_2\cap N_3=\emptyset$. Then \eqref{eq:k18n3} gives $|N_2\cap N_3|\ge 4$, which is impossible. So each $N_i$ contains either $w_4$ or $w_5$.
		\item\label{item:step2}%
                  Next, suppose there is some pair $i,j$ with $1\le i<j\le 3$ such that $N_i\cap N_j=\emptyset$. Without loss of generality, let $(i,j)=(1,2)$, and from the previous step we may assume $w_4\in N_1\setminus N_2$ and $w_5\in N_2\setminus N_1$. So $N_1\cap N_2\cap N_3=\emptyset$, and $|N_1\cap N_3|\le 1$, $|N_2\cap N_3|\le 1$. Then \eqref{eq:k18n3} gives
		\begin{align*}
			4\le |N_1\cap N_3|+|N_2\cap N_3|\le 1+1=2,
		\end{align*}
		a contradiction. So each pair among $N_i$'s intersects. 
		\item Then suppose $N_1\cap N_2\cap N_3=\emptyset$. Since 
		\begin{align*}
			N_1\cap N_2\cap N_3=(N_1\cap N_2)\cap (N_1\cap N_3),
		\end{align*}
		we may assume, without loss of generality, that $N_1\cap N_2=\{w_4\}$ and $N_1\cap N_3=\{w_5\}$, for they are both non-empty. So $w_4,w_5\in N_1$, $w_4\in N_2\setminus N_3$, and $w_5\in N_3\setminus N_2$. Therefore, $N_2\cap N_3=\emptyset$, which is a contradiction as has been shown in the step (\ref{item:step2}.)! We thus get $N_1\cap N_2\cap N_3\neq \emptyset$. 
		\item Finally, suppose there is some pair $i,j$ with $1\le i<j\le 3$ such that $|N_i\cap N_j|=1$. Without loss of generality, we assume $N_1\cap N_2 = \{w_4\}$ and $w_5\notin N_1$. So $N_1\cap N_3 = \{w_4\}$, and $N_1\cap N_2\cap N_3=\{w_4\}$. Therefore, the inequality \eqref{eq:k18n3} gives $|N_2\cap N_3|\ge 3$, which is impossible. So $N_1\cap N_2=N_1\cap N_3=N_2\cap N_3=\{w_4,w_5\}$, and thus $N_1\cap N_2\cap N_3=\{w_4,w_5\}$. 
	\end{enumerate}
\begin{figure}[tbp]
	\centering
	\definecolor{ccqqww}{rgb}{0.8,0.,0.4}
	\definecolor{zzttqq}{rgb}{0.6,0.2,0.}
	\definecolor{qqzzcc}{rgb}{0.,0.6,0.8}
	\definecolor{zzccqq}{rgb}{0.6,0.8,0.}
	\definecolor{aqaqaq}{rgb}{0.6274509803921569,0.6274509803921569,0.6274509803921569}
	\definecolor{ududff}{rgb}{0.30196078431372547,0.30196078431372547,1.}
	\begin{tikzpicture}[line cap=round,line join=round,>=triangle 45,x=0.8cm,y=0.8cm]
		\clip(2.599815564917047,-5.562614285932006) rectangle (14.916636718690448,3.768916268195108);
		\fill[line width=2.pt,color=aqaqaq,fill=aqaqaq,fill opacity=0.10000000149011612] (8.4,1.9) -- (5.132,-3.6504) -- (11.7628,-3.7956) -- cycle;
		\fill[line width=2.pt,color=zzccqq,fill=zzccqq,fill opacity=0.10000000149011612] (8.4,1.9) -- (7.2616,-2.4162) -- (6.342,0.6088) -- cycle;
		\fill[line width=2.pt,color=zzccqq,fill=zzccqq,fill opacity=0.10000000149011612] (5.132,-3.6504) -- (7.2616,-2.4162) -- (5.4708,-1.8112) -- cycle;
		\fill[line width=2.pt,color=zzccqq,fill=zzccqq,fill opacity=0.10000000149011612] (11.7628,-3.7956) -- (7.2616,-2.4162) -- (9.004,-4.4732) -- cycle;
		\fill[line width=2.pt,color=qqzzcc,fill=qqzzcc,fill opacity=0.10000000149011612] (8.4,1.9) -- (9.5848,-2.4404) -- (10.4076,0.4636) -- cycle;
		\fill[line width=2.pt,color=qqzzcc,fill=qqzzcc,fill opacity=0.10000000149011612] (5.132,-3.6504) -- (9.5848,-2.4404) -- (7.068,-4.4248) -- cycle;
		\fill[line width=2.pt,color=qqzzcc,fill=qqzzcc,fill opacity=0.10000000149011612] (11.7628,-3.7956) -- (9.5848,-2.4404) -- (11.424,-1.908) -- cycle;
		\fill[line width=2.pt,color=zzttqq,fill=zzttqq,fill opacity=0.10000000149011612] (8.4,1.9) -- (7.431,3.0772) -- (9.488,3.1014) -- cycle;
		\fill[line width=2.pt,color=zzttqq,fill=zzttqq,fill opacity=0.10000000149011612] (5.132,-3.6504) -- (3.4138,-3.1664) -- (4.3334,-5.054) -- cycle;
		\fill[line width=2.pt,color=zzttqq,fill=zzttqq,fill opacity=0.10000000149011612] (11.7628,-3.7956) -- (12.2226,-5.2476) -- (13.3842,-3.3358) -- cycle;
		\fill[line width=2.pt,color=ccqqww,fill=ccqqww,fill opacity=0.10000000149011612] (7.2616,-2.4162) -- (9.5848,-2.4404) -- (8.4958,-0.8432) -- cycle;
		\draw [line width=2.pt,color=aqaqaq] (8.4,1.9)-- (5.132,-3.6504);
		\draw [line width=2.pt,color=aqaqaq] (5.132,-3.6504)-- (11.7628,-3.7956);
		\draw [line width=2.pt,color=aqaqaq] (11.7628,-3.7956)-- (8.4,1.9);
		\draw [line width=2.pt,color=zzccqq] (8.4,1.9)-- (7.2616,-2.4162);
		\draw [line width=2.pt,color=zzccqq] (7.2616,-2.4162)-- (6.342,0.6088);
		\draw [line width=2.pt,color=zzccqq] (6.342,0.6088)-- (8.4,1.9);
		\draw [line width=2.pt,color=zzccqq] (5.132,-3.6504)-- (7.2616,-2.4162);
		\draw [line width=2.pt,color=zzccqq] (7.2616,-2.4162)-- (5.4708,-1.8112);
		\draw [line width=2.pt,color=zzccqq] (5.4708,-1.8112)-- (5.132,-3.6504);
		\draw [line width=2.pt,color=zzccqq] (11.7628,-3.7956)-- (7.2616,-2.4162);
		\draw [line width=2.pt,color=zzccqq] (7.2616,-2.4162)-- (9.004,-4.4732);
		\draw [line width=2.pt,color=zzccqq] (9.004,-4.4732)-- (11.7628,-3.7956);
		\draw [line width=2.pt,color=qqzzcc] (8.4,1.9)-- (9.5848,-2.4404);
		\draw [line width=2.pt,color=qqzzcc] (9.5848,-2.4404)-- (10.4076,0.4636);
		\draw [line width=2.pt,color=qqzzcc] (10.4076,0.4636)-- (8.4,1.9);
		\draw [line width=2.pt,color=qqzzcc] (5.132,-3.6504)-- (9.5848,-2.4404);
		\draw [line width=2.pt,color=qqzzcc] (9.5848,-2.4404)-- (7.068,-4.4248);
		\draw [line width=2.pt,color=qqzzcc] (7.068,-4.4248)-- (5.132,-3.6504);
		\draw [line width=2.pt,color=qqzzcc] (11.7628,-3.7956)-- (9.5848,-2.4404);
		\draw [line width=2.pt,color=qqzzcc] (9.5848,-2.4404)-- (11.424,-1.908);
		\draw [line width=2.pt,color=qqzzcc] (11.424,-1.908)-- (11.7628,-3.7956);
		\draw [line width=2.pt,color=zzttqq] (8.4,1.9)-- (7.431,3.0772);
		\draw [line width=2.pt,color=zzttqq] (7.431,3.0772)-- (9.488,3.1014);
		\draw [line width=2.pt,color=zzttqq] (9.488,3.1014)-- (8.4,1.9);
		\draw [line width=2.pt,color=zzttqq] (5.132,-3.6504)-- (3.4138,-3.1664);
		\draw [line width=2.pt,color=zzttqq] (3.4138,-3.1664)-- (4.3334,-5.054);
		\draw [line width=2.pt,color=zzttqq] (4.3334,-5.054)-- (5.132,-3.6504);
		\draw [line width=2.pt,color=zzttqq] (11.7628,-3.7956)-- (12.2226,-5.2476);
		\draw [line width=2.pt,color=zzttqq] (12.2226,-5.2476)-- (13.3842,-3.3358);
		\draw [line width=2.pt,color=zzttqq] (13.3842,-3.3358)-- (11.7628,-3.7956);
		\draw [line width=2.pt,color=ccqqww] (7.2616,-2.4162)-- (9.5848,-2.4404);
		\draw [line width=2.pt,color=ccqqww] (9.5848,-2.4404)-- (8.4958,-0.8432);
		\draw [line width=2.pt,color=ccqqww] (8.4958,-0.8432)-- (7.2616,-2.4162);
		\begin{scriptsize}
			\draw [color=ududff] (8.4,1.9) circle (2.5pt);
			\draw[color=ududff] (8.881633916077396,2.003718174536466) node {$w_1$};
			\draw [color=ududff] (5.132,-3.6504) circle (2.5pt);
			\draw[color=ududff] (4.907685121765494,-3.2394158185916237) node {$w_2$};
			\draw [color=ududff] (11.7628,-3.7956) circle (2.5pt);
			\draw[color=ududff] (12.045898539197904,-3.3861790404056055) node {$w_3$};
			\draw [color=ududff] (7.2616,-2.4162) circle (2.5pt);
			\draw[color=ududff] (7.241396577368768,-2.80886265107439) node {$w_4$};
			\draw [fill=ududff] (6.342,0.6088) circle (1.5pt);
			\draw [fill=ududff] (5.4708,-1.8112) circle (1.5pt);
			\draw [fill=ududff] (9.004,-4.4732) circle (1.5pt);
			\draw [color=ududff] (9.5848,-2.4404) circle (2.5pt);
			\draw[color=ududff] (9.167029169190327,-2.231441608276541) node {$w_5$};
			\draw [fill=ududff] (10.4076,0.4636) circle (1.5pt);
			\draw [fill=ududff] (7.068,-4.4248) circle (1.5pt);
			\draw [fill=ududff] (11.424,-1.908) circle (1.5pt);
			\draw [fill=ududff] (7.431,3.0772) circle (1.5pt);
			\draw [fill=ududff] (9.488,3.1014) circle (1.5pt);
			\draw [fill=ududff] (3.4138,-3.1664) circle (1.5pt);
			\draw [fill=ududff] (4.3334,-5.054) circle (1.5pt);
			\draw [fill=ududff] (12.2226,-5.2476) circle (1.5pt);
			\draw [fill=ududff] (13.3842,-3.3358) circle (1.5pt);
			\draw [fill=ududff] (8.4958,-0.8432) circle (1.5pt);
		\end{scriptsize}
	\end{tikzpicture}
	\caption{11 copies of $K_3$ in $K_{18}$}
	\label{fig:k18alpha11case}
\end{figure}
	Therefore, for $i\in \{1,2,3\}$ and $j\in\{4,5\}$, there is a 3-clique containing $w_i$ and $w_j$. Since $\alpha_{w_4}=\alpha_{w_5}=4$, each of $w_4$ and $w_5$ must be in an additional 3-clique. By counting the total number of vertices in the 3-cliques, we can easily conclude that the last 3-clique must be of the form $\{w_4,w_5,x\}$, where $\alpha_x=1$. This unique possible configuration (See Figure \ref{fig:k18alpha11case}) is then excluded by an exhaustive computer-assisted search. Hence $\alpha\neq 11$. 
\end{proof}

\begin{lemma}\label{lem:k18alphax1or4-13}
	In a $\{K_3,K_4\}$-decomposition of $K_{18}$, if $\alpha=13$, then for any element $x$, $\alpha_x\in \{1,4\}$.
\end{lemma}
\begin{proof}
	Suppose that a $\{K_3,K_4\}$-decomposition of $K_{18}$ exists with $\alpha=13$. For any vertex $x$ of $K_{18}$, from Lemma \ref{lem:alphax-k18} we know the number $\alpha_x$ of 3-cliques containing $x$ satisfies $\alpha_x\in \{1,4,7\}$.
	
	If there is some vertex $u$ with $\alpha_u=7$, by counting the total number of 3-cliques as in \eqref{eq:k18count3-cliques}, we deduce that there are three cases.
        Each of them we show to be impossible.
	\begin{itemize}
		\item Three vertices $u_1,u_2,u_3$ with $\alpha_{u_i}=7$ for $i\in\{1,2,3\}$, another vertex $w$ with $\alpha_w=4$ and the remaining 14 vertices only in 1 3-clique each. In particular, both $u_1$ and $u_2$ have 14 neighbours covered by 3-cliques with $u_1$ and $u_2$, respectively, so at least $2\cdot 14-18=10$ vertices among them overlap. But there are only four vertices $u_1,u_2,u_3$ and $w$ that are contained in more than one 3-clique. This is a contradiction.
		\item Two vertices $u_1, u_2$ with $\alpha_{u_1}=\alpha_{u_2}=7$, another three vertices $w_1,w_2,w_3$ with $\alpha_{w_i}=4$ for $i\in\{1,2,3\}$, and the remaining 13 vertices are only in 1 3-clique each. Then both $u_1$ and $u_2$ have 14 neighbours covered by 3-cliques with $u_1$ and $u_2$, respectively, so at least $2\cdot 14-18=10$ vertices among them overlap. But there are only five vertices $u_1,u_2,w_1,w_2,w_3$ that are contained in more than one 3-clique. This is a contradiction.
		\item One vertex $u$ with $\alpha_u=7$, another five vertices $w_1,w_2,w_3,w_4,w_5$ with $\alpha_{w_i}=4$ for $i\in\{1,2,3,4,5\}$ and the remaining 12 vertices are only in 1 3-clique each. Then $\beta_u=(17-2 \alpha_u)/3=1$, and we assume the unique 4-clique containing $u$ is $\{u,x,y,z\}$. If $w_i\in\{x,y,z\}$ for some $i\in\{1,2,3,4,5\}$, without loss of generality we may assume $w_1=x$, then to form 4 3-cliques containing $w_1$, 8 vertices are needed, but only 6 are available: $w_2,w_3,w_4,w_5,y,z$. Therefore, each $w_i$ shares a 3-clique with $u$ for $i\in\{1,2,3,4,5\}$. It is now necessary to form the remaining $13-7=6$ 3-cliques in the $K_8$-subgraph on $w_1,w_2,w_3,w_4,w_5,x,y,z$. The union of the 6 3-cliques should be a subgraph of $K_8$, with degree sequence $(6,6,6,6,6,2,2,2)$. The sequence is not graphic, since summing the first five numbers in the sequence, we get 
		\begin{align*}
			6+6+6+6+6 > 5(5-1)+2+2+2,
		\end{align*}
		which contradicts Theorem \ref{thm:erdos-gallai}.
	\end{itemize}
	Therefore, $\alpha_x\in\{1,4\}$ for each vertex $x$. 
\end{proof}

\begin{lemma}\label{lem:k18k3k4alpha13cases}
	In a $\{K_3,K_4\}$-decomposition of $K_{18}$, if $\alpha=13$, then there are 7 vertices $w_1,w_2,\cdots,w_7$ such that $\alpha_{w_i}=4$ for $i\in\{1,2,\cdots,7\}$, and there are at least two 3-cliques among the 7 vertices in the decomposition. 
\end{lemma}
\begin{proof}
	Suppose that a $\{K_3,K_4\}$-decomposition of $K_{18}$ exists with $\alpha=13$. By Lemma \ref{lem:k18alphax1or4-13}, each vertex is in one or four 3-cliques. By counting the total number of 3-cliques as in \eqref{eq:k18count3-cliques}, we deduce that there are 7 vertices $w_1,w_2,\cdots,w_7$ each in 4 3-cliques, and the remaining 11 vertices $y_1,y_2,\cdots,y_{11}$ are only in 1 3-clique each. 
	
	Denote by $t_i$ the number of 3-cliques in the decomposition with $i$ vertices among $w_1,w_2,\cdots,w_7$ and $3-i$ vertices among $y_1,y_2,\cdots,y_{11}$, for $i\in\{0,1,2,3\}$. Clearly, the integers $t_i$ satisfy:
	\begin{align*}
		t_0+t_1+t_2+t_3&=13,\\
		3t_3+2t_2+t_1&=\sum_{i=1}^7 \alpha_{w_i} = 28,\\
		3t_0+2t_1+t_2&=\sum_{j=1}^{11} \alpha_{y_j} = 11,\\
		t_0&\le \left\lfloor \frac{11}{3}\right\rfloor = 3,\\
		t_1&\le \left\lfloor \frac{11}{2}\right\rfloor = 5,\\
		t_2&\le 11,\\
		t_3&\le 7,\\
		t_2+3t_3&\le e(K_7) = 21,\\
		t_0,t_1,t_2,t_3&\ge 0.
	\end{align*}
	Solving the constraints above, we get 7 possible solutions:
	\begin{align*}
		(t_0,t_1,t_2,t_3)&=(0,4,3,6)\label{eq:c1} \tag{C1},\\
		(t_0,t_1,t_2,t_3)&=(1,1,6,5)\label{eq:c2} \tag{C2},\\
		(t_0,t_1,t_2,t_3)&=(0,3,5,5)\label{eq:c3} \tag{C3},\\
		(t_0,t_1,t_2,t_3)&=(1,0,8,4)\label{eq:c4} \tag{C4},\\
		(t_0,t_1,t_2,t_3)&=(0,2,7,4)\label{eq:c5} \tag{C5},\\
		(t_0,t_1,t_2,t_3)&=(0,1,9,3)\label{eq:c6} \tag{C6},\\
		(t_0,t_1,t_2,t_3)&=(0,0,11,2)\label{eq:c7} \tag{C7}.
	\end{align*}
	In particular, $t_3\ge 2$ as desired. 
\end{proof}

From now on, we consider a decomposition of $K_{18}$ into 13 3-cliques and 19 4-cliques, where seven vertices $w_1,\ldots,w_7$ each lie in four 3-cliques and the remaining eleven vertices $y_1,\ldots,y_{11}$ each lie in exactly one 3-clique. We shall deal with the cases \eqref{eq:c1}--\eqref{eq:c7} in the proof of Lemma \ref{lem:k18k3k4alpha13cases} one by one and show none is feasible. 

\begin{lemma}\label{lem:k18a13c1}
	In a $\{K_3,K_4\}$-decomposition of $K_{18}$, if $\alpha=13$, then the case \eqref{eq:c1} is impossible.
\end{lemma}
\begin{figure}[tbp]
	\centering
	\definecolor{zzttqq}{rgb}{0.6,0.2,0.}
	\definecolor{qqwuqq}{rgb}{0.,0.39215686274509803,0.}
	\definecolor{qqqqcc}{rgb}{0.,0.,0.8}
	\definecolor{ududff}{rgb}{0.30196078431372547,0.30196078431372547,1.}
	\begin{tikzpicture}[line cap=round,line join=round,>=triangle 45,x=0.5cm,y=0.5cm]
		\clip(1.8431311302969944,-8.727801435217188) rectangle (17.422408733306032,4.902186579189873);
		\fill[line width=2.pt,color=qqqqcc,fill=qqqqcc,fill opacity=0.10000000149011612] (9.76,3.32) -- (14.,1.38) -- (12.7416,3.5062) -- cycle;
		\fill[line width=2.pt,color=qqqqcc,fill=qqqqcc,fill opacity=0.10000000149011612] (9.76,3.32) -- (6.12,1.64) -- (7.1272,3.482) -- cycle;
		\fill[line width=2.pt,color=qqqqcc,fill=qqqqcc,fill opacity=0.10000000149011612] (6.12,1.64) -- (14.,1.38) -- (15.2584,-0.39) -- cycle;
		\fill[line width=2.pt,color=qqwuqq,fill=qqwuqq,fill opacity=0.10000000149011612] (14.3,-2.84) -- (16.2506,-2.1566) -- (15.7424,-4.4072) -- cycle;
		\fill[line width=2.pt,color=qqwuqq,fill=qqwuqq,fill opacity=0.10000000149011612] (11.8,-6.14) -- (14.242,-7.2386) -- (10.7088,-7.9646) -- cycle;
		\fill[line width=2.pt,color=qqwuqq,fill=qqwuqq,fill opacity=0.10000000149011612] (7.42,-6.04) -- (8.1436,-7.7468) -- (5.5542,-6.7062) -- cycle;
		\fill[line width=2.pt,color=qqwuqq,fill=qqwuqq,fill opacity=0.10000000149011612] (4.84,-2.5) -- (4.2474,-4.625) -- (3.11,-0.995) -- cycle;
		\fill[line width=2.pt,color=zzttqq,fill=zzttqq,fill opacity=0.10000000149011612] (9.76,3.32) -- (14.3,-2.84) -- (11.8,-6.14) -- cycle;
		\fill[line width=2.pt,color=zzttqq,fill=zzttqq,fill opacity=0.10000000149011612] (9.76,3.32) -- (7.42,-6.04) -- (4.84,-2.5) -- cycle;
		\fill[line width=2.pt,color=zzttqq,fill=zzttqq,fill opacity=0.10000000149011612] (14.,1.38) -- (14.3,-2.84) -- (7.42,-6.04) -- cycle;
		\fill[line width=2.pt,color=zzttqq,fill=zzttqq,fill opacity=0.10000000149011612] (14.,1.38) -- (11.8,-6.14) -- (4.84,-2.5) -- cycle;
		\fill[line width=2.pt,color=zzttqq,fill=zzttqq,fill opacity=0.10000000149011612] (6.12,1.64) -- (4.84,-2.5) -- (14.3,-2.84) -- cycle;
		\fill[line width=2.pt,color=zzttqq,fill=zzttqq,fill opacity=0.10000000149011612] (6.12,1.64) -- (7.42,-6.04) -- (11.8,-6.14) -- cycle;
		\draw [line width=2.pt,color=qqqqcc] (9.76,3.32)-- (14.,1.38);
		\draw [line width=2.pt,color=qqqqcc] (14.,1.38)-- (12.7416,3.5062);
		\draw [line width=2.pt,color=qqqqcc] (12.7416,3.5062)-- (9.76,3.32);
		\draw [line width=2.pt,color=qqqqcc] (9.76,3.32)-- (6.12,1.64);
		\draw [line width=2.pt,color=qqqqcc] (6.12,1.64)-- (7.1272,3.482);
		\draw [line width=2.pt,color=qqqqcc] (7.1272,3.482)-- (9.76,3.32);
		\draw [line width=2.pt,color=qqqqcc] (6.12,1.64)-- (14.,1.38);
		\draw [line width=2.pt,color=qqqqcc] (14.,1.38)-- (15.2584,-0.39);
		\draw [line width=2.pt,color=qqqqcc] (15.2584,-0.39)-- (6.12,1.64);
		\draw [line width=2.pt,color=qqwuqq] (14.3,-2.84)-- (16.2506,-2.1566);
		\draw [line width=2.pt,color=qqwuqq] (16.2506,-2.1566)-- (15.7424,-4.4072);
		\draw [line width=2.pt,color=qqwuqq] (15.7424,-4.4072)-- (14.3,-2.84);
		\draw [line width=2.pt,color=qqwuqq] (11.8,-6.14)-- (14.242,-7.2386);
		\draw [line width=2.pt,color=qqwuqq] (14.242,-7.2386)-- (10.7088,-7.9646);
		\draw [line width=2.pt,color=qqwuqq] (10.7088,-7.9646)-- (11.8,-6.14);
		\draw [line width=2.pt,color=qqwuqq] (7.42,-6.04)-- (8.1436,-7.7468);
		\draw [line width=2.pt,color=qqwuqq] (8.1436,-7.7468)-- (5.5542,-6.7062);
		\draw [line width=2.pt,color=qqwuqq] (5.5542,-6.7062)-- (7.42,-6.04);
		\draw [line width=2.pt,color=qqwuqq] (4.84,-2.5)-- (4.2474,-4.625);
		\draw [line width=2.pt,color=qqwuqq] (4.2474,-4.625)-- (3.11,-0.995);
		\draw [line width=2.pt,color=qqwuqq] (3.11,-0.995)-- (4.84,-2.5);
		\draw [line width=2.pt,color=zzttqq] (9.76,3.32)-- (14.3,-2.84);
		\draw [line width=2.pt,color=zzttqq] (14.3,-2.84)-- (11.8,-6.14);
		\draw [line width=2.pt,color=zzttqq] (11.8,-6.14)-- (9.76,3.32);
		\draw [line width=2.pt,color=zzttqq] (9.76,3.32)-- (7.42,-6.04);
		\draw [line width=2.pt,color=zzttqq] (7.42,-6.04)-- (4.84,-2.5);
		\draw [line width=2.pt,color=zzttqq] (4.84,-2.5)-- (9.76,3.32);
		\draw [line width=2.pt,color=zzttqq] (14.,1.38)-- (14.3,-2.84);
		\draw [line width=2.pt,color=zzttqq] (14.3,-2.84)-- (7.42,-6.04);
		\draw [line width=2.pt,color=zzttqq] (7.42,-6.04)-- (14.,1.38);
		\draw [line width=2.pt,color=zzttqq] (14.,1.38)-- (11.8,-6.14);
		\draw [line width=2.pt,color=zzttqq] (11.8,-6.14)-- (4.84,-2.5);
		\draw [line width=2.pt,color=zzttqq] (4.84,-2.5)-- (14.,1.38);
		\draw [line width=2.pt,color=zzttqq] (6.12,1.64)-- (4.84,-2.5);
		\draw [line width=2.pt,color=zzttqq] (4.84,-2.5)-- (14.3,-2.84);
		\draw [line width=2.pt,color=zzttqq] (14.3,-2.84)-- (6.12,1.64);
		\draw [line width=2.pt,color=zzttqq] (6.12,1.64)-- (7.42,-6.04);
		\draw [line width=2.pt,color=zzttqq] (7.42,-6.04)-- (11.8,-6.14);
		\draw [line width=2.pt,color=zzttqq] (11.8,-6.14)-- (6.12,1.64);
		\begin{scriptsize}
			\draw [fill=ududff] (9.76,3.32) circle (2.5pt);
			\draw[color=ududff] (9.984757539728044,3.8873361982539434) node {$w_1$};
			\draw [fill=ududff] (14.,1.38) circle (2.5pt);
			\draw[color=ududff] (14.208608834846407,1.9590562591781708) node {$w_2$};
			\draw [fill=ududff] (14.3,-2.84) circle (2.5pt);
			\draw[color=ududff] (14.814685015652086,-2.2647950359401876) node {$w_3$};
			\draw [fill=ududff] (11.8,-6.14) circle (2.5pt);
			\draw[color=ududff] (12.604860333045523,-5.870417788641512) node {$w_4$};
			\draw [fill=ududff] (7.42,-6.04) circle (2.5pt);
			\draw[color=ududff] (6.627970947524319,-5.878594934399808) node {$w_5$};
			\draw [fill=ududff] (4.84,-2.5) circle (2.5pt);
			\draw[color=ududff] (4.656931028756619,-1.628111237053942) node {$w_6$};
			\draw [fill=ududff] (6.12,1.64) circle (2.5pt);
			\draw[color=ududff] (5.642450988140469,1.7039172038227135) node {$w_7$};
			\draw [fill=ududff] (12.7416,3.5062) ++(-2.5pt,0 pt) -- ++(2.5pt,2.5pt)--++(2.5pt,-2.5pt)--++(-2.5pt,-2.5pt)--++(-2.5pt,2.5pt);
			\draw [fill=ududff] (7.1272,3.482) ++(-2.5pt,0 pt) -- ++(2.5pt,2.5pt)--++(2.5pt,-2.5pt)--++(-2.5pt,-2.5pt)--++(-2.5pt,2.5pt);
			\draw [fill=ududff] (15.2584,-0.39) ++(-2.5pt,0 pt) -- ++(2.5pt,2.5pt)--++(2.5pt,-2.5pt)--++(-2.5pt,-2.5pt)--++(-2.5pt,2.5pt);
			\draw [fill=ududff] (16.2506,-2.1566) ++(-2.5pt,0 pt) -- ++(2.5pt,2.5pt)--++(2.5pt,-2.5pt)--++(-2.5pt,-2.5pt)--++(-2.5pt,2.5pt);
			\draw [fill=ududff] (15.7424,-4.4072) ++(-2.5pt,0 pt) -- ++(2.5pt,2.5pt)--++(2.5pt,-2.5pt)--++(-2.5pt,-2.5pt)--++(-2.5pt,2.5pt);
			\draw [fill=ududff] (14.242,-7.2386) ++(-2.5pt,0 pt) -- ++(2.5pt,2.5pt)--++(2.5pt,-2.5pt)--++(-2.5pt,-2.5pt)--++(-2.5pt,2.5pt);
			\draw [fill=ududff] (10.7088,-7.9646) ++(-2.5pt,0 pt) -- ++(2.5pt,2.5pt)--++(2.5pt,-2.5pt)--++(-2.5pt,-2.5pt)--++(-2.5pt,2.5pt);
			\draw [fill=ududff] (8.1436,-7.7468) ++(-2.5pt,0 pt) -- ++(2.5pt,2.5pt)--++(2.5pt,-2.5pt)--++(-2.5pt,-2.5pt)--++(-2.5pt,2.5pt);
			\draw [fill=ududff] (5.5542,-6.7062) ++(-2.5pt,0 pt) -- ++(2.5pt,2.5pt)--++(2.5pt,-2.5pt)--++(-2.5pt,-2.5pt)--++(-2.5pt,2.5pt);
			\draw [fill=ududff] (4.2474,-4.625) ++(-2.5pt,0 pt) -- ++(2.5pt,2.5pt)--++(2.5pt,-2.5pt)--++(-2.5pt,-2.5pt)--++(-2.5pt,2.5pt);
			\draw [fill=ududff] (3.11,-0.995) ++(-2.5pt,0 pt) -- ++(2.5pt,2.5pt)--++(2.5pt,-2.5pt)--++(-2.5pt,-2.5pt)--++(-2.5pt,2.5pt);
		\end{scriptsize}
	\end{tikzpicture}
	\caption{13 3-cliques in $K_{18}$ for case \eqref{eq:c1}}
	\label{fig:k18alpha13-1}
\end{figure}
\begin{proof}
For contradiction, suppose case \eqref{eq:c1} is possible. Then six 3-cliques lie entirely among $w_1,\ldots,w_7$. Denote by $H$ the union of these 6 3-cliques and consider its complement $\overline{H}$ in the $K_7$-subgraph induced on vertices $w_1,\cdots,w_7$. Clearly, $\overline{H}$ is an even graph with $e(K_7)-6\cdot 3=3$ edges. Ignoring the isolated vertices in $\overline{H}$, we can easily deduce $\overline{H}\cong C_3$. The three edges in $\overline{H}$ provide the $t_2=3$ 3-cliques, and the four isolated vertices in $\overline{H}$ provide the $t_1=4$ 3-cliques, hence the distribution of 3-cliques in such a decomposition is unique, see Figure \ref{fig:k18alpha13-1}. This configuration is excluded by an exhaustive computer-assisted search on the remaining edges. 
\end{proof}

\begin{lemma}\label{lem:k18alpha13c2}
	In a $\{K_3,K_4\}$-decomposition of $K_{18}$, if $\alpha=13$, then the case \eqref{eq:c2} is impossible.
\end{lemma}
\begin{proof}
For contradiction, suppose case \eqref{eq:c2} is possible. Then five 3-cliques lie entirely among $w_1,\ldots,w_7$, so the complement of their union in the induced $K_7$ is an even graph with six edges. An isomorphism-free computer enumeration gives exactly three even graphs $H$ of order 7 and size 6, shown in Figure \ref{fig:nauty-7-6-e}, where vertex $i$ represents $w_{i+1}$ for $i\in\{0,1,\ldots,6\}$. The third graph, $C_3+C_3$, is excluded because its complement is not decomposable into 3-cliques; this is confirmed by an exhaustive computer-assisted search. 
	
	If $H$ is a bowtie as the first graph in Figure \ref{fig:nauty-7-6-e}, then the vertices $w_3$ and $w_4$ are each contained in three 3-cliques in $\overline{H}$, so they each require an additional 3-clique with two vertices among $y_1,\cdots,y_{11}$. This contradicts $t_1=1$ in \eqref{eq:c2}. 
	
	If $H\cong C_6$ is a cycle as the second graph in Figure \ref{fig:nauty-7-6-e}, we take the 5 3-cliques $\{w_1,w_2,w_3\}$, $\{w_1,w_4,w_7\}$, $\{w_2,w_4,w_6\}$, $\{w_3,w_4,w_5\}$, $\{w_5,w_6,w_7\}$ from $\overline{H}$, and we also notice that $w_4$ requires an additional 3-clique, which we can take, without loss of generality, to be $\{w_4,y_1,y_2\}$. For each edge in $H$, we form a 3-clique with a vertex among $y_1,\cdots,y_{11}$: $\{w_1,w_5,y_3\}$, $\{w_1,w_6,y_4\}$, $\{w_2,w_5,y_5\}$, $\{w_2,w_7,y_6\}$, $\{w_3,w_6,y_7\}$, $\{w_3,w_7,y_8\}$. Finally, we append a 3-clique $\{y_9,y_{10},y_{11}\}$. This configuration of 13 3-cliques (See Figure \ref{fig:k18alpha13-5}) is shown to be infeasible by an exhaustive computer-assisted search. 
	
	To sum up, the case \eqref{eq:c2} is impossible.
\end{proof}

\begin{figure}[tbp]
	\centering
	\definecolor{uuuuuu}{rgb}{0.26666666666666666,0.26666666666666666,0.26666666666666666}
	\begin{tikzpicture}[line cap=round,line join=round,>=triangle 45,x=0.6cm,y=0.6cm]
		\clip(-1.4252977007283407,1.5) rectangle (17.090576152895668,6.5);
		\draw [line width=2.pt] (1.99,5.613590426001814)-- (1.22,2.24);
		\draw [line width=2.pt] (0.6025079434302745,4.945409467760775)-- (1.99,5.613590426001814);
		\draw [line width=2.pt] (0.6025079434302745,4.945409467760775)-- (0.2598257051375503,3.444020483000766);
		\draw [line width=2.pt] (0.6025079434302745,4.945409467760775)-- (3.377492056569725,4.945409467760774);
		\draw [line width=2.pt] (0.2598257051375503,3.444020483000766)-- (3.377492056569725,4.945409467760774);
		\draw [line width=2.pt] (0.6025079434302745,4.945409467760775)-- (1.22,2.24);
		\draw [line width=2.pt] (7.793812862675348,5.455964700651117)-- (7.,2.18);
		\draw [line width=2.pt] (6.433681886039368,4.823158469332829)-- (9.136588489746625,4.787119714616732);
		\draw [line width=2.pt] (6.433681886039368,4.823158469332829)-- (9.45087133243746,3.32027742766487);
		\draw [line width=2.pt] (9.45087133243746,3.32027742766487)-- (6.08040192686126,3.3652170197392195);
		\draw [line width=2.pt] (6.08040192686126,3.3652170197392195)-- (7.793812862675348,5.455964700651117);
		\draw [line width=2.pt] (7.,2.18)-- (9.136588489746625,4.787119714616732);
		\draw [line width=2.pt] (13.68,5.477403288677162)-- (14.46,2.06);
		\draw [line width=2.pt] (14.46,2.06)-- (11.927355909100376,3.2796571126501264);
		\draw [line width=2.pt] (11.927355909100376,3.2796571126501264)-- (13.68,5.477403288677162);
		\draw [line width=2.pt] (12.274488566072225,4.800544655653772)-- (15.085511433927774,4.800544655653772);
		\draw [line width=2.pt] (15.085511433927774,4.800544655653772)-- (12.9,2.06);
		\draw [line width=2.pt] (12.9,2.06)-- (12.274488566072225,4.800544655653772);
		\begin{scriptsize}
			\draw [fill=uuuuuu] (1.22,2.24) circle (1.5pt);
			\draw[color=uuuuuu] (0.9962954592023646,2.170508596655724) node {$4$};
			\draw [fill=uuuuuu] (2.76,2.24) circle (1.5pt);
			\draw[color=uuuuuu] (3.0211515606857454,2.1332533172721746) node {$3$};
			\draw [fill=uuuuuu] (3.720174294862449,3.444020483000766) circle (1.5pt);
			\draw[color=uuuuuu] (3.9897888246580274,3.573790786769414) node {$2$};
			\draw [fill=uuuuuu] (3.377492056569725,4.945409467760774) circle (1.5pt);
			\draw[color=uuuuuu] (3.642072883744901,5.051583535650202) node {$1$};
			\draw [fill=uuuuuu] (1.99,5.613590426001814) circle (1.5pt);
			\draw[color=uuuuuu] (1.7538194733345338,5.771852270398822) node {$0$};
			\draw [fill=uuuuuu] (0.6025079434302745,4.945409467760775) circle (1.5pt);
			\draw[color=uuuuuu] (0.3877925626043925,5.1012572414949355) node {$6$};
			\draw [fill=uuuuuu] (0.2598257051375503,3.444020483000766) circle (1.5pt);
			\draw[color=uuuuuu] (0.02765819523008249,3.3875143898516673) node {$5$};
			\draw [fill=uuuuuu] (7.,2.18) circle (1.5pt);
			\draw[color=uuuuuu] (6.783282190113692,2.046324332043893) node {$4$};
			\draw [fill=uuuuuu] (8.5,2.16) circle (1.5pt);
			\draw[color=uuuuuu] (8.75846458575234,1.9842321997379777) node {$3$};
			\draw [fill=uuuuuu] (9.45087133243746,3.32027742766487) circle (1.5pt);
			\draw[color=uuuuuu] (9.71468342326344,3.300585404623386) node {$2$};
			\draw [fill=uuuuuu] (9.136588489746625,4.787119714616732) circle (1.5pt);
			\draw[color=uuuuuu] (9.316641188195044,5.014328256266653) node {$1$};
			\draw [fill=uuuuuu] (7.793812862675348,5.455964700651117) circle (1.5pt);
			\draw[color=uuuuuu] (7.93819585100372,5.747015417476456) node {$0$};
			\draw [fill=uuuuuu] (6.433681886039368,4.823158469332829) circle (1.5pt);
			\draw[color=uuuuuu] (6.149942440593352,5.114328256266653) node {$6$};
			\draw [fill=uuuuuu] (6.08040192686126,3.3652170197392195) circle (1.5pt);
			\draw[color=uuuuuu] (5.8394817790637745,3.2757485517010196) node {$5$};
			\draw [fill=uuuuuu] (12.9,2.06) circle (1.5pt);
			\draw[color=uuuuuu] (12.719290038559215,1.9221400674320621) node {$4$};
			\draw [fill=uuuuuu] (14.46,2.06) circle (1.5pt);
			\draw[color=uuuuuu] (14.719961875430765,1.884884788048513) node {$3$};
			\draw [fill=uuuuuu] (15.432644090899624,3.2796571126501273) circle (1.5pt);
			\draw[color=uuuuuu] (15.67552812463133,3.188819566472738) node {$2$};
			\draw [fill=uuuuuu] (15.085511433927774,4.800544655653772) circle (1.5pt);
			\draw[color=uuuuuu] (15.327812183718202,5.00190982980547) node {$1$};
			\draw [fill=uuuuuu] (13.68,5.477403288677162) circle (1.5pt);
			\draw[color=uuuuuu] (13.799693140682145,5.72217856455409) node {$0$};
			\draw [fill=uuuuuu] (12.274488566072225,4.800544655653772) circle (1.5pt);
			\draw[color=uuuuuu] (12.359155671184904,5.126746682727836) node {$6$};
			\draw [fill=uuuuuu] (11.927355909100376,3.2796571126501264) circle (1.5pt);
			\draw[color=uuuuuu] (11.688560642281017,3.176401140011555) node {$5$};
		\end{scriptsize}
	\end{tikzpicture}
	\caption{3 even graphs of order 7 and size 6}
	\label{fig:nauty-7-6-e}
\end{figure}

\begin{figure}
	\centering
	\definecolor{ccwwqq}{rgb}{0.8,0.4,0.}
	\definecolor{qqqqcc}{rgb}{0.,0.,0.8}
	\definecolor{qqwuqq}{rgb}{0.,0.39215686274509803,0.}
	\definecolor{ududff}{rgb}{0.30196078431372547,0.30196078431372547,1.}
	\definecolor{wqwqwq}{rgb}{0.3764705882352941,0.3764705882352941,0.3764705882352941}
	\definecolor{uuuuuu}{rgb}{0.26666666666666666,0.26666666666666666,0.26666666666666666}
	\begin{tikzpicture}[line cap=round,line join=round,>=triangle 45,x=2cm,y=2cm]
		\clip(-2.1671657301537617,-0.5350346502728847) rectangle (3.170541482665136,3.1822923572697652);
		\fill[line width=2.pt,color=wqwqwq,fill=wqwqwq,fill opacity=0.10000000149011612] (0.5,2.190643133767411) -- (1.400968867902419,1.756759394649853) -- (1.6234898018587334,0.7818314824680296) -- cycle;
		\fill[line width=2.pt,color=wqwqwq,fill=wqwqwq,fill opacity=0.10000000149011612] (0.,0.) -- (-0.6234898018587334,0.78183148246803) -- (-0.4009688679024189,1.7567593946498534) -- cycle;
		\fill[line width=2.pt,color=wqwqwq,fill=wqwqwq,fill opacity=0.10000000149011612] (0.,0.) -- (1.,0.) -- (1.6234898018587334,0.7818314824680296) -- cycle;
		\fill[line width=2.pt,color=wqwqwq,fill=wqwqwq,fill opacity=0.10000000149011612] (0.5,2.190643133767411) -- (1.,0.) -- (-0.4009688679024189,1.7567593946498534) -- cycle;
		\fill[line width=2.pt,color=wqwqwq,fill=wqwqwq,fill opacity=0.10000000149011612] (1.400968867902419,1.756759394649853) -- (1.,0.) -- (-0.6234898018587334,0.78183148246803) -- cycle;
		\fill[line width=2.pt,color=qqwuqq,fill=qqwuqq,fill opacity=0.10000000149011612] (1.,0.) -- (0.8755928502843471,-0.37664918434017003) -- (1.3798647595130393,-0.15442766501905125) -- cycle;
		\fill[line width=2.pt,color=qqqqcc,fill=qqqqcc,fill opacity=0.10000000149011612] (1.400968867902419,1.756759394649853) -- (-0.4009688679024189,1.7567593946498534) -- (1.1149083326301672,2.1190694172662408) -- cycle;
		\fill[line width=2.pt,color=qqqqcc,fill=qqqqcc,fill opacity=0.10000000149011612] (0.5,2.190643133767411) -- (0.,0.) -- (-0.5004711732041182,0.09343479883911965) -- cycle;
		\fill[line width=2.pt,color=ududff,fill=ududff,fill opacity=0.10000000149011612] (0.5,2.190643133767411) -- (-0.6234898018587334,0.78183148246803) -- (-0.6286758958893789,1.3498410811546757) -- cycle;
		\fill[line width=2.pt,color=ududff,fill=ududff,fill opacity=0.10000000149011612] (1.400968867902419,1.756759394649853) -- (0.,0.) -- (0.36277395954330416,-0.23135049863020776) -- cycle;
		\fill[line width=2.pt,color=ududff,fill=ududff,fill opacity=0.10000000149011612] (1.6234898018587334,0.7818314824680296) -- (-0.6234898018587334,0.78183148246803) -- (1.4567875931241958,0.21309254001202976) -- cycle;
		\fill[line width=2.pt,color=ududff,fill=ududff,fill opacity=0.10000000149011612] (1.6234898018587334,0.7818314824680296) -- (-0.4009688679024189,1.7567593946498534) -- (1.7046500569823666,1.4438578777905335) -- cycle;
		\fill[line width=2.pt,color=ccwwqq,fill=ccwwqq,fill opacity=0.10000000149011612] (-0.17093330413708208,2.0384566722231297) -- (0.22602225534569254,2.281462066000904) -- (-0.1829043686638885,2.503732967797565) -- cycle;
		\draw [line width=2.pt,color=wqwqwq] (0.5,2.190643133767411)-- (1.400968867902419,1.756759394649853);
		\draw [line width=2.pt,color=wqwqwq] (1.400968867902419,1.756759394649853)-- (1.6234898018587334,0.7818314824680296);
		\draw [line width=2.pt,color=wqwqwq] (1.6234898018587334,0.7818314824680296)-- (0.5,2.190643133767411);
		\draw [line width=2.pt,color=wqwqwq] (0.,0.)-- (-0.6234898018587334,0.78183148246803);
		\draw [line width=2.pt,color=wqwqwq] (-0.6234898018587334,0.78183148246803)-- (-0.4009688679024189,1.7567593946498534);
		\draw [line width=2.pt,color=wqwqwq] (-0.4009688679024189,1.7567593946498534)-- (0.,0.);
		\draw [line width=2.pt,color=wqwqwq] (0.,0.)-- (1.,0.);
		\draw [line width=2.pt,color=wqwqwq] (1.,0.)-- (1.6234898018587334,0.7818314824680296);
		\draw [line width=2.pt,color=wqwqwq] (1.6234898018587334,0.7818314824680296)-- (0.,0.);
		\draw [line width=2.pt,color=wqwqwq] (0.5,2.190643133767411)-- (1.,0.);
		\draw [line width=2.pt,color=wqwqwq] (1.,0.)-- (-0.4009688679024189,1.7567593946498534);
		\draw [line width=2.pt,color=wqwqwq] (-0.4009688679024189,1.7567593946498534)-- (0.5,2.190643133767411);
		\draw [line width=2.pt,color=wqwqwq] (1.400968867902419,1.756759394649853)-- (1.,0.);
		\draw [line width=2.pt,color=wqwqwq] (1.,0.)-- (-0.6234898018587334,0.78183148246803);
		\draw [line width=2.pt,color=wqwqwq] (-0.6234898018587334,0.78183148246803)-- (1.400968867902419,1.756759394649853);
		\draw [line width=2.pt,color=qqwuqq] (1.,0.)-- (0.8755928502843471,-0.37664918434017003);
		\draw [line width=2.pt,color=qqwuqq] (0.8755928502843471,-0.37664918434017003)-- (1.3798647595130393,-0.15442766501905125);
		\draw [line width=2.pt,color=qqwuqq] (1.3798647595130393,-0.15442766501905125)-- (1.,0.);
		\draw [line width=2.pt,color=qqqqcc] (1.400968867902419,1.756759394649853)-- (-0.4009688679024189,1.7567593946498534);
		\draw [line width=2.pt,color=qqqqcc] (-0.4009688679024189,1.7567593946498534)-- (1.1149083326301672,2.1190694172662408);
		\draw [line width=2.pt,color=qqqqcc] (1.1149083326301672,2.1190694172662408)-- (1.400968867902419,1.756759394649853);
		\draw [line width=2.pt,color=qqqqcc] (0.5,2.190643133767411)-- (0.,0.);
		\draw [line width=2.pt,color=qqqqcc] (0.,0.)-- (-0.5004711732041182,0.09343479883911965);
		\draw [line width=2.pt,color=qqqqcc] (-0.5004711732041182,0.09343479883911965)-- (0.5,2.190643133767411);
		\draw [line width=2.pt,color=ududff] (0.5,2.190643133767411)-- (-0.6234898018587334,0.78183148246803);
		\draw [line width=2.pt,color=ududff] (-0.6234898018587334,0.78183148246803)-- (-0.6286758958893789,1.3498410811546757);
		\draw [line width=2.pt,color=ududff] (-0.6286758958893789,1.3498410811546757)-- (0.5,2.190643133767411);
		\draw [line width=2.pt,color=ududff] (1.400968867902419,1.756759394649853)-- (0.,0.);
		\draw [line width=2.pt,color=ududff] (0.,0.)-- (0.36277395954330416,-0.23135049863020776);
		\draw [line width=2.pt,color=ududff] (0.36277395954330416,-0.23135049863020776)-- (1.400968867902419,1.756759394649853);
		\draw [line width=2.pt,color=ududff] (1.6234898018587334,0.7818314824680296)-- (-0.6234898018587334,0.78183148246803);
		\draw [line width=2.pt,color=ududff] (-0.6234898018587334,0.78183148246803)-- (1.4567875931241958,0.21309254001202976);
		\draw [line width=2.pt,color=ududff] (1.4567875931241958,0.21309254001202976)-- (1.6234898018587334,0.7818314824680296);
		\draw [line width=2.pt,color=ududff] (1.6234898018587334,0.7818314824680296)-- (-0.4009688679024189,1.7567593946498534);
		\draw [line width=2.pt,color=ududff] (-0.4009688679024189,1.7567593946498534)-- (1.7046500569823666,1.4438578777905335);
		\draw [line width=2.pt,color=ududff] (1.7046500569823666,1.4438578777905335)-- (1.6234898018587334,0.7818314824680296);
		\draw [line width=2.pt,color=ccwwqq] (-0.17093330413708208,2.0384566722231297)-- (0.22602225534569254,2.281462066000904);
		\draw [line width=2.pt,color=ccwwqq] (0.22602225534569254,2.281462066000904)-- (-0.1829043686638885,2.503732967797565);
		\draw [line width=2.pt,color=ccwwqq] (-0.1829043686638885,2.503732967797565)-- (-0.17093330413708208,2.0384566722231297);
		\begin{scriptsize}
			\draw [fill=uuuuuu] (0.,0.) circle (2.0pt);
			\draw[color=uuuuuu] (-0.06156114116360256,-0.15290940402956188) node {$w_5$};
			\draw [fill=uuuuuu] (1.,0.) circle (2.0pt);
			\draw[color=uuuuuu] (1.1635517933945184,0.02163820126752297) node {$w_4$};
			\draw [fill=uuuuuu] (1.6234898018587334,0.7818314824680296) circle (2.0pt);
			\draw[color=uuuuuu] (1.7909665008662134,0.7395188880239269) node {$w_3$};
			\draw [fill=uuuuuu] (1.400968867902419,1.756759394649853) circle (2.0pt);
			\draw[color=uuuuuu] (1.5662209340106809,1.8134163443979026) node {$w_2$};
			\draw [fill=uuuuuu] (0.5,2.190643133767411) circle (2.0pt);
			\draw[color=uuuuuu] (0.5672386665885508,2.294178680871007) node {$w_1$};
			\draw [fill=uuuuuu] (-0.4009688679024189,1.7567593946498534) circle (2.0pt);
			\draw[color=uuuuuu] (-0.43174360083357934,1.8334163443979026) node {$w_7$};
			\draw [fill=uuuuuu] (-0.6234898018587334,0.78183148246803) circle (2.0pt);
			\draw[color=uuuuuu] (-0.7564891676891119,0.8395188880239269) node {$w_6$};
			\draw [fill=ududff] (0.8755928502843471,-0.37664918434017003) circle (1.0pt);
			\draw [fill=ududff] (1.3798647595130393,-0.15442766501905125) circle (1.0pt);
			\draw [fill=ududff] (1.1149083326301672,2.1190694172662408) circle (1.0pt);
			\draw [fill=ududff] (-0.5004711732041182,0.09343479883911965) circle (1.0pt);
			\draw [fill=ududff] (-0.6286758958893789,1.3498410811546757) circle (1.0pt);
			\draw [fill=ududff] (0.36277395954330416,-0.23135049863020776) circle (1.0pt);
			\draw [fill=ududff] (1.4567875931241958,0.21309254001202976) circle (1.0pt);
			\draw [fill=ududff] (1.7046500569823666,1.4438578777905335) circle (1.0pt);
			\draw [fill=ududff] (-0.17093330413708208,2.0384566722231297) circle (1.0pt);
			\draw [fill=ududff] (0.22602225534569254,2.281462066000904) circle (1.0pt);
			\draw [fill=ududff] (-0.1829043686638885,2.503732967797565) circle (1.0pt);
		\end{scriptsize}
	\end{tikzpicture}
	\caption{13 3-cliques in $K_{18}$ for case \eqref{eq:c2}}
	\label{fig:k18alpha13-5}
\end{figure}

\begin{lemma}\label{lem:k18a13c3}
	In a $\{K_3,K_4\}$-decomposition of $K_{18}$, if $\alpha=13$, then the case \eqref{eq:c3} is impossible.
\end{lemma}
\begin{proof}
For contradiction, suppose case \eqref{eq:c3} is possible. Since $t_3=5$, as in case \eqref{eq:c2}, it remains to consider the first two graphs $H$ in Figure \ref{fig:nauty-7-6-e}. Their complements in the $K_7$ on $w_1,\ldots,w_7$ give the union of the five 3-cliques lying entirely in this vertex set. 
	
	If $H$ is a bowtie as the first graph in Figure \ref{fig:nauty-7-6-e}, we take the 5 3-cliques $\{w_1,w_2,w_3\}$, $\{w_1,w_4,w_6\}$, $\{w_2,w_4,w_5\}$, $\{w_3,w_5,w_6\}$, $\{w_3,w_4,w_7\}$ from $\overline{H}$, and $w_3$ and $w_4$ each need an additional 3-clique, which we take, without loss of generality, to be $\{w_3,y_1,y_2\}$ and $\{w_4, y_3,y_4\}$. Since $t_2=5$ in \eqref{eq:c3}, we need to select 5 edges out of 6 in $H$ to form 3-cliques with vertices among $y_1,\cdots,y_{11}$. By counting the 3-cliques containing each vertex, we deduce that precisely 1 edge incident to $w_7$ must not be chosen. Therefore, without loss of generality, we select all edges except $\{w_6,w_7\}$, and they form the 3-cliques $\{w_1,w_5,y_5\}$, $\{w_1,w_7,y_6\}$, $\{w_2,w_6,y_7\}$, $\{w_2,w_7,y_8\}$, $\{w_5,w_7,y_9\}$. We append a final 3-clique for $w_6$: $\{w_6,y_{10},y_{11}\}$. This configuration (See Figure \ref{fig:k18alpha13-6} up) is then excluded by an exhaustive computer-assisted search. 
	
	If $H\cong C_6$ is a cycle as the second graph in Figure \ref{fig:nauty-7-6-e}, we take the 5 3-cliques $\{w_1,w_2,w_3\}$, $\{w_1,w_4,w_7\}$, $\{w_2,w_4,w_6\}$, $\{w_3,w_4,w_5\}$, $\{w_5,w_6,w_7\}$ from $\overline{H}$, and we also notice that $w_4$ requires an additional 3-clique, which we can take, without loss of generality, to be $\{w_4,y_1,y_2\}$. Since $t_2=5$ in \eqref{eq:c3}, we need to select 5 edges out of 6 in $H$ to form 3-cliques with vertices among $y_1,\cdots,y_{11}$. As $H\cong C_6$ is a cycle, we can take, without loss of generality, all edges of $H$ except $\{w_3,w_7\}$. Therefore, the 3-cliques can be taken to be $\{w_1,w_5,y_3\}$, $\{w_1,w_6,y_4\}$, $\{w_2,w_5,y_5\}$, $\{w_2,w_7,y_6\}$ and $\{w_3,w_6,y_7\}$. We append two 3-cliques for $w_3$ and $w_7$ respectively: $\{w_3,y_8,y_9\}$, $\{w_7,y_{10},y_{11}\}$. This configuration (See Figure \ref{fig:k18alpha13-6} down) is then excluded by an exhaustive computer-assisted search. 
	
	To sum up, the case \eqref{eq:c3} is impossible.
\end{proof}

\begin{figure}
	\centering
	\definecolor{qqwuqq}{rgb}{0.,0.39215686274509803,0.}
	\definecolor{ududff}{rgb}{0.30196078431372547,0.30196078431372547,1.}
	\definecolor{wqwqwq}{rgb}{0.3764705882352941,0.3764705882352941,0.3764705882352941}
	\definecolor{uuuuuu}{rgb}{0.26666666666666666,0.26666666666666666,0.26666666666666666}
	\begin{tikzpicture}[line cap=round,line join=round,>=triangle 45,x=2.0cm,y=2.0cm]
		\clip(-1.5,-0.5) rectangle (2.5,2.5);
		\fill[line width=2.pt,color=wqwqwq,fill=wqwqwq,fill opacity=0.10000000149011612] (0.5,2.190643133767411) -- (1.400968867902419,1.756759394649853) -- (1.6234898018587334,0.7818314824680296) -- cycle;
		\fill[line width=2.pt,color=wqwqwq,fill=wqwqwq,fill opacity=0.10000000149011612] (0.5,2.190643133767411) -- (1.,0.) -- (-0.6234898018587334,0.78183148246803) -- cycle;
		\fill[line width=2.pt,color=wqwqwq,fill=wqwqwq,fill opacity=0.10000000149011612] (1.400968867902419,1.756759394649853) -- (1.,0.) -- (0.,0.) -- cycle;
		\fill[line width=2.pt,color=wqwqwq,fill=wqwqwq,fill opacity=0.10000000149011612] (1.6234898018587334,0.7818314824680296) -- (0.,0.) -- (-0.6234898018587334,0.78183148246803) -- cycle;
		\fill[line width=2.pt,color=wqwqwq,fill=wqwqwq,fill opacity=0.10000000149011612] (1.6234898018587334,0.7818314824680296) -- (1.,0.) -- (-0.4009688679024189,1.7567593946498534) -- cycle;
		\fill[line width=2.pt,color=ududff,fill=ududff,fill opacity=0.10000000149011612] (0.5,2.190643133767411) -- (0.,0.) -- (-0.5817730439603593,0.10515097107276046) -- cycle;
		\fill[line width=2.pt,color=ududff,fill=ududff,fill opacity=0.10000000149011612] (0.5,2.190643133767411) -- (-0.4009688679024189,1.7567593946498534) -- (-0.07703195839730917,2.22712369813375) -- cycle;
		\fill[line width=2.pt,color=ududff,fill=ududff,fill opacity=0.10000000149011612] (1.400968867902419,1.756759394649853) -- (-0.6234898018587334,0.78183148246803) -- (1.1693694977889983,2.196221182691114) -- cycle;
		\fill[line width=2.pt,color=ududff,fill=ududff,fill opacity=0.10000000149011612] (1.400968867902419,1.756759394649853) -- (-0.4009688679024189,1.7567593946498534) -- (-0.6744805902882665,1.3618532657399485) -- cycle;
		\fill[line width=2.pt,color=ududff,fill=ududff,fill opacity=0.10000000149011612] (0.,0.) -- (-0.4009688679024189,1.7567593946498534) -- (-0.5714722054794807,1.1352348191606194) -- cycle;
		\fill[line width=2.pt,color=qqwuqq,fill=qqwuqq,fill opacity=0.10000000149011612] (1.6234898018587334,0.7818314824680296) -- (2.,1.) -- (1.88012735296962,0.4759811563843897) -- cycle;
		\fill[line width=2.pt,color=qqwuqq,fill=qqwuqq,fill opacity=0.10000000149011612] (1.,0.) -- (0.8088401509582482,-0.30688256816238313) -- (1.2517762056360269,-0.18327250639184006) -- cycle;
		\fill[line width=2.pt,color=qqwuqq,fill=qqwuqq,fill opacity=0.10000000149011612] (-0.6234898018587334,0.78183148246803) -- (-1.,1.) -- (-0.9938065831955023,0.5789895411931756) -- cycle;
		\draw [line width=2.pt,color=wqwqwq] (0.5,2.190643133767411)-- (1.400968867902419,1.756759394649853);
		\draw [line width=2.pt,color=wqwqwq] (1.400968867902419,1.756759394649853)-- (1.6234898018587334,0.7818314824680296);
		\draw [line width=2.pt,color=wqwqwq] (1.6234898018587334,0.7818314824680296)-- (0.5,2.190643133767411);
		\draw [line width=2.pt,color=wqwqwq] (0.5,2.190643133767411)-- (1.,0.);
		\draw [line width=2.pt,color=wqwqwq] (1.,0.)-- (-0.6234898018587334,0.78183148246803);
		\draw [line width=2.pt,color=wqwqwq] (-0.6234898018587334,0.78183148246803)-- (0.5,2.190643133767411);
		\draw [line width=2.pt,color=wqwqwq] (1.400968867902419,1.756759394649853)-- (1.,0.);
		\draw [line width=2.pt,color=wqwqwq] (1.,0.)-- (0.,0.);
		\draw [line width=2.pt,color=wqwqwq] (0.,0.)-- (1.400968867902419,1.756759394649853);
		\draw [line width=2.pt,color=wqwqwq] (1.6234898018587334,0.7818314824680296)-- (0.,0.);
		\draw [line width=2.pt,color=wqwqwq] (0.,0.)-- (-0.6234898018587334,0.78183148246803);
		\draw [line width=2.pt,color=wqwqwq] (-0.6234898018587334,0.78183148246803)-- (1.6234898018587334,0.7818314824680296);
		\draw [line width=2.pt,color=wqwqwq] (1.6234898018587334,0.7818314824680296)-- (1.,0.);
		\draw [line width=2.pt,color=wqwqwq] (1.,0.)-- (-0.4009688679024189,1.7567593946498534);
		\draw [line width=2.pt,color=wqwqwq] (-0.4009688679024189,1.7567593946498534)-- (1.6234898018587334,0.7818314824680296);
		\draw [line width=2.pt,color=ududff] (0.5,2.190643133767411)-- (0.,0.);
		\draw [line width=2.pt,color=ududff] (0.,0.)-- (-0.5817730439603593,0.10515097107276046);
		\draw [line width=2.pt,color=ududff] (-0.5817730439603593,0.10515097107276046)-- (0.5,2.190643133767411);
		\draw [line width=2.pt,color=ududff] (0.5,2.190643133767411)-- (-0.4009688679024189,1.7567593946498534);
		\draw [line width=2.pt,color=ududff] (-0.4009688679024189,1.7567593946498534)-- (-0.07703195839730917,2.22712369813375);
		\draw [line width=2.pt,color=ududff] (-0.07703195839730917,2.22712369813375)-- (0.5,2.190643133767411);
		\draw [line width=2.pt,color=ududff] (1.400968867902419,1.756759394649853)-- (-0.6234898018587334,0.78183148246803);
		\draw [line width=2.pt,color=ududff] (-0.6234898018587334,0.78183148246803)-- (1.1693694977889983,2.196221182691114);
		\draw [line width=2.pt,color=ududff] (1.1693694977889983,2.196221182691114)-- (1.400968867902419,1.756759394649853);
		\draw [line width=2.pt,color=ududff] (1.400968867902419,1.756759394649853)-- (-0.4009688679024189,1.7567593946498534);
		\draw [line width=2.pt,color=ududff] (-0.4009688679024189,1.7567593946498534)-- (-0.6744805902882665,1.3618532657399485);
		\draw [line width=2.pt,color=ududff] (-0.6744805902882665,1.3618532657399485)-- (1.400968867902419,1.756759394649853);
		\draw [line width=2.pt,color=ududff] (0.,0.)-- (-0.4009688679024189,1.7567593946498534);
		\draw [line width=2.pt,color=ududff] (-0.4009688679024189,1.7567593946498534)-- (-0.5714722054794807,1.1352348191606194);
		\draw [line width=2.pt,color=ududff] (-0.5714722054794807,1.1352348191606194)-- (0.,0.);
		\draw [line width=2.pt,color=qqwuqq] (1.6234898018587334,0.7818314824680296)-- (2.,1.);
		\draw [line width=2.pt,color=qqwuqq] (2.,1.)-- (1.88012735296962,0.4759811563843897);
		\draw [line width=2.pt,color=qqwuqq] (1.88012735296962,0.4759811563843897)-- (1.6234898018587334,0.7818314824680296);
		\draw [line width=2.pt,color=qqwuqq] (1.,0.)-- (0.8088401509582482,-0.30688256816238313);
		\draw [line width=2.pt,color=qqwuqq] (0.8088401509582482,-0.30688256816238313)-- (1.2517762056360269,-0.18327250639184006);
		\draw [line width=2.pt,color=qqwuqq] (1.2517762056360269,-0.18327250639184006)-- (1.,0.);
		\draw [line width=2.pt,color=qqwuqq] (-0.6234898018587334,0.78183148246803)-- (-1.,1.);
		\draw [line width=2.pt,color=qqwuqq] (-1.,1.)-- (-0.9938065831955023,0.5789895411931756);
		\draw [line width=2.pt,color=qqwuqq] (-0.9938065831955023,0.5789895411931756)-- (-0.6234898018587334,0.78183148246803);
		\begin{scriptsize}
			\draw [fill=uuuuuu] (0.,0.) circle (2.0pt);
			\draw[color=uuuuuu] (-0.11363824390786975,-0.10129145302749142) node {$w_5$};
			\draw [fill=uuuuuu] (1.,0.) circle (2.0pt);
			\draw[color=uuuuuu] (1.1618267724165396,0.02853892487619718) node {$w_4$};
			\draw [fill=uuuuuu] (1.6234898018587334,0.7818314824680296) circle (2.0pt);
			\draw[color=uuuuuu] (1.6968686964604684,0.9756143881552705) node {$w_3$};
			\draw [fill=uuuuuu] (1.400968867902419,1.756759394649853) circle (2.0pt);
			\draw[color=uuuuuu] (1.4689451226998385,1.9120542500533242) node {$w_2$};
			\draw [fill=uuuuuu] (0.5,2.190643133767411) circle (2.0pt);
			\draw[color=uuuuuu] (0.5699628552777084,2.342816586526429) node {$w_1$};
			\draw [fill=uuuuuu] (-0.4009688679024189,1.7567593946498534) circle (2.0pt);
			\draw[color=uuuuuu] (-0.5444005803809737,1.7996814666255578) node {$w_7$};
			\draw [fill=uuuuuu] (-0.6234898018587334,0.78183148246803) circle (2.0pt);
			\draw[color=uuuuuu] (-0.6942309582846621,0.6372248351099321) node {$w_6$};
			\draw [fill=ududff] (-0.5817730439603593,0.10515097107276046) circle (1.0pt);
			\draw [fill=ududff] (-0.07703195839730917,2.22712369813375) circle (1.0pt);
			\draw [fill=ududff] (1.1693694977889983,2.196221182691114) circle (1.0pt);
			\draw [fill=ududff] (-0.6744805902882665,1.3618532657399485) circle (1.0pt);
			\draw [fill=ududff] (-0.5714722054794807,1.1352348191606194) circle (1.0pt);
			\draw [fill=ududff] (2.,1.) circle (1.0pt);
			\draw [fill=ududff] (1.88012735296962,0.4759811563843897) circle (1.0pt);
			\draw [fill=ududff] (0.8088401509582482,-0.30688256816238313) circle (1.0pt);
			\draw [fill=ududff] (1.2517762056360269,-0.18327250639184006) circle (1.0pt);
			\draw [fill=ududff] (-1.,1.) circle (1.0pt);
			\draw [fill=ududff] (-0.9938065831955023,0.5789895411931756) circle (1.0pt);
		\end{scriptsize}
	\end{tikzpicture}
	\definecolor{qqqqcc}{rgb}{0.,0.,0.8}
	\definecolor{qqwuqq}{rgb}{0.,0.39215686274509803,0.}
	\definecolor{ududff}{rgb}{0.30196078431372547,0.30196078431372547,1.}
	\definecolor{wqwqwq}{rgb}{0.3764705882352941,0.3764705882352941,0.3764705882352941}
	\definecolor{uuuuuu}{rgb}{0.26666666666666666,0.26666666666666666,0.26666666666666666}
	\begin{tikzpicture}[line cap=round,line join=round,>=triangle 45,x=2.0cm,y=2.0cm]
		\clip(-1,-0.5) rectangle (2,2.5);
		\fill[line width=2.pt,color=wqwqwq,fill=wqwqwq,fill opacity=0.10000000149011612] (0.5,2.190643133767411) -- (1.400968867902419,1.756759394649853) -- (1.6234898018587334,0.7818314824680296) -- cycle;
		\fill[line width=2.pt,color=wqwqwq,fill=wqwqwq,fill opacity=0.10000000149011612] (0.,0.) -- (-0.6234898018587334,0.78183148246803) -- (-0.4009688679024189,1.7567593946498534) -- cycle;
		\fill[line width=2.pt,color=wqwqwq,fill=wqwqwq,fill opacity=0.10000000149011612] (0.,0.) -- (1.,0.) -- (1.6234898018587334,0.7818314824680296) -- cycle;
		\fill[line width=2.pt,color=wqwqwq,fill=wqwqwq,fill opacity=0.10000000149011612] (0.5,2.190643133767411) -- (1.,0.) -- (-0.4009688679024189,1.7567593946498534) -- cycle;
		\fill[line width=2.pt,color=wqwqwq,fill=wqwqwq,fill opacity=0.10000000149011612] (1.400968867902419,1.756759394649853) -- (1.,0.) -- (-0.6234898018587334,0.78183148246803) -- cycle;
		\fill[line width=2.pt,color=qqwuqq,fill=qqwuqq,fill opacity=0.10000000149011612] (1.,0.) -- (0.8755928502843471,-0.37664918434017003) -- (1.3798647595130393,-0.15442766501905125) -- cycle;
		\fill[line width=2.pt,color=qqqqcc,fill=qqqqcc,fill opacity=0.10000000149011612] (1.400968867902419,1.756759394649853) -- (-0.4009688679024189,1.7567593946498534) -- (1.1149083326301672,2.1190694172662408) -- cycle;
		\fill[line width=2.pt,color=qqqqcc,fill=qqqqcc,fill opacity=0.10000000149011612] (0.5,2.190643133767411) -- (0.,0.) -- (-0.5004711732041182,0.09343479883911965) -- cycle;
		\fill[line width=2.pt,color=ududff,fill=ududff,fill opacity=0.10000000149011612] (0.5,2.190643133767411) -- (-0.6234898018587334,0.78183148246803) -- (-0.6286758958893789,1.3498410811546757) -- cycle;
		\fill[line width=2.pt,color=ududff,fill=ududff,fill opacity=0.10000000149011612] (1.400968867902419,1.756759394649853) -- (0.,0.) -- (0.36277395954330416,-0.23135049863020776) -- cycle;
		\fill[line width=2.pt,color=ududff,fill=ududff,fill opacity=0.10000000149011612] (1.6234898018587334,0.7818314824680296) -- (-0.6234898018587334,0.78183148246803) -- (1.4567875931241958,0.21309254001202976) -- cycle;
		\fill[line width=2.pt,color=qqwuqq,fill=qqwuqq,fill opacity=0.10000000149011612] (1.6234898018587334,0.7818314824680296) -- (1.9063476691027654,0.5977183384311374) -- (1.8501612773888825,1.2251330459028333) -- cycle;
		\fill[line width=2.pt,color=qqwuqq,fill=qqwuqq,fill opacity=0.10000000149011612] (-0.4009688679024189,1.7567593946498534) -- (-0.2474640132627547,2.0960221174680234) -- (-0.7156839442117808,1.796361361660646) -- cycle;
		\draw [line width=2.pt,color=wqwqwq] (0.5,2.190643133767411)-- (1.400968867902419,1.756759394649853);
		\draw [line width=2.pt,color=wqwqwq] (1.400968867902419,1.756759394649853)-- (1.6234898018587334,0.7818314824680296);
		\draw [line width=2.pt,color=wqwqwq] (1.6234898018587334,0.7818314824680296)-- (0.5,2.190643133767411);
		\draw [line width=2.pt,color=wqwqwq] (0.,0.)-- (-0.6234898018587334,0.78183148246803);
		\draw [line width=2.pt,color=wqwqwq] (-0.6234898018587334,0.78183148246803)-- (-0.4009688679024189,1.7567593946498534);
		\draw [line width=2.pt,color=wqwqwq] (-0.4009688679024189,1.7567593946498534)-- (0.,0.);
		\draw [line width=2.pt,color=wqwqwq] (0.,0.)-- (1.,0.);
		\draw [line width=2.pt,color=wqwqwq] (1.,0.)-- (1.6234898018587334,0.7818314824680296);
		\draw [line width=2.pt,color=wqwqwq] (1.6234898018587334,0.7818314824680296)-- (0.,0.);
		\draw [line width=2.pt,color=wqwqwq] (0.5,2.190643133767411)-- (1.,0.);
		\draw [line width=2.pt,color=wqwqwq] (1.,0.)-- (-0.4009688679024189,1.7567593946498534);
		\draw [line width=2.pt,color=wqwqwq] (-0.4009688679024189,1.7567593946498534)-- (0.5,2.190643133767411);
		\draw [line width=2.pt,color=wqwqwq] (1.400968867902419,1.756759394649853)-- (1.,0.);
		\draw [line width=2.pt,color=wqwqwq] (1.,0.)-- (-0.6234898018587334,0.78183148246803);
		\draw [line width=2.pt,color=wqwqwq] (-0.6234898018587334,0.78183148246803)-- (1.400968867902419,1.756759394649853);
		\draw [line width=2.pt,color=qqwuqq] (1.,0.)-- (0.8755928502843471,-0.37664918434017003);
		\draw [line width=2.pt,color=qqwuqq] (0.8755928502843471,-0.37664918434017003)-- (1.3798647595130393,-0.15442766501905125);
		\draw [line width=2.pt,color=qqwuqq] (1.3798647595130393,-0.15442766501905125)-- (1.,0.);
		\draw [line width=2.pt,color=qqqqcc] (1.400968867902419,1.756759394649853)-- (-0.4009688679024189,1.7567593946498534);
		\draw [line width=2.pt,color=qqqqcc] (-0.4009688679024189,1.7567593946498534)-- (1.1149083326301672,2.1190694172662408);
		\draw [line width=2.pt,color=qqqqcc] (1.1149083326301672,2.1190694172662408)-- (1.400968867902419,1.756759394649853);
		\draw [line width=2.pt,color=qqqqcc] (0.5,2.190643133767411)-- (0.,0.);
		\draw [line width=2.pt,color=qqqqcc] (0.,0.)-- (-0.5004711732041182,0.09343479883911965);
		\draw [line width=2.pt,color=qqqqcc] (-0.5004711732041182,0.09343479883911965)-- (0.5,2.190643133767411);
		\draw [line width=2.pt,color=ududff] (0.5,2.190643133767411)-- (-0.6234898018587334,0.78183148246803);
		\draw [line width=2.pt,color=ududff] (-0.6234898018587334,0.78183148246803)-- (-0.6286758958893789,1.3498410811546757);
		\draw [line width=2.pt,color=ududff] (-0.6286758958893789,1.3498410811546757)-- (0.5,2.190643133767411);
		\draw [line width=2.pt,color=ududff] (1.400968867902419,1.756759394649853)-- (0.,0.);
		\draw [line width=2.pt,color=ududff] (0.,0.)-- (0.36277395954330416,-0.23135049863020776);
		\draw [line width=2.pt,color=ududff] (0.36277395954330416,-0.23135049863020776)-- (1.400968867902419,1.756759394649853);
		\draw [line width=2.pt,color=ududff] (1.6234898018587334,0.7818314824680296)-- (-0.6234898018587334,0.78183148246803);
		\draw [line width=2.pt,color=ududff] (-0.6234898018587334,0.78183148246803)-- (1.4567875931241958,0.21309254001202976);
		\draw [line width=2.pt,color=ududff] (1.4567875931241958,0.21309254001202976)-- (1.6234898018587334,0.7818314824680296);
		\draw [line width=2.pt,color=qqwuqq] (1.6234898018587334,0.7818314824680296)-- (1.9063476691027654,0.5977183384311374);
		\draw [line width=2.pt,color=qqwuqq] (1.9063476691027654,0.5977183384311374)-- (1.8501612773888825,1.2251330459028333);
		\draw [line width=2.pt,color=qqwuqq] (1.8501612773888825,1.2251330459028333)-- (1.6234898018587334,0.7818314824680296);
		\draw [line width=2.pt,color=qqwuqq] (-0.4009688679024189,1.7567593946498534)-- (-0.2474640132627547,2.0960221174680234);
		\draw [line width=2.pt,color=qqwuqq] (-0.2474640132627547,2.0960221174680234)-- (-0.7156839442117808,1.796361361660646);
		\draw [line width=2.pt,color=qqwuqq] (-0.7156839442117808,1.796361361660646)-- (-0.4009688679024189,1.7567593946498534);
		\begin{scriptsize}
			\draw [fill=uuuuuu] (0.,0.) circle (2.0pt);
			\draw[color=uuuuuu] (-0.11636243259702736,-0.09992935868291314) node {$w_5$};
			\draw [fill=uuuuuu] (1.,0.) circle (2.0pt);
			\draw[color=uuuuuu] (1.229102583727382,0.04990101922077546) node {$w_4$};
			\draw [fill=uuuuuu] (1.6234898018587334,0.7818314824680296) circle (2.0pt);
			\draw[color=uuuuuu] (1.6941445077713108,0.5785869294545104) node {$w_3$};
			\draw [fill=uuuuuu] (1.400968867902419,1.756759394649853) circle (2.0pt);
			\draw[color=uuuuuu] (1.4662209340106809,1.9134163443979022) node {$w_2$};
			\draw [fill=uuuuuu] (0.5,2.190643133767411) circle (2.0pt);
			\draw[color=uuuuuu] (0.5672386665885508,2.344178680871007) node {$w_1$};
			\draw [fill=uuuuuu] (-0.4009688679024189,1.7567593946498534) circle (2.0pt);
			\draw[color=uuuuuu] (-0.5283959718321704,1.6793063789233889) node {$w_7$};
			\draw [fill=uuuuuu] (-0.6234898018587334,0.78183148246803) circle (2.0pt);
			\draw[color=uuuuuu] (-0.7718703359256639,0.705408922549413) node {$w_6$};
			\draw [fill=ududff] (0.8755928502843471,-0.37664918434017003) circle (1pt);
			\draw [fill=ududff] (1.3798647595130393,-0.15442766501905125) circle (1pt);
			\draw [fill=ududff] (1.1149083326301672,2.1190694172662408) circle (1pt);
			\draw [fill=ududff] (-0.5004711732041182,0.09343479883911965) circle (1pt);
			\draw [fill=ududff] (-0.6286758958893789,1.3498410811546757) circle (1pt);
			\draw [fill=ududff] (0.36277395954330416,-0.23135049863020776) circle (1pt);
			\draw [fill=ududff] (1.4567875931241958,0.21309254001202976) circle (1pt);
			\draw [fill=ududff] (1.9063476691027654,0.5977183384311374) circle (1pt);
			\draw [fill=ududff] (1.8501612773888825,1.2251330459028333) circle (1pt);
			\draw [fill=ududff] (-0.2474640132627547,2.0960221174680234) circle (1pt);
			\draw [fill=ududff] (-0.7156839442117808,1.796361361660646) circle (1pt);
		\end{scriptsize}
	\end{tikzpicture}
	\caption{13 3-cliques in $K_{18}$ for case \eqref{eq:c3}}
	\label{fig:k18alpha13-6}
\end{figure}

\begin{lemma}\label{lem:k18alpha13c4}
In a $\{K_3,K_4\}$-decomposition of $K_{18}$, if $\alpha=13$, then the case \eqref{eq:c4} is impossible.
\end{lemma}
\begin{proof}
For contradiction, suppose case \eqref{eq:c4} is possible. Then four 3-cliques lie entirely among $w_1,\ldots,w_7$, so the complement of their union in the induced $K_7$ is an even graph with nine edges. An isomorphism-free computer enumeration gives exactly six even graphs $G_1,G_2,\ldots,G_6$ of order 7 and size 9, shown in Figure \ref{fig:nauty-7-9-e}, where vertex $i$ represents $w_{i+1}$ for $i\in\{0,1,\ldots,6\}$. An exhaustive computer-assisted search shows that the complements of $G_1$, $G_2$, and $G_6$ are not decomposable into 3-cliques. Hence, up to isomorphism, 
	\begin{align*}
		H\in \{G_3,G_4,G_5\}.
	\end{align*}
	If $H\cong G_3$, then $w_4$ is already contained in three 3-cliques in $\overline{H}$, so it needs another 3-clique formed with two vertices among $y_1,\cdots,y_{11}$, which is impossible for $t_1=0$ in \eqref{eq:c4}. 
	
	If $H\cong G_4$, then $w_7$ is incident to 6 edges in~$H$, but it is only contained in 4 3-cliques. Therefore, at least 2 of the 6 edges in~$H$ do not contribute to 3-cliques. In other words, $|E(H)| = 9$ and $t_2\le 9-2=7$. However, $t_2=8$ in \eqref{eq:c4}, a contradiction.
	
	If $H\cong G_5$, we first take the 4 3-cliques in $\overline{H}$: $\{w_1,w_2,w_7\}$, $\{w_1,w_3,w_5\}$, $\{w_2,w_3,w_4\}$, $\{w_4,w_5,w_6\}$. Since $w_6$ and $w_7$ are each in 1 3-clique in $\overline{H}$, we need to take in $H$ 3 edges adjacent to $w_6$ and 3 adjacent to $w_7$ to form 3-cliques. Since $t_2=8$ in $\eqref{eq:c4}$, it is obvious that all edges in $H$ except $\{w_6,w_7\}$ contribute to 3-cliques, which we may take as: $\{w_1,w_4,y_1\}$, $\{w_1,w_6,y_2\}$, $\{w_2,w_5,y_3\}$, $\{w_2,w_6,y_4\}$, $\{w_3,w_6,y_5\}$, $\{w_3,w_7,y_6\}$, $\{w_4,w_7,y_7\}$ and $\{w_5,w_7,y_8\}$. We append a final 3-clique out of the $K_7$: $\{y_9,y_{10},y_{11}\}$. This configuration (See Figure \ref{fig:k18alpha13-7}) is then shown infeasible by an exhaustive computer-assisted search. 
\end{proof}

\begin{figure}
	\centering
	\definecolor{uuuuuu}{rgb}{0.26666666666666666,0.26666666666666666,0.26666666666666666}
	\begin{tikzpicture}[line cap=round,line join=round,>=triangle 45,x=0.6cm,y=0.6cm]
		\clip(-10.,-6.) rectangle (10.,6.);
		\draw [line width=2.pt] (-8.801937735804838,4.513518789299707)-- (-5.198062264195162,4.513518789299706);
		\draw [line width=2.pt] (-8.801937735804838,4.513518789299707)-- (-6.,1.);
		\draw [line width=2.pt] (-8.801937735804838,4.513518789299707)-- (-4.753020396282533,2.563662964936059);
		\draw [line width=2.pt] (-8.801937735804838,4.513518789299707)-- (-8.,1.);
		\draw [line width=2.pt] (-7.,5.381286267534822)-- (-8.,1.);
		\draw [line width=2.pt] (-7.,5.381286267534822)-- (-9.246979603717467,2.56366296493606);
		\draw [line width=2.pt] (-5.198062264195162,4.513518789299706)-- (-9.246979603717467,2.56366296493606);
		\draw [line width=2.pt] (-4.753020396282533,2.563662964936059)-- (-9.246979603717467,2.56366296493606);
		\draw [line width=2.pt] (-6.,1.)-- (-9.246979603717467,2.56366296493606);
		\draw [line width=2.pt] (-1.8019377358048378,4.513518789299707)-- (1.8019377358048383,4.513518789299706);
		\draw [line width=2.pt] (-1.8019377358048378,4.513518789299707)-- (2.2469796037174667,2.563662964936059);
		\draw [line width=2.pt] (-1.8019377358048378,4.513518789299707)-- (1.,1.);
		\draw [line width=2.pt] (-1.8019377358048378,4.513518789299707)-- (-2.2469796037174667,2.56366296493606);
		\draw [line width=2.pt] (-2.2469796037174667,2.56366296493606)-- (0.,5.381286267534822);
		\draw [line width=2.pt] (0.,5.381286267534822)-- (-1.,1.);
		\draw [line width=2.pt] (-1.,1.)-- (1.8019377358048383,4.513518789299706);
		\draw [line width=2.pt] (2.2469796037174667,2.563662964936059)-- (-2.2469796037174667,2.56366296493606);
		\draw [line width=2.pt] (1.,1.)-- (-2.2469796037174667,2.56366296493606);
		\draw [line width=2.pt] (4.753020396282533,2.56366296493606)-- (5.198062264195162,4.513518789299707);
		\draw [line width=2.pt] (4.753020396282533,2.56366296493606)-- (6.,1.);
		\draw [line width=2.pt] (5.198062264195162,4.513518789299707)-- (8.801937735804838,4.513518789299706);
		\draw [line width=2.pt] (8.801937735804838,4.513518789299706)-- (6.,1.);
		\draw [line width=2.pt] (4.753020396282533,2.56366296493606)-- (9.246979603717467,2.563662964936059);
		\draw [line width=2.pt] (9.246979603717467,2.563662964936059)-- (5.198062264195162,4.513518789299707);
		\draw [line width=2.pt] (4.753020396282533,2.56366296493606)-- (7.,5.381286267534822);
		\draw [line width=2.pt] (7.,5.381286267534822)-- (6.,1.);
		\draw [line width=2.pt] (6.,1.)-- (5.198062264195162,4.513518789299707);
		\draw [line width=2.pt] (-8.801937735804838,-1.4864812107002932)-- (-7.,-0.6187137324651775);
		\draw [line width=2.pt] (-8.801937735804838,-1.4864812107002932)-- (-5.198062264195162,-1.4864812107002938);
		\draw [line width=2.pt] (-8.801937735804838,-1.4864812107002932)-- (-4.753020396282533,-3.436337035063941);
		\draw [line width=2.pt] (-8.801937735804838,-1.4864812107002932)-- (-9.246979603717467,-3.43633703506394);
		\draw [line width=2.pt] (-8.801937735804838,-1.4864812107002932)-- (-6.,-5.);
		\draw [line width=2.pt] (-8.801937735804838,-1.4864812107002932)-- (-8.,-5.);
		\draw [line width=2.pt] (-8.,-5.)-- (-5.198062264195162,-1.4864812107002938);
		\draw [line width=2.pt] (-6.,-5.)-- (-7.,-0.6187137324651775);
		\draw [line width=2.pt] (-9.246979603717467,-3.43633703506394)-- (-4.753020396282533,-3.436337035063941);
		\draw [line width=2.pt] (-1.8019377358048378,-1.4864812107002932)-- (-2.2469796037174667,-3.43633703506394);
		\draw [line width=2.pt] (-1.8019377358048378,-1.4864812107002932)-- (2.2469796037174667,-3.436337035063941);
		\draw [line width=2.pt] (-2.2469796037174667,-3.43633703506394)-- (2.2469796037174667,-3.436337035063941);
		\draw [line width=2.pt] (-1.8019377358048378,-1.4864812107002932)-- (-1.,-5.);
		\draw [line width=2.pt] (-1.,-5.)-- (1.8019377358048383,-1.4864812107002938);
		\draw [line width=2.pt] (1.8019377358048383,-1.4864812107002938)-- (-2.2469796037174667,-3.43633703506394);
		\draw [line width=2.pt] (1.,-5.)-- (-1.8019377358048378,-1.4864812107002932);
		\draw [line width=2.pt] (0.,-0.6187137324651775)-- (1.,-5.);
		\draw [line width=2.pt] (-2.2469796037174667,-3.43633703506394)-- (0.,-0.6187137324651775);
		\draw [line width=2.pt] (5.198062264195162,-1.4864812107002932)-- (4.753020396282533,-3.43633703506394);
		\draw [line width=2.pt] (4.753020396282533,-3.43633703506394)-- (8.,-5.);
		\draw [line width=2.pt] (5.198062264195162,-1.4864812107002932)-- (6.,-5.);
		\draw [line width=2.pt] (6.,-5.)-- (8.801937735804838,-1.4864812107002938);
		\draw [line width=2.pt] (8.801937735804838,-1.4864812107002938)-- (5.198062264195162,-1.4864812107002932);
		\draw [line width=2.pt] (8.,-5.)-- (7.,-0.6187137324651775);
		\draw [line width=2.pt] (7.,-0.6187137324651775)-- (4.753020396282533,-3.43633703506394);
		\draw [line width=2.pt] (5.198062264195162,-1.4864812107002932)-- (9.246979603717467,-3.436337035063941);
		\draw [line width=2.pt] (9.246979603717467,-3.436337035063941)-- (4.753020396282533,-3.43633703506394);
		\begin{scriptsize}
			\draw [fill=uuuuuu] (-8.,1.) circle (2.0pt);
			\draw[color=uuuuuu] (-8.22316292125253,0.7810334052207423) node {$4$};
			\draw [fill=uuuuuu] (-6.,1.) circle (2.0pt);
			\draw[color=uuuuuu] (-5.7629113015872385,0.8365595487978362) node {$3$};
			\draw [fill=uuuuuu] (-4.753020396282533,2.563662964936059) circle (2.0pt);
			\draw[color=uuuuuu] (-4.476207061023311,2.5951036558639826) node {$2$};
			\draw [fill=uuuuuu] (-5.198062264195162,4.513518789299706) circle (2.0pt);
			\draw[color=uuuuuu] (-4.975809999314079,4.637856911546881) node {$1$};
			\draw [fill=uuuuuu] (-7.,5.381286267534822) circle (2.0pt);
			\draw[color=uuuuuu] (-6.878695608900055,5.703641218859697) node {$0$};
			\draw[color=red] (-5.878695608900055,5.703641218859697) node {$G_1$};
			\draw [fill=uuuuuu] (-8.801937735804838,4.513518789299707) circle (2.0pt);
			\draw[color=uuuuuu] (-9.081316510679876,4.673383055123974) node {$6$};
			\draw [fill=uuuuuu] (-9.246979603717467,2.56366296493606) circle (2.0pt);
			\draw[color=uuuuuu] (-9.57868252075919,2.4352360097670602) node {$5$};
			\draw [fill=uuuuuu] (-1.,1.) circle (2.0pt);
			\draw[color=uuuuuu] (-1.2478018519306748,0.7655072616436485) node {$4$};
			\draw [fill=uuuuuu] (1.,1.) circle (2.0pt);
			\draw[color=uuuuuu] (1.2567912702544437,0.9253749077405709) node {$3$};
			\draw [fill=uuuuuu] (2.2469796037174667,2.563662964936059) circle (2.0pt);
			\draw[color=uuuuuu] (2.5002062954527298,2.541814440498342) node {$2$};
			\draw [fill=uuuuuu] (1.8019377358048383,4.513518789299706) circle (2.0pt);
			\draw[color=uuuuuu] (1.9673141417963214,4.797724557643803) node {$1$};
			\draw [fill=uuuuuu] (0.,5.381286267534822) circle (2.0pt);
			\draw[color=uuuuuu] (0.15548081936453348,5.756930434225337) node {$0$};
			\draw[color=red] (1.15548081936453348,5.756930434225337) node {$G_2$};
			\draw [fill=uuuuuu] (-1.8019377358048378,4.513518789299707) circle (2.0pt);
			\draw[color=uuuuuu] (-1.851746292741271,4.904302988375084) node {$6$};
			\draw [fill=uuuuuu] (-2.2469796037174667,2.56366296493606) circle (2.0pt);
			\draw[color=uuuuuu] (-2.6155583796487893,2.470762153344154) node {$5$};
			\draw [fill=uuuuuu] (6.,1.) circle (2.0pt);
			\draw[color=uuuuuu] (5.8096637916995555,0.7710334052207423) node {$4$};
			\draw [fill=uuuuuu] (8.,1.) circle (2.0pt);
			\draw[color=uuuuuu] (8.196626195999205,0.7887964770092893) node {$3$};
			\draw [fill=uuuuuu] (9.246979603717467,2.563662964936059) circle (2.0pt);
			\draw[color=uuuuuu] (9.44004122119749,2.387472937978513) node {$2$};
			\draw [fill=uuuuuu] (8.801937735804838,4.513518789299706) circle (2.0pt);
			\draw[color=uuuuuu] (9.019253641849458,4.81548762943235) node {$1$};
			\draw [fill=uuuuuu] (7.,5.381286267534822) circle (2.0pt);
			\draw[color=uuuuuu] (7.118604960474935,5.703641218859697) node {$0$};
			\draw[color=blue] (8.118604960474935,5.703641218859697) node {$G_3$};
			\draw [fill=uuuuuu] (5.198062264195162,4.513518789299707) circle (2.0pt);
			\draw[color=uuuuuu] (4.9337471304836615,4.797724557643803) node {$6$};
			\draw [fill=uuuuuu] (4.753020396282533,2.56366296493606) circle (2.0pt);
			\draw[color=uuuuuu] (4.489670335769987,2.293131435458685) node {$5$};
			\draw [fill=uuuuuu] (-8.,-5.) circle (2.0pt);
			\draw[color=uuuuuu] (-8.317504423772357,-5.131832572153933) node {$4$};
			\draw [fill=uuuuuu] (-6.,-5.) circle (2.0pt);
			\draw[color=uuuuuu] (-5.765148229798691,-5.131832572153933) node {$3$};
			\draw [fill=uuuuuu] (-4.753020396282533,-3.436337035063941) circle (2.0pt);
			\draw[color=uuuuuu] (-4.509496276388953,-3.2311838907794117) node {$2$};
			\draw [fill=uuuuuu] (-5.198062264195162,-1.4864812107002938) circle (2.0pt);
			\draw[color=uuuuuu] (-4.986928463419806,-1.161786027413694) node {$1$};
			\draw [fill=uuuuuu] (-7.,-0.6187137324651775) circle (2.0pt);
			\draw[color=uuuuuu] (-6.8342879294286885,-0.23810629440925335) node {$0$};
			\draw[color=blue] (-5.8342879294286885,-0.23810629440925335) node {$G_4$};
			\draw [fill=uuuuuu] (-8.801937735804838,-1.4864812107002932) circle (2.0pt);
			\draw[color=uuuuuu] (-9.125724190151242,-1.1440229556251469) node {$6$};
			\draw [fill=uuuuuu] (-9.246979603717467,-3.43633703506394) circle (2.0pt);
			\draw[color=uuuuuu] (-9.480985625922182,-3.630853006021718) node {$5$};
			\draw [fill=uuuuuu] (-1.,-5.) circle (2.0pt);
			\draw[color=uuuuuu] (-1.2311572442478544,-5.122951036259661) node {$4$};
			\draw [fill=uuuuuu] (1.,-5.) circle (2.0pt);
			\draw[color=uuuuuu] (1.2501466625716235,-5.0874248926825665) node {$3$};
			\draw [fill=uuuuuu] (2.2469796037174667,-3.436337035063941) circle (2.0pt);
			\draw[color=uuuuuu] (2.562377046712644,-3.293354642039326) node {$2$};
			\draw [fill=uuuuuu] (1.8019377358048383,-1.4864812107002938) circle (2.0pt);
			\draw[color=uuuuuu] (1.9761956776905951,-1.1084968120480532) node {$1$};
			\draw [fill=uuuuuu] (0.,-0.6187137324651775) circle (2.0pt);
			\draw[color=uuuuuu] (0.16436235525880694,-0.23810629440925335) node {$0$};
			\draw[color=blue] (1.16436235525880694,-0.23810629440925335) node {$G_5$};
			\draw [fill=uuuuuu] (-1.8019377358048378,-1.4864812107002932) circle (2.0pt);
			\draw[color=uuuuuu] (-1.629707895384434,-1.1084968120480532) node {$6$};
			\draw [fill=uuuuuu] (-2.2469796037174667,-3.43633703506394) circle (2.0pt);
			\draw[color=uuuuuu] (-2.5356245566003284,-3.5242745752904363) node {$5$};
			\draw [fill=uuuuuu] (6.,-5.) circle (2.0pt);
			\draw[color=uuuuuu] (5.741966896862547,-5.105187964471114) node {$4$};
			\draw [fill=uuuuuu] (8.,-5.) circle (2.0pt);
			\draw[color=uuuuuu] (8.207744660104932,-5.158477179836754) node {$3$};
			\draw [fill=uuuuuu] (9.246979603717467,-3.436337035063941) circle (2.0pt);
			\draw[color=uuuuuu] (9.525501187823046,-3.5065115035018897) node {$2$};
			\draw [fill=uuuuuu] (8.801937735804838,-1.4864812107002938) circle (2.0pt);
			\draw[color=uuuuuu] (8.97484596237809,-1.1084968120480532) node {$1$};
			\draw [fill=uuuuuu] (7.,-0.6187137324651775) circle (2.0pt);
			\draw[color=uuuuuu] (7.163012639946301,-0.23810629440925335) node {$0$};
			\draw[color=red] (8.163012639946301,-0.23810629440925335) node {$G_6$};
			\draw [fill=uuuuuu] (5.198062264195162,-1.4864812107002932) circle (2.0pt);
			\draw[color=uuuuuu] (5.368942389303061,-1.1084968120480532) node {$6$};
			\draw [fill=uuuuuu] (4.753020396282533,-3.43633703506394) circle (2.0pt);
			\draw[color=uuuuuu] (4.480788799875714,-3.6130899342331713) node {$5$};
		\end{scriptsize}
	\end{tikzpicture}
	\caption{6 even graphs of order 7 and size 9}
	\label{fig:nauty-7-9-e}
\end{figure}

\begin{figure}
	\centering
	\definecolor{ccwwqq}{rgb}{0.8,0.4,0.}
	\definecolor{xdxdff}{rgb}{0.49019607843137253,0.49019607843137253,1.}
	\definecolor{ududff}{rgb}{0.30196078431372547,0.30196078431372547,1.}
	\definecolor{wqwqwq}{rgb}{0.3764705882352941,0.3764705882352941,0.3764705882352941}
	\definecolor{uuuuuu}{rgb}{0.26666666666666666,0.26666666666666666,0.26666666666666666}
	\begin{tikzpicture}[line cap=round,line join=round,>=triangle 45,x=1.0cm,y=1.0cm]
		\clip(-3.,-3.) rectangle (3.,3.);
		\fill[line width=2.pt,color=wqwqwq,fill=wqwqwq,fill opacity=0.10000000149011612] (0.,2.3812862675348225) -- (1.8019377358048383,1.5135187892997062) -- (-1.8019377358048378,1.5135187892997068) -- cycle;
		\fill[line width=2.pt,color=wqwqwq,fill=wqwqwq,fill opacity=0.10000000149011612] (1.,-2.) -- (-1.,-2.) -- (-2.2469796037174667,-0.43633703506393995) -- cycle;
		\fill[line width=2.pt,color=wqwqwq,fill=wqwqwq,fill opacity=0.10000000149011612] (1.8019377358048383,1.5135187892997062) -- (2.2469796037174667,-0.43633703506394084) -- (1.,-2.) -- cycle;
		\fill[line width=2.pt,color=wqwqwq,fill=wqwqwq,fill opacity=0.10000000149011612] (0.,2.3812862675348225) -- (2.2469796037174667,-0.43633703506394084) -- (-1.,-2.) -- cycle;
		\fill[line width=2.pt,color=ududff,fill=ududff,fill opacity=0.10000000149011612] (0.,2.3812862675348225) -- (1.,-2.) -- (2.0403248844414823,-1.4866531267173806) -- cycle;
		\fill[line width=2.pt,color=ududff,fill=ududff,fill opacity=0.10000000149011612] (0.,2.3812862675348225) -- (-2.2469796037174667,-0.43633703506393995) -- (-2.520584158655396,0.5630411158648467) -- cycle;
		\fill[line width=2.pt,color=ududff,fill=ududff,fill opacity=0.10000000149011612] (1.8019377358048383,1.5135187892997062) -- (-2.2469796037174667,-0.43633703506393995) -- (-2.2762497456323487,-1.364485920205857) -- cycle;
		\fill[line width=2.pt,color=ududff,fill=ududff,fill opacity=0.10000000149011612] (1.8019377358048383,1.5135187892997062) -- (-1.,-2.) -- (0.,-2.3961201085253885) -- cycle;
		\fill[line width=2.pt,color=ududff,fill=ududff,fill opacity=0.10000000149011612] (2.2469796037174667,-0.43633703506394084) -- (-2.2469796037174667,-0.43633703506393995) -- (2.5018454423739045,0.5494669818080108) -- cycle;
		\fill[line width=2.pt,color=ududff,fill=ududff,fill opacity=0.10000000149011612] (2.2469796037174667,-0.43633703506394084) -- (-1.8019377358048378,1.5135187892997068) -- (-0.986707010232934,2.3412526773103552) -- cycle;
		\fill[line width=2.pt,color=ududff,fill=ududff,fill opacity=0.10000000149011612] (1.,-2.) -- (-1.8019377358048378,1.5135187892997068) -- (-2.3984169521438723,1.1467288803087923) -- cycle;
		\fill[line width=2.pt,color=ududff,fill=ududff,fill opacity=0.10000000149011612] (-1.,-2.) -- (-1.8019377358048378,1.5135187892997068) -- (-1.8283033217567626,-1.7309875397404275) -- cycle;
		\fill[line width=2.pt,color=ccwwqq,fill=ccwwqq,fill opacity=0.10000000149011612] (-1.9873728985517802,2.810593478276708) -- (-2.6165501296618636,1.9618082787913598) -- (-1.5713418072854959,1.9346600106776881) -- cycle;
		\draw [line width=2.pt,color=wqwqwq] (0.,2.3812862675348225)-- (1.8019377358048383,1.5135187892997062);
		\draw [line width=2.pt,color=wqwqwq] (1.8019377358048383,1.5135187892997062)-- (-1.8019377358048378,1.5135187892997068);
		\draw [line width=2.pt,color=wqwqwq] (-1.8019377358048378,1.5135187892997068)-- (0.,2.3812862675348225);
		\draw [line width=2.pt,color=wqwqwq] (1.,-2.)-- (-1.,-2.);
		\draw [line width=2.pt,color=wqwqwq] (-1.,-2.)-- (-2.2469796037174667,-0.43633703506393995);
		\draw [line width=2.pt,color=wqwqwq] (-2.2469796037174667,-0.43633703506393995)-- (1.,-2.);
		\draw [line width=2.pt,color=wqwqwq] (1.8019377358048383,1.5135187892997062)-- (2.2469796037174667,-0.43633703506394084);
		\draw [line width=2.pt,color=wqwqwq] (2.2469796037174667,-0.43633703506394084)-- (1.,-2.);
		\draw [line width=2.pt,color=wqwqwq] (1.,-2.)-- (1.8019377358048383,1.5135187892997062);
		\draw [line width=2.pt,color=wqwqwq] (0.,2.3812862675348225)-- (2.2469796037174667,-0.43633703506394084);
		\draw [line width=2.pt,color=wqwqwq] (2.2469796037174667,-0.43633703506394084)-- (-1.,-2.);
		\draw [line width=2.pt,color=wqwqwq] (-1.,-2.)-- (0.,2.3812862675348225);
		\draw [line width=2.pt,color=ududff] (0.,2.3812862675348225)-- (1.,-2.);
		\draw [line width=2.pt,color=ududff] (1.,-2.)-- (2.0403248844414823,-1.4866531267173806);
		\draw [line width=2.pt,color=ududff] (2.0403248844414823,-1.4866531267173806)-- (0.,2.3812862675348225);
		\draw [line width=2.pt,color=ududff] (0.,2.3812862675348225)-- (-2.2469796037174667,-0.43633703506393995);
		\draw [line width=2.pt,color=ududff] (-2.2469796037174667,-0.43633703506393995)-- (-2.520584158655396,0.5630411158648467);
		\draw [line width=2.pt,color=ududff] (-2.520584158655396,0.5630411158648467)-- (0.,2.3812862675348225);
		\draw [line width=2.pt,color=ududff] (1.8019377358048383,1.5135187892997062)-- (-2.2469796037174667,-0.43633703506393995);
		\draw [line width=2.pt,color=ududff] (-2.2469796037174667,-0.43633703506393995)-- (-2.2762497456323487,-1.364485920205857);
		\draw [line width=2.pt,color=ududff] (-2.2762497456323487,-1.364485920205857)-- (1.8019377358048383,1.5135187892997062);
		\draw [line width=2.pt,color=ududff] (1.8019377358048383,1.5135187892997062)-- (-1.,-2.);
		\draw [line width=2.pt,color=ududff] (-1.,-2.)-- (0.,-2.3961201085253885);
		\draw [line width=2.pt,color=ududff] (0.,-2.3961201085253885)-- (1.8019377358048383,1.5135187892997062);
		\draw [line width=2.pt,color=ududff] (2.2469796037174667,-0.43633703506394084)-- (-2.2469796037174667,-0.43633703506393995);
		\draw [line width=2.pt,color=ududff] (-2.2469796037174667,-0.43633703506393995)-- (2.5018454423739045,0.5494669818080108);
		\draw [line width=2.pt,color=ududff] (2.5018454423739045,0.5494669818080108)-- (2.2469796037174667,-0.43633703506394084);
		\draw [line width=2.pt,color=ududff] (2.2469796037174667,-0.43633703506394084)-- (-1.8019377358048378,1.5135187892997068);
		\draw [line width=2.pt,color=ududff] (-1.8019377358048378,1.5135187892997068)-- (-0.986707010232934,2.3412526773103552);
		\draw [line width=2.pt,color=ududff] (-0.986707010232934,2.3412526773103552)-- (2.2469796037174667,-0.43633703506394084);
		\draw [line width=2.pt,color=ududff] (1.,-2.)-- (-1.8019377358048378,1.5135187892997068);
		\draw [line width=2.pt,color=ududff] (-1.8019377358048378,1.5135187892997068)-- (-2.3984169521438723,1.1467288803087923);
		\draw [line width=2.pt,color=ududff] (-2.3984169521438723,1.1467288803087923)-- (1.,-2.);
		\draw [line width=2.pt,color=ududff] (-1.,-2.)-- (-1.8019377358048378,1.5135187892997068);
		\draw [line width=2.pt,color=ududff] (-1.8019377358048378,1.5135187892997068)-- (-1.8283033217567626,-1.7309875397404275);
		\draw [line width=2.pt,color=ududff] (-1.8283033217567626,-1.7309875397404275)-- (-1.,-2.);
		\draw [line width=2.pt,color=ccwwqq] (-1.9873728985517802,2.810593478276708)-- (-2.6165501296618636,1.9618082787913598);
		\draw [line width=2.pt,color=ccwwqq] (-2.6165501296618636,1.9618082787913598)-- (-1.5713418072854959,1.9346600106776881);
		\draw [line width=2.pt,color=ccwwqq] (-1.5713418072854959,1.9346600106776881)-- (-1.9873728985517802,2.810593478276708);
		\begin{scriptsize}
			\draw [fill=uuuuuu] (-1.,-2.) circle (2.0pt);
			\draw[color=uuuuuu] (-1.1464082757992877,-2.349228707459076) node {$w_5$};
			\draw [fill=uuuuuu] (1.,-2.) circle (2.0pt);
			\draw[color=uuuuuu] (1.1012157819837811,-2.349228707459076) node {$w_4$};
			\draw [fill=uuuuuu] (2.2469796037174667,-0.43633703506394084) circle (2.0pt);
			\draw[color=uuuuuu] (2.45392592680573,-0.558213924983421) node {$w_3$};
			\draw [fill=uuuuuu] (1.8019377358048383,1.5135187892997062) circle (2.0pt);
			\draw[color=uuuuuu] (1.8840960717733892,1.7525456400407449) node {$w_2$};
			\draw [fill=uuuuuu] (0.,2.3812862675348225) circle (2.0pt);
			\draw[color=uuuuuu] (0.09104766612622209,2.61118853851967) node {$w_1$};
			\draw [fill=uuuuuu] (-1.8019377358048378,1.5135187892997068) circle (2.0pt);
			\draw[color=uuuuuu] (-1.9671698699335542,1.7146643356960865) node {$w_7$};
			\draw [fill=uuuuuu] (-2.2469796037174667,-0.43633703506393995) circle (2.0pt);
			\draw[color=uuuuuu] (-2.5101352322069923,-0.4698242148458845) node {$w_6$};
			\draw [fill=ududff] (2.0403248844414823,-1.4866531267173806) circle (1pt);
			\draw [fill=ududff] (-2.520584158655396,0.5630411158648467) circle (1pt);
			\draw [fill=ududff] (-2.2762497456323487,-1.364485920205857) circle (1pt);
			\draw [fill=xdxdff] (0.,-2.3961201085253885) circle (1pt);
			\draw [fill=ududff] (2.5018454423739045,0.5494669818080108) circle (1pt);
			\draw [fill=ududff] (-0.986707010232934,2.3412526773103552) circle (1pt);
			\draw [fill=ududff] (-2.3984169521438723,1.1467288803087923) circle (1pt);
			\draw [fill=ududff] (-1.8283033217567626,-1.7309875397404275) circle (1pt);
			\draw [fill=ududff] (-1.9873728985517802,2.810593478276708) circle (1pt);
			\draw [fill=ududff] (-2.6165501296618636,1.9618082787913598) circle (1pt);
			\draw [fill=ududff] (-1.5713418072854959,1.9346600106776881) circle (1pt);
		\end{scriptsize}
	\end{tikzpicture}
	\caption{13 3-cliques in $K_{18}$ for case \eqref{eq:c4}}
	\label{fig:k18alpha13-7}
\end{figure}

\begin{lemma}\label{lem:k18a13c5-1}
	In a $\{K_3,K_4\}$-decomposition of $K_{18}$, if $\alpha=13$ with the case \eqref{eq:c5} holds, then up to isomorphism $H\notin \{G_3,G_4\}$.
\end{lemma}
\begin{proof}
For contradiction, suppose case \eqref{eq:c5} is possible. Since $t_3=4$, as in case \eqref{eq:c4}, it remains to consider $H\in\{G_3,G_4,G_5\}$ from Figure \ref{fig:nauty-7-9-e}. The complement $\overline H$ in the $K_7$ on $w_1,\ldots,w_7$ gives the union of the four 3-cliques lying entirely in this vertex set. 
	
	If $H\cong G_3$, we first take the 4 3-cliques in $\overline{H}$: $\{w_1,w_2,w_3\}$, $\{w_1,w_4,w_7\}$, $\{w_2,w_4,w_6\}$, $\{w_3,w_4,w_5\}$, and we form another 3-clique for $w_4$: $\{w_4,y_1,y_2\}$. Notice that, $w_5,w_6,w_7$ each have 4 incident edges in $H$ and are each contained in 1 3-clique in $\overline{H}$. Since $t_2=7$ in \eqref{eq:c5}, we deduce that precisely 2 edges among $\{w_5,w_6\}$, $\{w_5,w_7\}$, and $\{w_6,w_7\}$ must not form a 3-clique. Without loss of generality, we assume that all edges of $H$ except $\{w_5,w_6\}$ and $\{w_6,w_7\}$ contribute to 3-cliques, which may be taken as: $\{w_1,w_5, y_3\}$, $\{w_1,w_6,y_4\}$, $\{w_2,w_5,y_5\}$, $\{w_2,w_7,y_6\}$, $\{w_3,w_6,y_7\}$, $\{w_3,w_7,y_8\}$, $\{w_5,w_7,y_9\}$. We finally append a 3-clique for $w_6$: $\{w_6,y_{10}, y_{11}\}$. This configuration (See Figure \ref{fig:k18alpha13-8}) is then shown infeasible by an exhaustive computer-assisted search. 
	
	If $H\cong G_4$, we first take the 4 3-cliques in $\overline{H}$: $\{w_1,w_2,w_3\}$, $\{w_1,w_5,w_6\}$, $\{w_2,w_4,w_6\}$ and, $\{w_3,w_4,w_5\}$. Since $\deg_{H}(w_7)=6$, we deduce that precisely 2 out of the 6 edges incident to $w_7$ must not contribute to 3-cliques. Up to isomorphism, there are two cases. 
	\begin{enumerate}
		\item These 2 edges are chosen in the same 3-clique in $H$. Without loss of generality, we assume that all edges except $\{w_3,w_7\}$ and $\{w_6,w_7\}$ contribute to 3-cliques, which we take as: $\{w_1,w_4,y_1\}$, $\{w_1,w_7,y_2\}$, $\{w_2,w_5,y_3\}$, $\{w_2,w_7,y_4\}$, $\{w_3,w_6,y_5\}$, $\{w_4,w_7,y_6\}$, $\{w_5,w_7,y_7\}$. We finally append two 3-cliques for $w_3$ and $w_6$ respectively: $\{w_3,y_8,y_9\}$, $\{w_6,y_{10}, y_{11}\}$. This configuration (See Figure \ref{fig:k18alpha13-9} left) is then shown infeasible by an exhaustive computer-assisted search. 
		\item These 2 edges are chosen from different 3-cliques in $H$. Without loss of generality, we assume that all edges except $\{w_1,w_7\}$ and $\{w_2,w_7\}$ contribute to 3-cliques, which we take as: $\{w_1,w_4,y_1\}$, $\{w_2,w_5,y_2\}$, $\{w_3,w_6,y_3\}$, $\{w_3,w_7,y_4\}$, $\{w_4,w_7,y_5\}$, $\{w_5,w_7,y_6\}$, $\{w_6,w_7,y_7\}$. We finally append two 3-cliques for $w_1$ and $w_2$ respectively: $\{w_1,y_8,y_9\}$, $\{w_2,y_{10}, y_{11}\}$. This configuration (See Figure \ref{fig:k18alpha13-9} right) is then shown infeasible by an exhaustive computer-assisted search. 
	\end{enumerate} 
This concludes the proof.
\end{proof}

\begin{figure}
	\centering
	\definecolor{wwzzqq}{rgb}{0.4,0.6,0.}
	\definecolor{ududff}{rgb}{0.30196078431372547,0.30196078431372547,1.}
	\definecolor{wqwqwq}{rgb}{0.3764705882352941,0.3764705882352941,0.3764705882352941}
	\definecolor{uuuuuu}{rgb}{0.26666666666666666,0.26666666666666666,0.26666666666666666}
	\begin{tikzpicture}[line cap=round,line join=round,>=triangle 45,x=1.0cm,y=1.0cm]
		\clip(-3.,-3.) rectangle (3.,3.);
		\fill[line width=2.pt,color=wqwqwq,fill=wqwqwq,fill opacity=0.10000000149011612] (0.,2.3812862675348225) -- (1.8019377358048383,1.5135187892997062) -- (2.2469796037174667,-0.43633703506394084) -- cycle;
		\fill[line width=2.pt,color=wqwqwq,fill=wqwqwq,fill opacity=0.10000000149011612] (0.,2.3812862675348225) -- (1.,-2.) -- (-1.8019377358048378,1.5135187892997068) -- cycle;
		\fill[line width=2.pt,color=wqwqwq,fill=wqwqwq,fill opacity=0.10000000149011612] (1.8019377358048383,1.5135187892997062) -- (1.,-2.) -- (-2.2469796037174667,-0.43633703506393995) -- cycle;
		\fill[line width=2.pt,color=wqwqwq,fill=wqwqwq,fill opacity=0.10000000149011612] (2.2469796037174667,-0.43633703506394084) -- (1.,-2.) -- (-1.,-2.) -- cycle;
		\fill[line width=2.pt,color=wwzzqq,fill=wwzzqq,fill opacity=0.10000000149011612] (1.,-2.) -- (1.9347714866423569,-2.3291730838033864) -- (0.5381933877484443,-2.7008430617348305) -- cycle;
		\fill[line width=2.pt,color=wwzzqq,fill=wwzzqq,fill opacity=0.10000000149011612] (-2.2469796037174667,-0.43633703506393995) -- (-2.682946420990741,-1.2254258766130377) -- (-2.727997327406674,0.08105040944900749) -- cycle;
		\fill[line width=2.pt,color=ududff,fill=ududff,fill opacity=0.10000000149011612] (0.,2.3812862675348225) -- (-1.,-2.) -- (-1.877335068956929,-1.7265066747873177) -- cycle;
		\fill[line width=2.pt,color=ududff,fill=ududff,fill opacity=0.10000000149011612] (0.,2.3812862675348225) -- (-2.2469796037174667,-0.43633703506393995) -- (-2.3692438929721042,0.5438417437442564) -- cycle;
		\fill[line width=2.pt,color=ududff,fill=ududff,fill opacity=0.10000000149011612] (1.8019377358048383,1.5135187892997062) -- (-1.,-2.) -- (-0.3259303162936841,-2.398151415269575) -- cycle;
		\fill[line width=2.pt,color=ududff,fill=ududff,fill opacity=0.10000000149011612] (1.8019377358048383,1.5135187892997062) -- (-1.8019377358048378,1.5135187892997068) -- (1.1577256148632524,2.2134462241542714) -- cycle;
		\fill[line width=2.pt,color=ududff,fill=ududff,fill opacity=0.10000000149011612] (2.2469796037174667,-0.43633703506394084) -- (-1.8019377358048378,1.5135187892997068) -- (2.342023421823769,0.717117857276017) -- cycle;
		\fill[line width=2.pt,color=ududff,fill=ududff,fill opacity=0.10000000149011612] (2.2469796037174667,-0.43633703506394084) -- (-2.2469796037174667,-0.43633703506393995) -- (1.9378020146576045,-1.410363233072212) -- cycle;
		\fill[line width=2.pt,color=ududff,fill=ududff,fill opacity=0.10000000149011612] (-1.,-2.) -- (-1.8019377358048378,1.5135187892997068) -- (-2.,-1.) -- cycle;
		\draw [line width=2.pt,color=wqwqwq] (0.,2.3812862675348225)-- (1.8019377358048383,1.5135187892997062);
		\draw [line width=2.pt,color=wqwqwq] (1.8019377358048383,1.5135187892997062)-- (2.2469796037174667,-0.43633703506394084);
		\draw [line width=2.pt,color=wqwqwq] (2.2469796037174667,-0.43633703506394084)-- (0.,2.3812862675348225);
		\draw [line width=2.pt,color=wqwqwq] (0.,2.3812862675348225)-- (1.,-2.);
		\draw [line width=2.pt,color=wqwqwq] (1.,-2.)-- (-1.8019377358048378,1.5135187892997068);
		\draw [line width=2.pt,color=wqwqwq] (-1.8019377358048378,1.5135187892997068)-- (0.,2.3812862675348225);
		\draw [line width=2.pt,color=wqwqwq] (1.8019377358048383,1.5135187892997062)-- (1.,-2.);
		\draw [line width=2.pt,color=wqwqwq] (1.,-2.)-- (-2.2469796037174667,-0.43633703506393995);
		\draw [line width=2.pt,color=wqwqwq] (-2.2469796037174667,-0.43633703506393995)-- (1.8019377358048383,1.5135187892997062);
		\draw [line width=2.pt,color=wqwqwq] (2.2469796037174667,-0.43633703506394084)-- (1.,-2.);
		\draw [line width=2.pt,color=wqwqwq] (1.,-2.)-- (-1.,-2.);
		\draw [line width=2.pt,color=wqwqwq] (-1.,-2.)-- (2.2469796037174667,-0.43633703506394084);
		\draw [line width=2.pt,color=wwzzqq] (1.,-2.)-- (1.9347714866423569,-2.3291730838033864);
		\draw [line width=2.pt,color=wwzzqq] (1.9347714866423569,-2.3291730838033864)-- (0.5381933877484443,-2.7008430617348305);
		\draw [line width=2.pt,color=wwzzqq] (0.5381933877484443,-2.7008430617348305)-- (1.,-2.);
		\draw [line width=2.pt,color=wwzzqq] (-2.2469796037174667,-0.43633703506393995)-- (-2.682946420990741,-1.2254258766130377);
		\draw [line width=2.pt,color=wwzzqq] (-2.682946420990741,-1.2254258766130377)-- (-2.727997327406674,0.08105040944900749);
		\draw [line width=2.pt,color=wwzzqq] (-2.727997327406674,0.08105040944900749)-- (-2.2469796037174667,-0.43633703506393995);
		\draw [line width=2.pt,color=ududff] (0.,2.3812862675348225)-- (-1.,-2.);
		\draw [line width=2.pt,color=ududff] (-1.,-2.)-- (-1.877335068956929,-1.7265066747873177);
		\draw [line width=2.pt,color=ududff] (-1.877335068956929,-1.7265066747873177)-- (0.,2.3812862675348225);
		\draw [line width=2.pt,color=ududff] (0.,2.3812862675348225)-- (-2.2469796037174667,-0.43633703506393995);
		\draw [line width=2.pt,color=ududff] (-2.2469796037174667,-0.43633703506393995)-- (-2.3692438929721042,0.5438417437442564);
		\draw [line width=2.pt,color=ududff] (-2.3692438929721042,0.5438417437442564)-- (0.,2.3812862675348225);
		\draw [line width=2.pt,color=ududff] (1.8019377358048383,1.5135187892997062)-- (-1.,-2.);
		\draw [line width=2.pt,color=ududff] (-1.,-2.)-- (-0.3259303162936841,-2.398151415269575);
		\draw [line width=2.pt,color=ududff] (-0.3259303162936841,-2.398151415269575)-- (1.8019377358048383,1.5135187892997062);
		\draw [line width=2.pt,color=ududff] (1.8019377358048383,1.5135187892997062)-- (-1.8019377358048378,1.5135187892997068);
		\draw [line width=2.pt,color=ududff] (-1.8019377358048378,1.5135187892997068)-- (1.1577256148632524,2.2134462241542714);
		\draw [line width=2.pt,color=ududff] (1.1577256148632524,2.2134462241542714)-- (1.8019377358048383,1.5135187892997062);
		\draw [line width=2.pt,color=ududff] (2.2469796037174667,-0.43633703506394084)-- (-1.8019377358048378,1.5135187892997068);
		\draw [line width=2.pt,color=ududff] (-1.8019377358048378,1.5135187892997068)-- (2.342023421823769,0.717117857276017);
		\draw [line width=2.pt,color=ududff] (2.342023421823769,0.717117857276017)-- (2.2469796037174667,-0.43633703506394084);
		\draw [line width=2.pt,color=ududff] (2.2469796037174667,-0.43633703506394084)-- (-2.2469796037174667,-0.43633703506393995);
		\draw [line width=2.pt,color=ududff] (-2.2469796037174667,-0.43633703506393995)-- (1.9378020146576045,-1.410363233072212);
		\draw [line width=2.pt,color=ududff] (1.9378020146576045,-1.410363233072212)-- (2.2469796037174667,-0.43633703506394084);
		\draw [line width=2.pt,color=ududff] (-1.,-2.)-- (-1.8019377358048378,1.5135187892997068);
		\draw [line width=2.pt,color=ududff] (-1.8019377358048378,1.5135187892997068)-- (-2.,-1.);
		\draw [line width=2.pt,color=ududff] (-2.,-1.)-- (-1.,-2.);
		\begin{scriptsize}
			\draw [fill=uuuuuu] (-1.,-2.) circle (2pt);
			\draw[color=uuuuuu] (-1.2236046021134125,-2.12588086860149) node {$w_5$};
			\draw [fill=uuuuuu] (1.,-2.) circle (2pt);
			\draw[color=uuuuuu] (1.37702825296974,-1.8904898825782417) node {$w_4$};
			\draw [fill=uuuuuu] (2.2469796037174667,-0.43633703506394084) circle (2pt);
			\draw[color=uuuuuu] (2.4939831830859815,-0.35425397379494056) node {$w_3$};
			\draw [fill=uuuuuu] (1.8019377358048383,1.5135187892997062) circle (2pt);
			\draw[color=uuuuuu] (1.8840892169244328,1.7394869018855266) node {$w_2$};
			\draw [fill=uuuuuu] (0.,2.3812862675348225) circle (2pt);
			\draw[color=uuuuuu] (0.08768432358911762,2.6067168503922287) node {$w_1$};
			\draw [fill=uuuuuu] (-1.8019377358048378,1.5135187892997068) circle (2pt);
			\draw[color=uuuuuu] (-2.0184455513557347,1.689930904828001) node {$w_7$};
			\draw [fill=uuuuuu] (-2.2469796037174667,-0.43633703506393995) circle (2pt);
			\draw[color=uuuuuu] (-2.178614535907746,-0.7135349524620029) node {$w_6$};
			\draw [fill=ududff] (1.9347714866423569,-2.3291730838033864) circle (1pt);
			\draw [fill=ududff] (0.5381933877484443,-2.7008430617348305) circle (1pt);
			\draw [fill=ududff] (-2.682946420990741,-1.2254258766130377) circle (1pt);
			\draw [fill=ududff] (-2.727997327406674,0.08105040944900749) circle (1pt);
			\draw [fill=ududff] (-1.877335068956929,-1.7265066747873177) circle (1pt);
			\draw [fill=ududff] (-2.3692438929721042,0.5438417437442564) circle (1pt);
			\draw [fill=ududff] (-0.3259303162936841,-2.398151415269575) circle (1pt);
			\draw [fill=ududff] (1.1577256148632524,2.2134462241542714) circle (1pt);
			\draw [fill=ududff] (2.342023421823769,0.717117857276017) circle (1pt);
			\draw [fill=ududff] (1.9378020146576045,-1.410363233072212) circle (1pt);
			\draw [fill=ududff] (-2.,-1.) circle (1pt);
		\end{scriptsize}
	\end{tikzpicture}
	\caption{13 3-cliques in $K_{18}$ for case \eqref{eq:c5}}
	\label{fig:k18alpha13-8}
\end{figure}

\begin{figure}
	\centering
	\definecolor{xdxdff}{rgb}{0.49019607843137253,0.49019607843137253,1.}
	\definecolor{wwzzqq}{rgb}{0.4,0.6,0.}
	\definecolor{qqwwzz}{rgb}{0.,0.4,0.6}
	\definecolor{ududff}{rgb}{0.30196078431372547,0.30196078431372547,1.}
	\definecolor{wqwqwq}{rgb}{0.3764705882352941,0.3764705882352941,0.3764705882352941}
	\definecolor{uuuuuu}{rgb}{0.26666666666666666,0.26666666666666666,0.26666666666666666}
	\begin{tikzpicture}[line cap=round,line join=round,>=triangle 45,x=1.0cm,y=1.0cm]
		\clip(-7.,-3.) rectangle (7.,3.5);
		\fill[line width=2.pt,color=wqwqwq,fill=wqwqwq,fill opacity=0.10000000149011612] (-4.,2.3812862675348225) -- (-2.198062264195162,1.5135187892997062) -- (-1.7530203962825333,-0.43633703506394084) -- cycle;
		\fill[line width=2.pt,color=wqwqwq,fill=wqwqwq,fill opacity=0.10000000149011612] (-4.,2.3812862675348225) -- (-5.,-2.) -- (-6.246979603717467,-0.43633703506393995) -- cycle;
		\fill[line width=2.pt,color=wqwqwq,fill=wqwqwq,fill opacity=0.10000000149011612] (-2.198062264195162,1.5135187892997062) -- (-3.,-2.) -- (-6.246979603717467,-0.43633703506393995) -- cycle;
		\fill[line width=2.pt,color=wqwqwq,fill=wqwqwq,fill opacity=0.10000000149011612] (-1.7530203962825333,-0.43633703506394084) -- (-3.,-2.) -- (-5.,-2.) -- cycle;
		\fill[line width=2.pt,color=wqwqwq,fill=wqwqwq,fill opacity=0.10000000149011612] (4.,2.3812862675348225) -- (5.801937735804838,1.5135187892997062) -- (6.246979603717467,-0.43633703506394084) -- cycle;
		\fill[line width=2.pt,color=wqwqwq,fill=wqwqwq,fill opacity=0.10000000149011612] (4.,2.3812862675348225) -- (3.,-2.) -- (1.7530203962825333,-0.43633703506393995) -- cycle;
		\fill[line width=2.pt,color=wqwqwq,fill=wqwqwq,fill opacity=0.10000000149011612] (5.801937735804838,1.5135187892997062) -- (5.,-2.) -- (1.7530203962825333,-0.43633703506393995) -- cycle;
		\fill[line width=2.pt,color=wqwqwq,fill=wqwqwq,fill opacity=0.10000000149011612] (6.246979603717467,-0.43633703506394084) -- (5.,-2.) -- (3.,-2.) -- cycle;
		\fill[line width=2.pt,color=qqwwzz,fill=qqwwzz,fill opacity=0.10000000149011612] (-4.,2.3812862675348225) -- (-3.,-2.) -- (-2.829565117941744,2.409429679766491) -- cycle;
		\fill[line width=2.pt,color=qqwwzz,fill=qqwwzz,fill opacity=0.10000000149011612] (-4.,2.3812862675348225) -- (-5.801937735804838,1.5135187892997068) -- (-5.147974069221822,2.3261335497803803) -- cycle;
		\fill[line width=2.pt,color=qqwwzz,fill=qqwwzz,fill opacity=0.10000000149011612] (-2.198062264195162,1.5135187892997062) -- (-5.,-2.) -- (-4.051241691071366,-2.3939804827658895) -- cycle;
		\fill[line width=2.pt,color=qqwwzz,fill=qqwwzz,fill opacity=0.10000000149011612] (-2.198062264195162,1.5135187892997062) -- (-5.801937735804838,1.5135187892997068) -- (-6.550125590654683,0.9378647166785362) -- cycle;
		\fill[line width=2.pt,color=qqwwzz,fill=qqwwzz,fill opacity=0.10000000149011612] (-1.7530203962825333,-0.43633703506394084) -- (-6.246979603717467,-0.43633703506393995) -- (-1.635653921474159,0.5491494434100199) -- cycle;
		\fill[line width=2.pt,color=qqwwzz,fill=qqwwzz,fill opacity=0.10000000149011612] (-3.,-2.) -- (-5.801937735804838,1.5135187892997068) -- (-2.0104865064116564,-1.366661546270525) -- cycle;
		\fill[line width=2.pt,color=qqwwzz,fill=qqwwzz,fill opacity=0.10000000149011612] (-5.,-2.) -- (-5.801937735804838,1.5135187892997068) -- (-5.967052680751909,-1.561019182904783) -- cycle;
		\fill[line width=2.pt,color=wwzzqq,fill=wwzzqq,fill opacity=0.10000000149011612] (-1.7530203962825333,-0.43633703506394084) -- (-1.,0.) -- (-1.,-1.) -- cycle;
		\fill[line width=2.pt,color=wwzzqq,fill=wwzzqq,fill opacity=0.10000000149011612] (-6.246979603717467,-0.43633703506393995) -- (-6.786131292281997,0.) -- (-6.786131292281997,-1.0195943379950638) -- cycle;
		\fill[line width=2.pt,color=qqwwzz,fill=qqwwzz,fill opacity=0.10000000149011612] (4.,2.3812862675348225) -- (5.,-2.) -- (5.139097984062835,2.3677816147734356) -- cycle;
		\fill[line width=2.pt,color=qqwwzz,fill=qqwwzz,fill opacity=0.10000000149011612] (5.801937735804838,1.5135187892997062) -- (3.,-2.) -- (6.416305310516531,0.7018590150512227) -- cycle;
		\fill[line width=2.pt,color=qqwwzz,fill=qqwwzz,fill opacity=0.10000000149011612] (6.246979603717467,-0.43633703506394084) -- (1.7530203962825333,-0.43633703506393995) -- (6.013707348916996,-1.4777230529186725) -- cycle;
		\fill[line width=2.pt,color=qqwwzz,fill=qqwwzz,fill opacity=0.10000000149011612] (6.246979603717467,-0.43633703506394084) -- (2.198062264195162,1.5135187892997068) -- (2.820689032782757,2.3122508614493618) -- cycle;
		\fill[line width=2.pt,color=qqwwzz,fill=qqwwzz,fill opacity=0.10000000149011612] (5.,-2.) -- (2.198062264195162,1.5135187892997068) -- (1.6128951479841538,0.7157417033822412) -- cycle;
		\fill[line width=2.pt,color=qqwwzz,fill=qqwwzz,fill opacity=0.10000000149011612] (3.,-2.) -- (2.198062264195162,1.5135187892997068) -- (4.097896359236453,-2.3662151061038528) -- cycle;
		\fill[line width=2.pt,color=qqwwzz,fill=qqwwzz,fill opacity=0.10000000149011612] (1.7530203962825333,-0.43633703506393995) -- (2.198062264195162,1.5135187892997068) -- (2.2098507462179464,-1.616549936228857) -- cycle;
		\fill[line width=2.pt,color=wwzzqq,fill=wwzzqq,fill opacity=0.10000000149011612] (4.,2.3812862675348225) -- (3.3898792543545126,2.8259103296970443) -- (4.4449635675119135,3.006385278000284) -- cycle;
		\fill[line width=2.pt,color=wwzzqq,fill=wwzzqq,fill opacity=0.10000000149011612] (5.801937735804838,1.5135187892997062) -- (6.055355413910052,2.3261335497803803) -- (6.666193700474863,1.7291779515465873) -- cycle;
		\draw [line width=2.pt,color=wqwqwq] (-4.,2.3812862675348225)-- (-2.198062264195162,1.5135187892997062);
		\draw [line width=2.pt,color=wqwqwq] (-2.198062264195162,1.5135187892997062)-- (-1.7530203962825333,-0.43633703506394084);
		\draw [line width=2.pt,color=wqwqwq] (-1.7530203962825333,-0.43633703506394084)-- (-4.,2.3812862675348225);
		\draw [line width=2.pt,color=wqwqwq] (-4.,2.3812862675348225)-- (-5.,-2.);
		\draw [line width=2.pt,color=wqwqwq] (-5.,-2.)-- (-6.246979603717467,-0.43633703506393995);
		\draw [line width=2.pt,color=wqwqwq] (-6.246979603717467,-0.43633703506393995)-- (-4.,2.3812862675348225);
		\draw [line width=2.pt,color=wqwqwq] (-2.198062264195162,1.5135187892997062)-- (-3.,-2.);
		\draw [line width=2.pt,color=wqwqwq] (-3.,-2.)-- (-6.246979603717467,-0.43633703506393995);
		\draw [line width=2.pt,color=wqwqwq] (-6.246979603717467,-0.43633703506393995)-- (-2.198062264195162,1.5135187892997062);
		\draw [line width=2.pt,color=wqwqwq] (-1.7530203962825333,-0.43633703506394084)-- (-3.,-2.);
		\draw [line width=2.pt,color=wqwqwq] (-3.,-2.)-- (-5.,-2.);
		\draw [line width=2.pt,color=wqwqwq] (-5.,-2.)-- (-1.7530203962825333,-0.43633703506394084);
		\draw [line width=2.pt,color=wqwqwq] (4.,2.3812862675348225)-- (5.801937735804838,1.5135187892997062);
		\draw [line width=2.pt,color=wqwqwq] (5.801937735804838,1.5135187892997062)-- (6.246979603717467,-0.43633703506394084);
		\draw [line width=2.pt,color=wqwqwq] (6.246979603717467,-0.43633703506394084)-- (4.,2.3812862675348225);
		\draw [line width=2.pt,color=wqwqwq] (4.,2.3812862675348225)-- (3.,-2.);
		\draw [line width=2.pt,color=wqwqwq] (3.,-2.)-- (1.7530203962825333,-0.43633703506393995);
		\draw [line width=2.pt,color=wqwqwq] (1.7530203962825333,-0.43633703506393995)-- (4.,2.3812862675348225);
		\draw [line width=2.pt,color=wqwqwq] (5.801937735804838,1.5135187892997062)-- (5.,-2.);
		\draw [line width=2.pt,color=wqwqwq] (5.,-2.)-- (1.7530203962825333,-0.43633703506393995);
		\draw [line width=2.pt,color=wqwqwq] (1.7530203962825333,-0.43633703506393995)-- (5.801937735804838,1.5135187892997062);
		\draw [line width=2.pt,color=wqwqwq] (6.246979603717467,-0.43633703506394084)-- (5.,-2.);
		\draw [line width=2.pt,color=wqwqwq] (5.,-2.)-- (3.,-2.);
		\draw [line width=2.pt,color=wqwqwq] (3.,-2.)-- (6.246979603717467,-0.43633703506394084);
		\draw [line width=2.pt,color=qqwwzz] (-4.,2.3812862675348225)-- (-3.,-2.);
		\draw [line width=2.pt,color=qqwwzz] (-3.,-2.)-- (-2.829565117941744,2.409429679766491);
		\draw [line width=2.pt,color=qqwwzz] (-2.829565117941744,2.409429679766491)-- (-4.,2.3812862675348225);
		\draw [line width=2.pt,color=qqwwzz] (-4.,2.3812862675348225)-- (-5.801937735804838,1.5135187892997068);
		\draw [line width=2.pt,color=qqwwzz] (-5.801937735804838,1.5135187892997068)-- (-5.147974069221822,2.3261335497803803);
		\draw [line width=2.pt,color=qqwwzz] (-5.147974069221822,2.3261335497803803)-- (-4.,2.3812862675348225);
		\draw [line width=2.pt,color=qqwwzz] (-2.198062264195162,1.5135187892997062)-- (-5.,-2.);
		\draw [line width=2.pt,color=qqwwzz] (-5.,-2.)-- (-4.051241691071366,-2.3939804827658895);
		\draw [line width=2.pt,color=qqwwzz] (-4.051241691071366,-2.3939804827658895)-- (-2.198062264195162,1.5135187892997062);
		\draw [line width=2.pt,color=qqwwzz] (-2.198062264195162,1.5135187892997062)-- (-5.801937735804838,1.5135187892997068);
		\draw [line width=2.pt,color=qqwwzz] (-5.801937735804838,1.5135187892997068)-- (-6.550125590654683,0.9378647166785362);
		\draw [line width=2.pt,color=qqwwzz] (-6.550125590654683,0.9378647166785362)-- (-2.198062264195162,1.5135187892997062);
		\draw [line width=2.pt,color=qqwwzz] (-1.7530203962825333,-0.43633703506394084)-- (-6.246979603717467,-0.43633703506393995);
		\draw [line width=2.pt,color=qqwwzz] (-6.246979603717467,-0.43633703506393995)-- (-1.635653921474159,0.5491494434100199);
		\draw [line width=2.pt,color=qqwwzz] (-1.635653921474159,0.5491494434100199)-- (-1.7530203962825333,-0.43633703506394084);
		\draw [line width=2.pt,color=qqwwzz] (-3.,-2.)-- (-5.801937735804838,1.5135187892997068);
		\draw [line width=2.pt,color=qqwwzz] (-5.801937735804838,1.5135187892997068)-- (-2.0104865064116564,-1.366661546270525);
		\draw [line width=2.pt,color=qqwwzz] (-2.0104865064116564,-1.366661546270525)-- (-3.,-2.);
		\draw [line width=2.pt,color=qqwwzz] (-5.,-2.)-- (-5.801937735804838,1.5135187892997068);
		\draw [line width=2.pt,color=qqwwzz] (-5.801937735804838,1.5135187892997068)-- (-5.967052680751909,-1.561019182904783);
		\draw [line width=2.pt,color=qqwwzz] (-5.967052680751909,-1.561019182904783)-- (-5.,-2.);
		\draw [line width=2.pt,color=wwzzqq] (-1.7530203962825333,-0.43633703506394084)-- (-1.,0.);
		\draw [line width=2.pt,color=wwzzqq] (-1.,0.)-- (-1.,-1.);
		\draw [line width=2.pt,color=wwzzqq] (-1.,-1.)-- (-1.7530203962825333,-0.43633703506394084);
		\draw [line width=2.pt,color=wwzzqq] (-6.246979603717467,-0.43633703506393995)-- (-6.786131292281997,0.);
		\draw [line width=2.pt,color=wwzzqq] (-6.786131292281997,0.)-- (-6.786131292281997,-1.0195943379950638);
		\draw [line width=2.pt,color=wwzzqq] (-6.786131292281997,-1.0195943379950638)-- (-6.246979603717467,-0.43633703506393995);
		\draw [line width=2.pt,color=qqwwzz] (4.,2.3812862675348225)-- (5.,-2.);
		\draw [line width=2.pt,color=qqwwzz] (5.,-2.)-- (5.139097984062835,2.3677816147734356);
		\draw [line width=2.pt,color=qqwwzz] (5.139097984062835,2.3677816147734356)-- (4.,2.3812862675348225);
		\draw [line width=2.pt,color=qqwwzz] (5.801937735804838,1.5135187892997062)-- (3.,-2.);
		\draw [line width=2.pt,color=qqwwzz] (3.,-2.)-- (6.416305310516531,0.7018590150512227);
		\draw [line width=2.pt,color=qqwwzz] (6.416305310516531,0.7018590150512227)-- (5.801937735804838,1.5135187892997062);
		\draw [line width=2.pt,color=qqwwzz] (6.246979603717467,-0.43633703506394084)-- (1.7530203962825333,-0.43633703506393995);
		\draw [line width=2.pt,color=qqwwzz] (1.7530203962825333,-0.43633703506393995)-- (6.013707348916996,-1.4777230529186725);
		\draw [line width=2.pt,color=qqwwzz] (6.013707348916996,-1.4777230529186725)-- (6.246979603717467,-0.43633703506394084);
		\draw [line width=2.pt,color=qqwwzz] (6.246979603717467,-0.43633703506394084)-- (2.198062264195162,1.5135187892997068);
		\draw [line width=2.pt,color=qqwwzz] (2.198062264195162,1.5135187892997068)-- (2.820689032782757,2.3122508614493618);
		\draw [line width=2.pt,color=qqwwzz] (2.820689032782757,2.3122508614493618)-- (6.246979603717467,-0.43633703506394084);
		\draw [line width=2.pt,color=qqwwzz] (5.,-2.)-- (2.198062264195162,1.5135187892997068);
		\draw [line width=2.pt,color=qqwwzz] (2.198062264195162,1.5135187892997068)-- (1.6128951479841538,0.7157417033822412);
		\draw [line width=2.pt,color=qqwwzz] (1.6128951479841538,0.7157417033822412)-- (5.,-2.);
		\draw [line width=2.pt,color=qqwwzz] (3.,-2.)-- (2.198062264195162,1.5135187892997068);
		\draw [line width=2.pt,color=qqwwzz] (2.198062264195162,1.5135187892997068)-- (4.097896359236453,-2.3662151061038528);
		\draw [line width=2.pt,color=qqwwzz] (4.097896359236453,-2.3662151061038528)-- (3.,-2.);
		\draw [line width=2.pt,color=qqwwzz] (1.7530203962825333,-0.43633703506393995)-- (2.198062264195162,1.5135187892997068);
		\draw [line width=2.pt,color=qqwwzz] (2.198062264195162,1.5135187892997068)-- (2.2098507462179464,-1.616549936228857);
		\draw [line width=2.pt,color=qqwwzz] (2.2098507462179464,-1.616549936228857)-- (1.7530203962825333,-0.43633703506393995);
		\draw [line width=2.pt,color=wwzzqq] (4.,2.3812862675348225)-- (3.3898792543545126,2.8259103296970443);
		\draw [line width=2.pt,color=wwzzqq] (3.3898792543545126,2.8259103296970443)-- (4.4449635675119135,3.006385278000284);
		\draw [line width=2.pt,color=wwzzqq] (4.4449635675119135,3.006385278000284)-- (4.,2.3812862675348225);
		\draw [line width=2.pt,color=wwzzqq] (5.801937735804838,1.5135187892997062)-- (6.055355413910052,2.3261335497803803);
		\draw [line width=2.pt,color=wwzzqq] (6.055355413910052,2.3261335497803803)-- (6.666193700474863,1.7291779515465873);
		\draw [line width=2.pt,color=wwzzqq] (6.666193700474863,1.7291779515465873)-- (5.801937735804838,1.5135187892997062);
		\begin{scriptsize}
			\draw [fill=uuuuuu] (-5.,-2.) circle (2.5pt);
			\draw[color=uuuuuu] (-5.2800516364817325,-2.1578803410819027) node {$w_5$};
			\draw [fill=uuuuuu] (-3.,-2.) circle (2.5pt);
			\draw[color=uuuuuu] (-2.883377062234357,-2.2081602272549246) node {$w_4$};
			\draw [fill=uuuuuu] (-1.7530203962825333,-0.43633703506394084) circle (2.5pt);
			\draw[color=uuuuuu] (-1.7101797181971805,-0.6997636420642677) node {$w_3$};
			\draw [fill=uuuuuu] (-2.198062264195162,1.5135187892997062) circle (2.5pt);
			\draw[color=uuuuuu] (-2.0453789593506597,1.7807107424714788) node {$w_2$};
			\draw [fill=uuuuuu] (-4.,2.3812862675348225) circle (2.5pt);
			\draw[color=uuuuuu] (-3.922494709810142,2.7192686177012204) node {$w_1$};
			\draw [fill=uuuuuu] (-5.801937735804838,1.5135187892997068) circle (2.5pt);
			\draw[color=uuuuuu] (-6.0677698531924085,1.7471908183561309) node {$w_7$};
			\draw [fill=uuuuuu] (-6.246979603717467,-0.43633703506393995) circle (2.5pt);
			\draw[color=uuuuuu] (-6.3359292461151915,-0.7997636420642677) node {$w_6$};
			\draw [fill=uuuuuu] (3.,-2.) circle (2.5pt);
			\draw[color=uuuuuu] (2.7814901132594385,-2.2249201893125985) node {$w_5$};
			\draw [fill=uuuuuu] (5.,-2.) circle (2.5pt);
			\draw[color=uuuuuu] (5.228444573679836,-2.2416801513702724) node {$w_4$};
			\draw [fill=uuuuuu] (6.246979603717467,-0.43633703506394084) circle (2.5pt);
			\draw[color=uuuuuu] (6.488401879774686,-0.565683945602876) node {$w_3$};
			\draw [fill=uuuuuu] (5.801937735804838,1.5135187892997062) circle (2.5pt);
			\draw[color=uuuuuu] (6.116722562736555,1.4022714633756735) node {$w_2$};
			\draw [fill=uuuuuu] (4.,2.3812862675348225) circle (2.5pt);
			\draw[color=uuuuuu] (3.6362481782008103,2.417589300663089) node {$w_1$};
			\draw [fill=uuuuuu] (2.198062264195162,1.5135187892997068) circle (2.5pt);
			\draw[color=uuuuuu] (1.8764521621450452,1.7807107424714788) node {$w_7$};
			\draw [fill=uuuuuu] (1.7530203962825333,-0.43633703506393995) circle (2.5pt);
			\draw[color=uuuuuu] (1.4239331865878486,-0.565683945602876) node {$w_6$};
			\draw [fill=ududff] (-2.829565117941744,2.409429679766491) circle (2.5pt);
			\draw [fill=ududff] (-5.147974069221822,2.3261335497803803) circle (2.5pt);
			\draw [fill=ududff] (-4.051241691071366,-2.3939804827658895) circle (2.5pt);
			\draw [fill=ududff] (-6.550125590654683,0.9378647166785362) circle (2.5pt);
			\draw [fill=ududff] (-1.635653921474159,0.5491494434100199) circle (2.5pt);
			\draw [fill=ududff] (-2.0104865064116564,-1.366661546270525) circle (2.5pt);
			\draw [fill=ududff] (-5.967052680751909,-1.561019182904783) circle (2.5pt);
			\draw [fill=ududff] (-1.,0.) circle (2.5pt);
			\draw [fill=ududff] (-1.,-1.) circle (2.5pt);
			\draw [fill=xdxdff] (-6.786131292281997,0.) circle (2.5pt);
			\draw [fill=ududff] (-6.786131292281997,-1.0195943379950638) circle (2.5pt);
			\draw [fill=ududff] (5.139097984062835,2.3677816147734356) circle (2.5pt);
			\draw [fill=ududff] (6.416305310516531,0.7018590150512227) circle (2.5pt);
			\draw [fill=ududff] (6.013707348916996,-1.4777230529186725) circle (2.5pt);
			\draw [fill=ududff] (2.820689032782757,2.3122508614493618) circle (2.5pt);
			\draw [fill=ududff] (1.6128951479841538,0.7157417033822412) circle (2.5pt);
			\draw [fill=ududff] (4.097896359236453,-2.3662151061038528) circle (2.5pt);
			\draw [fill=ududff] (2.2098507462179464,-1.616549936228857) circle (2.5pt);
			\draw [fill=ududff] (3.3898792543545126,2.8259103296970443) circle (2.5pt);
			\draw [fill=ududff] (4.4449635675119135,3.006385278000284) circle (2.5pt);
			\draw [fill=ududff] (6.055355413910052,2.3261335497803803) circle (2.5pt);
			\draw [fill=ududff] (6.666193700474863,1.7291779515465873) circle (2.5pt);
		\end{scriptsize}
	\end{tikzpicture}
	\caption{13 3-cliques in $K_{18}$ for case \eqref{eq:c5}}
	\label{fig:k18alpha13-9}
\end{figure}

\begin{lemma}\label{lem:k18a13c5-2}
In a $\{K_3,K_4\}$-decomposition of $K_{18}$, if $\alpha=13$, then the case \eqref{eq:c5} is impossible.
\end{lemma}
\begin{proof}
	Following Lemma \ref{lem:k18a13c5-1}, we continue to exclude the possibility that $H\cong G_5$. Assume that $H\cong G_5$, we first take the 4 3-cliques in $\overline{H}$: $\{w_1,w_2,w_7\}$, $\{w_1,w_3,w_5\}$, $\{w_2,w_3,w_4\}$, and $\{w_4,w_5,w_6\}$. Since $t_1=2$ in \eqref{eq:c5}, we discuss about which vertices in $H$ are contained in these 2 3-cliques. 
	\begin{itemize}
		\item If at least one of the degree-4 vertices $w_6,w_7$ in $H$ is contained in such a 3-clique, then, since they cannot be both, without loss of generality, we assume there is a 3-clique $\{w_6, y_1,y_2\}$. Then the other vertex contained in the $t_1=2$ 3-cliques can be one of several possibilities. 
		\begin{enumerate}
			\item It is $w_3$, the unique common neighbour of $w_6$ and $w_7$. This 3-clique can be taken to be $\{w_3,y_3,y_4\}$. We then deduce that all edges of $H$ except $\{w_3,w_6\}$ and $\{w_6,w_7\}$ contribute to the $t_2=7$ 3-cliques, which we may take to be $\{w_1,w_4,y_5\}$, $\{w_1,w_6,y_6\}$, $\{w_2,w_5,y_7\}$, $\{w_2,w_6,y_8\}$, $\{w_3,w_7,y_9\}$, $\{w_4,w_7,y_{10}\}$, $\{w_5,w_7,y_{11}\}$. This configuration (See Figure \ref{fig:k18alpha13-10} upper left) is then excluded by an exhaustive computer-assisted search. 
			\item It is adjacent to $w_6$ but not adjacent to $w_7$, and without loss of generality we assume it is $w_1$, and take a 3-clique $\{w_1,y_3,y_4\}$. We then deduce that all edges of $H$ except $\{w_1,w_6\}$ and $\{w_6,w_7\}$ contribute to the $t_2=7$ 3-cliques, which we may take to be $\{w_1,w_4,y_5\}$, $\{w_2,w_5,y_6\}$, $\{w_2,w_6,y_7\}$, $\{w_3,w_6,y_8\}$, $\{w_3,w_7,y_9\}$, $\{w_4,w_7,y_{10}\}$, $\{w_5,w_7,y_{11}\}$. This configuration (See Figure \ref{fig:k18alpha13-10} upper right) is then excluded by an exhaustive computer-assisted search. 
		\end{enumerate}
		\item Else if $w_3$ is one of them, we may take the 3-clique $\{w_3,y_1,y_2\}$. Then the other vertex contained in the $t_1=2$ 3-cliques has 2 possibilities. 
		\begin{enumerate}
			\item It coincides with $w_3$, i.e. $w_3$ is contained in another 3-clique $\{w_3, y_3, y_4\}$, and we can deduce that all edges of $H$ except $\{w_3,w_6\}$ and $\{w_3,w_7\}$ contribute to the $t_2=7$ 3-cliques, which we may take to be $\{w_1,w_4,y_5\}$, $\{w_1,w_6,y_6\}$, $\{w_2,w_5,y_7\}$, $\{w_2,w_6,y_8\}$, $\{w_4,w_7,y_9\}$, $\{w_5,w_7,y_{10}\}$, $\{w_6,w_7,y_{11}\}$. This configuration (See Figure \ref{fig:k18alpha13-10} lower left) is then excluded by an exhaustive computer-assisted search. 
			\item It is one of the neighbours $w_6$, $w_7$ of $w_3$. This case, up to isomorphism, has been already considered. 
		\end{enumerate}
		\item Otherwise, they must be two adjacent vertices among $\{w_1,w_2,w_4,w_5\}$. Without loss of generality, we take two 3-cliques $\{w_1,y_1,y_2\}$ and $\{w_4,y_3,y_4\}$. Then we deduce that all edges of $H$ except $\{w_1,w_4\}$ and $\{w_6,w_7\}$ contribute to the $t_2=7$ 3-cliques, which we may take to be $\{w_1,w_6,y_5\}$, $\{w_2,w_5,y_6\}$, $\{w_2,w_6,y_7\}$, $\{w_3,w_6,y_8\}$, $\{w_3,w_7,y_9\}$, $\{w_4,w_7,y_{10}\}$, $\{w_5,w_7,y_{11}\}$. This configuration (See Figure \ref{fig:k18alpha13-10} lower right) is then excluded by an exhaustive computer-assisted search. 
	\end{itemize}
        This concludes the proof.
\end{proof}

\begin{figure}
	\centering
	\definecolor{qqzzcc}{rgb}{0.,0.6,0.8}
	\definecolor{wqwqwq}{rgb}{0.3764705882352941,0.3764705882352941,0.3764705882352941}
	\definecolor{qqccqq}{rgb}{0.,0.8,0.}
	\definecolor{wwzzqq}{rgb}{0.4,0.6,0.}
	\definecolor{uuuuuu}{rgb}{0.26666666666666666,0.26666666666666666,0.26666666666666666}
	\begin{tikzpicture}[line cap=round,line join=round,>=triangle 45,x=0.3cm,y=0.3cm]
		\clip(-12.,-12.) rectangle (12.,12.);
		\fill[line width=2.pt,color=wqwqwq,fill=wqwqwq,fill opacity=0.10000000149011612] (-6.,10.762572535069644) -- (-2.3961245283903234,9.027037578599412) -- (-9.603875471609676,9.027037578599414) -- cycle;
		\fill[line width=2.pt,color=wqwqwq,fill=wqwqwq,fill opacity=0.10000000149011612] (-6.,10.762572535069644) -- (-1.5060407925650665,5.127325929872118) -- (-8.,2.) -- cycle;
		\fill[line width=2.pt,color=wqwqwq,fill=wqwqwq,fill opacity=0.10000000149011612] (-2.3961245283903234,9.027037578599412) -- (-1.5060407925650665,5.127325929872118) -- (-4.,2.) -- cycle;
		\fill[line width=2.pt,color=wqwqwq,fill=wqwqwq,fill opacity=0.10000000149011612] (-4.,2.) -- (-8.,2.) -- (-10.493959207434933,5.12732592987212) -- cycle;
		\fill[line width=2.pt,color=wqwqwq,fill=wqwqwq,fill opacity=0.10000000149011612] (2.3961245283903243,9.027037578599414) -- (6.,10.762572535069644) -- (9.603875471609676,9.027037578599412) -- cycle;
		\fill[line width=2.pt,color=wqwqwq,fill=wqwqwq,fill opacity=0.10000000149011612] (9.603875471609676,9.027037578599412) -- (10.493959207434933,5.127325929872118) -- (8.,2.) -- cycle;
		\fill[line width=2.pt,color=wqwqwq,fill=wqwqwq,fill opacity=0.10000000149011612] (6.,10.762572535069644) -- (10.493959207434933,5.127325929872118) -- (4.,2.) -- cycle;
		\fill[line width=2.pt,color=wqwqwq,fill=wqwqwq,fill opacity=0.10000000149011612] (8.,2.) -- (4.,2.) -- (1.5060407925650665,5.12732592987212) -- cycle;
		\fill[line width=2.pt,color=wqwqwq,fill=wqwqwq,fill opacity=0.10000000149011612] (-6.,-1.237427464930355) -- (-9.603875471609676,-2.9729624214005863) -- (-2.3961245283903234,-2.9729624214005876) -- cycle;
		\fill[line width=2.pt,color=wqwqwq,fill=wqwqwq,fill opacity=0.10000000149011612] (-2.3961245283903234,-2.9729624214005876) -- (-1.5060407925650665,-6.872674070127882) -- (-4.,-10.) -- cycle;
		\fill[line width=2.pt,color=wqwqwq,fill=wqwqwq,fill opacity=0.10000000149011612] (-6.,-1.237427464930355) -- (-1.5060407925650665,-6.872674070127882) -- (-8.,-10.) -- cycle;
		\fill[line width=2.pt,color=wqwqwq,fill=wqwqwq,fill opacity=0.10000000149011612] (-10.493959207434933,-6.87267407012788) -- (-8.,-10.) -- (-4.,-10.) -- cycle;
		\fill[line width=2.pt,color=wqwqwq,fill=wqwqwq,fill opacity=0.10000000149011612] (2.3961245283903243,-2.9729624214005863) -- (6.,-1.237427464930355) -- (9.603875471609676,-2.9729624214005876) -- cycle;
		\fill[line width=2.pt,color=wqwqwq,fill=wqwqwq,fill opacity=0.10000000149011612] (9.603875471609676,-2.9729624214005876) -- (10.493959207434933,-6.872674070127882) -- (8.,-10.) -- cycle;
		\fill[line width=2.pt,color=wqwqwq,fill=wqwqwq,fill opacity=0.10000000149011612] (6.,-1.237427464930355) -- (10.493959207434933,-6.872674070127882) -- (4.,-10.) -- cycle;
		\fill[line width=2.pt,color=wqwqwq,fill=wqwqwq,fill opacity=0.10000000149011612] (1.5060407925650665,-6.87267407012788) -- (4.,-10.) -- (8.,-10.) -- cycle;
		\draw [line width=2.pt,color=wqwqwq] (-6.,10.762572535069644)-- (-2.3961245283903234,9.027037578599412);
		\draw [line width=2.pt,color=wqwqwq] (-2.3961245283903234,9.027037578599412)-- (-9.603875471609676,9.027037578599414);
		\draw [line width=2.pt,color=wqwqwq] (-9.603875471609676,9.027037578599414)-- (-6.,10.762572535069644);
		\draw [line width=2.pt,color=wqwqwq] (-6.,10.762572535069644)-- (-1.5060407925650665,5.127325929872118);
		\draw [line width=2.pt,color=wqwqwq] (-1.5060407925650665,5.127325929872118)-- (-8.,2.);
		\draw [line width=2.pt,color=wqwqwq] (-8.,2.)-- (-6.,10.762572535069644);
		\draw [line width=2.pt,color=wqwqwq] (-2.3961245283903234,9.027037578599412)-- (-1.5060407925650665,5.127325929872118);
		\draw [line width=2.pt,color=wqwqwq] (-1.5060407925650665,5.127325929872118)-- (-4.,2.);
		\draw [line width=2.pt,color=wqwqwq] (-4.,2.)-- (-2.3961245283903234,9.027037578599412);
		\draw [line width=2.pt,color=wqwqwq] (-4.,2.)-- (-8.,2.);
		\draw [line width=2.pt,color=wqwqwq] (-8.,2.)-- (-10.493959207434933,5.12732592987212);
		\draw [line width=2.pt,color=wqwqwq] (-10.493959207434933,5.12732592987212)-- (-4.,2.);
		\draw [line width=2.pt,color=wqwqwq] (2.3961245283903243,9.027037578599414)-- (6.,10.762572535069644);
		\draw [line width=2.pt,color=wqwqwq] (6.,10.762572535069644)-- (9.603875471609676,9.027037578599412);
		\draw [line width=2.pt,color=wqwqwq] (9.603875471609676,9.027037578599412)-- (2.3961245283903243,9.027037578599414);
		\draw [line width=2.pt,color=wqwqwq] (9.603875471609676,9.027037578599412)-- (10.493959207434933,5.127325929872118);
		\draw [line width=2.pt,color=wqwqwq] (10.493959207434933,5.127325929872118)-- (8.,2.);
		\draw [line width=2.pt,color=wqwqwq] (8.,2.)-- (9.603875471609676,9.027037578599412);
		\draw [line width=2.pt,color=wqwqwq] (6.,10.762572535069644)-- (10.493959207434933,5.127325929872118);
		\draw [line width=2.pt,color=wqwqwq] (10.493959207434933,5.127325929872118)-- (4.,2.);
		\draw [line width=2.pt,color=wqwqwq] (4.,2.)-- (6.,10.762572535069644);
		\draw [line width=2.pt,color=wqwqwq] (8.,2.)-- (4.,2.);
		\draw [line width=2.pt,color=wqwqwq] (4.,2.)-- (1.5060407925650665,5.12732592987212);
		\draw [line width=2.pt,color=wqwqwq] (1.5060407925650665,5.12732592987212)-- (8.,2.);
		\draw [line width=2.pt,color=wqwqwq] (-6.,-1.237427464930355)-- (-9.603875471609676,-2.9729624214005863);
		\draw [line width=2.pt,color=wqwqwq] (-9.603875471609676,-2.9729624214005863)-- (-2.3961245283903234,-2.9729624214005876);
		\draw [line width=2.pt,color=wqwqwq] (-2.3961245283903234,-2.9729624214005876)-- (-6.,-1.237427464930355);
		\draw [line width=2.pt,color=wqwqwq] (-2.3961245283903234,-2.9729624214005876)-- (-1.5060407925650665,-6.872674070127882);
		\draw [line width=2.pt,color=wqwqwq] (-1.5060407925650665,-6.872674070127882)-- (-4.,-10.);
		\draw [line width=2.pt,color=wqwqwq] (-4.,-10.)-- (-2.3961245283903234,-2.9729624214005876);
		\draw [line width=2.pt,color=wqwqwq] (-6.,-1.237427464930355)-- (-1.5060407925650665,-6.872674070127882);
		\draw [line width=2.pt,color=wqwqwq] (-1.5060407925650665,-6.872674070127882)-- (-8.,-10.);
		\draw [line width=2.pt,color=wqwqwq] (-8.,-10.)-- (-6.,-1.237427464930355);
		\draw [line width=2.pt,color=wqwqwq] (-10.493959207434933,-6.87267407012788)-- (-8.,-10.);
		\draw [line width=2.pt,color=wqwqwq] (-8.,-10.)-- (-4.,-10.);
		\draw [line width=2.pt,color=wqwqwq] (-4.,-10.)-- (-10.493959207434933,-6.87267407012788);
		\draw [line width=2.pt,color=wqwqwq] (2.3961245283903243,-2.9729624214005863)-- (6.,-1.237427464930355);
		\draw [line width=2.pt,color=wqwqwq] (6.,-1.237427464930355)-- (9.603875471609676,-2.9729624214005876);
		\draw [line width=2.pt,color=wqwqwq] (9.603875471609676,-2.9729624214005876)-- (2.3961245283903243,-2.9729624214005863);
		\draw [line width=2.pt,color=wqwqwq] (9.603875471609676,-2.9729624214005876)-- (10.493959207434933,-6.872674070127882);
		\draw [line width=2.pt,color=wqwqwq] (10.493959207434933,-6.872674070127882)-- (8.,-10.);
		\draw [line width=2.pt,color=wqwqwq] (8.,-10.)-- (9.603875471609676,-2.9729624214005876);
		\draw [line width=2.pt,color=wqwqwq] (6.,-1.237427464930355)-- (10.493959207434933,-6.872674070127882);
		\draw [line width=2.pt,color=wqwqwq] (10.493959207434933,-6.872674070127882)-- (4.,-10.);
		\draw [line width=2.pt,color=wqwqwq] (4.,-10.)-- (6.,-1.237427464930355);
		\draw [line width=2.pt,color=wqwqwq] (1.5060407925650665,-6.87267407012788)-- (4.,-10.);
		\draw [line width=2.pt,color=wqwqwq] (4.,-10.)-- (8.,-10.);
		\draw [line width=2.pt,color=wqwqwq] (8.,-10.)-- (1.5060407925650665,-6.87267407012788);
		\draw [line width=2.pt,color=qqzzcc] (-6.,10.762572535069644)-- (-4.,2.);
		\draw [line width=2.pt,color=qqzzcc] (-6.,10.762572535069644)-- (-10.493959207434933,5.12732592987212);
		\draw [line width=2.pt,color=qqzzcc] (-2.3961245283903234,9.027037578599412)-- (-8.,2.);
		\draw [line width=2.pt,color=qqzzcc] (-2.3961245283903234,9.027037578599412)-- (-10.493959207434933,5.12732592987212);
		\draw [line width=2.pt,color=qqzzcc] (-1.5060407925650665,5.127325929872118)-- (-9.603875471609676,9.027037578599414);
		\draw [line width=2.pt,color=qqzzcc] (-4.,2.)-- (-9.603875471609676,9.027037578599414);
		\draw [line width=2.pt,color=qqzzcc] (-8.,2.)-- (-9.603875471609676,9.027037578599414);
		\draw [line width=2.pt,color=qqzzcc] (6.,10.762572535069644)-- (8.,2.);
		\draw [line width=2.pt,color=qqzzcc] (9.603875471609676,9.027037578599412)-- (4.,2.);
		\draw [line width=2.pt,color=qqzzcc] (9.603875471609676,9.027037578599412)-- (1.5060407925650665,5.12732592987212);
		\draw [line width=2.pt,color=qqzzcc] (10.493959207434933,5.127325929872118)-- (1.5060407925650665,5.12732592987212);
		\draw [line width=2.pt,color=qqzzcc] (10.493959207434933,5.127325929872118)-- (2.3961245283903243,9.027037578599414);
		\draw [line width=2.pt,color=qqzzcc] (8.,2.)-- (2.3961245283903243,9.027037578599414);
		\draw [line width=2.pt,color=qqzzcc] (4.,2.)-- (2.3961245283903243,9.027037578599414);
		\draw [line width=2.pt,color=qqzzcc] (-6.,-1.237427464930355)-- (-4.,-10.);
		\draw [line width=2.pt,color=qqzzcc] (-6.,-1.237427464930355)-- (-10.493959207434933,-6.87267407012788);
		\draw [line width=2.pt,color=qqzzcc] (-2.3961245283903234,-2.9729624214005876)-- (-8.,-10.);
		\draw [line width=2.pt,color=qqzzcc] (-2.3961245283903234,-2.9729624214005876)-- (-10.493959207434933,-6.87267407012788);
		\draw [line width=2.pt,color=qqzzcc] (-4.,-10.)-- (-9.603875471609676,-2.9729624214005863);
		\draw [line width=2.pt,color=qqzzcc] (-8.,-10.)-- (-9.603875471609676,-2.9729624214005863);
		\draw [line width=2.pt,color=qqzzcc] (-10.493959207434933,-6.87267407012788)-- (-9.603875471609676,-2.9729624214005863);
		\draw [line width=2.pt,color=qqzzcc] (6.,-1.237427464930355)-- (1.5060407925650665,-6.87267407012788);
		\draw [line width=2.pt,color=qqzzcc] (9.603875471609676,-2.9729624214005876)-- (4.,-10.);
		\draw [line width=2.pt,color=qqzzcc] (9.603875471609676,-2.9729624214005876)-- (1.5060407925650665,-6.87267407012788);
		\draw [line width=2.pt,color=qqzzcc] (10.493959207434933,-6.872674070127882)-- (1.5060407925650665,-6.87267407012788);
		\draw [line width=2.pt,color=qqzzcc] (10.493959207434933,-6.872674070127882)-- (2.3961245283903243,-2.9729624214005863);
		\draw [line width=2.pt,color=qqzzcc] (8.,-10.)-- (2.3961245283903243,-2.9729624214005863);
		\draw [line width=2.pt,color=qqzzcc] (4.,-10.)-- (2.3961245283903243,-2.9729624214005863);
		\begin{scriptsize}
			\draw [fill=uuuuuu] (-8.,2.) circle (2.5pt);
			\draw[color=uuuuuu] (-8.400575883233115,1.3944657098239397) node {$w_5$};
			\draw [fill=uuuuuu] (-4.,2.) circle (2.5pt);
			\draw[color=uuuuuu] (-3.85819423966104,1.397112216678323) node {$w_4$};
			\draw [fill=wwzzqq] (-1.5060407925650665,5.127325929872118) circle (2.5pt);
			\draw[color=wwzzqq] (-0.7037646279527632,5.304194590891306) node {$w_3$};
			\draw [fill=uuuuuu] (-2.3961245283903234,9.027037578599412) circle (2.5pt);
			\draw[color=uuuuuu] (-2.1693600388457814,9.613633586075068) node {$w_2$};
			\draw [fill=uuuuuu] (-6.,10.762572535069644) circle (2.5pt);
			\draw[color=uuuuuu] (-5.809088919913148,11.506292604230097) node {$w_1$};
			\draw [fill=uuuuuu] (-9.603875471609676,9.027037578599414) circle (2.5pt);
			\draw[color=uuuuuu] (-10.34999239291067,9.380690937686756) node {$w_7$};
			\draw [fill=wwzzqq] (-10.593959207434933,5.12732592987212) circle (2.5pt);
			\draw[color=wwzzqq] (-11.437652001795773,5.362430252988385) node {$w_6$};
			\draw [fill=uuuuuu] (4.,2.) circle (2.5pt);
			\draw[color=uuuuuu] (3.5086170156193113,1.3935835408724786) node {$w_5$};
			\draw [fill=uuuuuu] (8.,2.) circle (2.5pt);
			\draw[color=uuuuuu] (8.633355280162164,1.3962300477268618) node {$w_4$};
			\draw [fill=uuuuuu] (10.493959207434933,5.127325929872118) circle (2.5pt);
			\draw[color=uuuuuu] (11.0945461019482,5.508019408231079) node {$w_3$};
			\draw [fill=uuuuuu] (9.603875471609676,9.027037578599412) circle (2.5pt);
			\draw[color=uuuuuu] (9.885422015249338,9.555397923977988) node {$w_2$};
			\draw [fill=wwzzqq] (6.,10.762572535069644) circle (2.5pt);
			\draw[color=wwzzqq] (6.2456931341819715,11.360703448987403) node {$w_1$};
			\draw [fill=uuuuuu] (2.3961245283903243,9.027037578599414) circle (2.5pt);
			\draw[color=uuuuuu] (1.7680184769011303,9.52628009292945) node {$w_7$};
			\draw [fill=wwzzqq] (1.5060407925650665,5.12732592987212) circle (2.5pt);
			\draw[color=wwzzqq] (0.7488943902022677,5.624490732425235) node {$w_6$};
			\draw [fill=uuuuuu] (-8.,-10.) circle (2.5pt);
			\draw[color=uuuuuu] (-8.54616503847581,-10.519434175319717) node {$w_5$};
			\draw [fill=uuuuuu] (-4.,-10.) circle (2.5pt);
			\draw[color=uuuuuu] (-3.6834872533698064,-10.590316344271177) node {$w_4$};
			\draw [fill=red] (-1.5060407925650665,-6.872674070127882) circle (2.5pt);
			\draw[color=red] (-0.9628824590013022,-6.4011734906213436) node {$w_3$};
			\draw [fill=uuuuuu] (-2.3961245283903234,-2.9729624214005876) circle (2.5pt);
			\draw[color=uuuuuu] (-2.023770883603087,-2.3829128059229716) node {$w_2$};
			\draw [fill=uuuuuu] (-6.,-1.237427464930355) circle (2.5pt);
			\draw[color=uuuuuu] (-5.634381933621915,-0.5484894498650188) node {$w_1$};
			\draw [fill=uuuuuu] (-9.603875471609676,-2.9729624214005863) circle (2.5pt);
			\draw[color=uuuuuu] (-10.089410084048373,-2.3537949748744325) node {$w_7$};
			\draw [fill=uuuuuu] (-10.493959207434933,-6.87267407012788) circle (2.5pt);
			\draw[color=uuuuuu] (-11.337652001795773,-6.837940956349428) node {$w_6$};
			\draw [fill=uuuuuu] (4.,-10.) circle (2.5pt);
			\draw[color=uuuuuu] (3.595970508764928,-10.548552006368256) node {$w_5$};
			\draw [fill=wwzzqq] (8.,-10.) circle (2.5pt);
			\draw[color=wwzzqq] (8.080116490239924,-10.579434175319717) node {$w_4$};
			\draw [fill=uuuuuu] (10.493959207434933,-6.872674070127882) circle (2.5pt);
			\draw[color=uuuuuu] (10.99366393299674,-6.242937828524265) node {$w_3$};
			\draw [fill=uuuuuu] (9.603875471609676,-2.9729624214005876) circle (2.5pt);
			\draw[color=uuuuuu] (9.929098761822146,-2.2518825662045465) node {$w_2$};
			\draw [fill=wwzzqq] (6.,-1.237427464930355) circle (2.5pt);
			\draw[color=wwzzqq] (6.26025204970624,-0.41745921014659365) node {$w_1$};
			\draw [fill=uuuuuu] (2.3961245283903243,-2.9729624214005863) circle (2.5pt);
			\draw[color=uuuuuu] (1.7469882371827052,-2.2810003972530852) node {$w_7$};
			\draw [fill=uuuuuu] (1.5060407925650665,-6.87267407012788) circle (2.5pt);
			\draw[color=uuuuuu] (0.8152176436294593,-6.328378912999996) node {$w_6$};
		\end{scriptsize}
	\end{tikzpicture}
	\caption{The 13 3-cliques in $K_{18}$ for case \eqref{eq:c5}. Each blue edge represents a 3-clique with one unshown vertex; each green vertex represents one 3-clique with two unshown vertices, and each red vertex represents two such 3-cliques.}
	\label{fig:k18alpha13-10}
\end{figure}

\begin{lemma}\label{lem:k18a13c6}
	In a $\{K_3,K_4\}$-decomposition of $K_{18}$, if $\alpha=13$, then the case \eqref{eq:c6} is impossible.
\end{lemma}
\begin{proof}
For contradiction, suppose case \eqref{eq:c6} is possible. Then three 3-cliques lie entirely among $w_1,\ldots,w_7$, so the complement of their union in the induced $K_7$ is an even graph with twelve edges. An isomorphism-free computer enumeration gives exactly six even graphs $H$ of order 7 and size 12, shown in Figure \ref{fig:nauty-7-12-e}, where vertex $i$ represents $w_{i+1}$ for $i\in\{0,1,\ldots,6\}$. An exhaustive computer-assisted search excludes $H_4,H_5,H_6$. 
	
	If $H\cong H_1$, we first take the 3 3-cliques in $\overline{H}$: $\{w_1,w_2,w_3\}$, $\{w_2,w_4,w_6\}$, $\{w_3,w_4,w_5\}$. Since $t_2=9$ in \eqref{eq:c6}, we then need to remove 3 edges from $H_1$ that do not contribute to 3-cliques. Since $\deg_{H}(w_7)=6$, we need to remove at least 2 edges incident to $w_7$. Similarly, since $\deg_H(w_1)=\deg_H(w_5)=\deg_H(w_6)=4$, for each of them we need to remove at least 1 incident edge. Since $t_1=1$ in \eqref{eq:c6}, there are 2 cases.
	\begin{itemize}
		\item If the 3 edges $\{w_1,w_7\}$, $\{w_5,w_7\}$ and $\{w_6,w_7\}$ incident to $w_7$ are removed, then the remaining 9 edges generate the 3-cliques: $\{w_1,w_4,y_1\}$, $\{w_1,w_5,y_2\}$, $\{w_1,w_6,y_3\}$, $\{w_2,w_5,y_4\}$, $\{w_2,w_7,y_5\}$, $\{w_3,w_6,y_6\}$, $\{w_3,w_7,y_7\}$, $\{w_4,w_7,y_8\}$, $\{w_5,w_6,y_9\}$. We finally append one 3-clique for $w_7$: $\{w_7,y_{10}, y_{11}\}$. This configuration (See Figure \ref{fig:k18a13c6} upper left) is then excluded by an exhaustive computer-assisted search. 
		\item Otherwise, only 2 edges incident to $w_7$ are removed, and without loss of generality, we assume $\{w_1,w_7\}$ remains, which forces us to remove another edge incident to $w_1$. There are 2 possibilities.
		\begin{enumerate}
			\item If $\{w_1,w_4\}$ is removed, then we have to form a 3-clique for $w_4$: $\{w_4,y_1,y_2\}$. The remaining 9 edges generate the 3-cliques: $\{w_1,w_5,y_3\}$, $\{w_1,w_6,y_4\}$, $\{w_1,w_7,y_5\}$, $\{w_2,w_5,y_6\}$, $\{w_2,w_7,y_7\}$, $\{w_3,w_6,y_8\}$, $\{w_3,w_7,y_9\}$, $\{w_4,w_7,y_{10}\}$, $\{w_5,w_6,y_{11}\}$. This configuration (See Figure \ref{fig:k18a13c6} upper middle) is then excluded by an exhaustive computer-assisted search. 
			\item Otherwise, without loss of generality, we may assume $\{w_1,w_5\}$ is removed, then we have to form a 3-clique for $w_5$: $\{w_5,y_1,y_2\}$. The remaining 9 edges generate the 3-cliques: $\{w_1,w_4,y_3\}$, $\{w_1,w_6,y_4\}$, $\{w_1,w_7,y_5\}$, $\{w_2,w_5,y_6\}$, $\{w_2,w_7,y_7\}$, $\{w_3,w_6,y_8\}$, $\{w_3,w_7,y_9\}$, $\{w_4,w_7,y_{10}\}$, $\{w_5,w_6,y_{11}\}$. This configuration (See Figure \ref{fig:k18a13c6} upper right) is then excluded by an exhaustive computer-assisted search. 
		\end{enumerate}
	\end{itemize}
	
	If $H\cong H_2$, we first take the 3 3-cliques in $\overline{H}$: $\{w_1,w_2,w_3\}$, $\{w_3,w_4,w_5\}$, $\{w_3,w_6,w_7\}$. Since $t_2=9$ in \eqref{eq:c6}, we then need to remove 3 edges from $H_1$ that do not contribute to 3-cliques. Now that $H_2$ is actually a 4-regular graph of order 6, ignoring the isolated vertex $w_3$, we deduce that for each vertex other than $w_3$, precisely one of its incident edges must be removed. In other words, the 3 removed edges form a perfect matching of the 4-regular graph. Since a $K_6$ with 2 different perfect matchings removed is always isomorphic to $K_6-C_6$, we may take an arbitrary perfect matching: $\{w_1,w_4\}$, $\{w_2,w_6\}$, and $\{w_5,w_7\}$. Then the remaining 9 edges generate the 3-cliques: $\{w_1,w_5,y_1\}$, $\{w_1,w_6,y_2\}$, $\{w_1,w_7,y_3\}$, $\{w_2,w_4,y_4\}$, $\{w_2,w_5,y_5\}$, $\{w_2,w_7,y_6\}$, $\{w_4,w_6,y_7\}$, $\{w_4,w_7,y_8\}$, $\{w_5,w_6,y_9\}$. We finally append a 3-clique for $w_3$: $\{w_3,y_{10}, y_{11}\}$. This configuration (See Figure \ref{fig:k18a13c6} centre) is then excluded by an exhaustive computer-assisted search. 
	
	If $H\cong H_3$, we first take the 3 3-cliques in $\overline{H}$: $\{w_1,w_2,w_3\}$, $\{w_3,w_4,w_7\}$, $\{w_4,w_5,w_6\}$. Since $t_1=1$ in \eqref{eq:c6}, we need to choose one vertex that is contained in this 3-clique. Inspecting the structure of $H_3$, there are 3 cases for this vertex. 
	\begin{itemize}
		\item If this vertex is $w_7$, we take the 3-clique $\{w_7,y_1,y_2\}$, and 2 edges incident to $w_7$ are to be removed, of which one leads to either $w_1$ or $w_2$ and the other to either $w_5$ or $w_6$. This is because, for each vertex among $w_1,w_2,w_5,w_6$, precisely 1 incident edge is to remove. Without loss of generality, we remove $\{w_1,w_7\}$ and $\{w_6,w_7\}$, and the 3rd edges to remove is clearly $\{w_2,w_5\}$. Then the remaining 9 edges generate the 3-cliques: $\{w_1,w_4,y_3\}$, $\{w_1,w_5,y_4\}$, $\{w_1,w_6,y_5\}$, $\{w_2,w_4,y_6\}$, $\{w_2,w_6,y_7\}$, $\{w_2,w_7,y_8\}$, $\{w_3,w_5,y_9\}$, $\{w_3,w_6,y_{10}\}$, $\{w_5,w_7,y_{11}\}$. This configuration (See Figure \ref{fig:k18a13c6} lower left) is then excluded by an exhaustive computer-assisted search. 
		\item Else if this vertex is one of $w_1,w_2,w_5,w_6$, without loss of generality we assume it to be $w_1$ and take the 3-clique $\{w_1,y_1,y_2\}$. Then 2 edges incident to $w_1$ are to be removed, of which one is $\{w_1,w_7\}$ and the other is either $\{w_1,w_5\}$ or $\{w_1,w_6\}$, since each of $w_2,w_5,w_6,w_7$ must have precisely 1 incident edge removed. Without loss of generality, we remove $\{w_1,w_5\}$, $\{w_1,w_7\}$ and consequently also $\{w_2,w_6\}$. Then the remaining 9 edges generate the 3-cliques: $\{w_1,w_4,y_3\}$, $\{w_1,w_6,y_4\}$, $\{w_2,w_4,y_5\}$, $\{w_2,w_5,y_6\}$, $\{w_2,w_7,y_7\}$, $\{w_3,w_5,y_8\}$, $\{w_3,w_6,y_9\}$, $\{w_5,w_7,y_{10}\}$, $\{w_6,w_7,y_{11}\}$. This configuration (See Figure \ref{fig:k18a13c6} lower middle) is then excluded by an exhaustive computer-assisted search. 
		\item Otherwise, this vertex is $w_3$ or $w_4$. Without loss of generality, we assume it to be $w_3$, and take the 3-clique $\{w_3, y_1, y_2\}$. Then 1 edge incident to $w_3$ is to be removed, which we assume to be $\{w_3,w_6\}$. Then the other 2 edges to be removed can be assumed, without loss of generality, to be $\{w_1,w_7\}$ and $\{w_2,w_5\}$. Then the remaining 9 edges generate the 3-cliques: $\{w_1,w_4,y_3\}$, $\{w_1,w_5,y_4\}$, $\{w_1,w_6,y_5\}$, $\{w_2,w_4,y_6\}$, $\{w_2,w_6,y_7\}$, $\{w_2,w_7,y_8\}$, $\{w_3,w_5,y_9\}$, $\{w_5,w_7,y_{10}\}$, $\{w_6,w_7,y_{11}\}$. This configuration (See Figure \ref{fig:k18a13c6} lower right) is then excluded by an exhaustive computer-assisted search. 
	\end{itemize}
        This concludes the proof.
\end{proof}

\begin{figure}
	\centering
	\definecolor{ffttww}{rgb}{1.,0.2,0.4}
	\definecolor{ududff}{rgb}{0.30196078431372547,0.30196078431372547,1.}
	\definecolor{uuuuuu}{rgb}{0.26666666666666666,0.26666666666666666,0.26666666666666666}

	\caption{6 even graphs of order 7 and size 12}
	\label{fig:nauty-7-12-e}
\end{figure}

\begin{figure}
	\centering
	\definecolor{qqzzcc}{rgb}{0.,0.6,0.8}
	\definecolor{wqwqwq}{rgb}{0.3764705882352941,0.3764705882352941,0.3764705882352941}
	\definecolor{wwzzqq}{rgb}{0.4,0.6,0.}
	\definecolor{uuuuuu}{rgb}{0.26666666666666666,0.26666666666666666,0.26666666666666666}
	%
	\caption{13 3-cliques in $K_{18}$ for the case \eqref{eq:c6}; each blue edge stands for a 3-clique with an unshown vertex, each green vertex stands for a 3-clique with two unshown vertices}
	\label{fig:k18a13c6}
\end{figure}

\begin{figure}[tbp]
	\centering
	\definecolor{yqyqyq}{rgb}{0.5019607843137255,0.5019607843137255,0.5019607843137255}
	\definecolor{ffwwqq}{rgb}{1.,0.4,0.}
	\definecolor{ffqqqq}{rgb}{1.,0.,0.}
	\definecolor{zzttqq}{rgb}{0.6,0.2,0.}
	\definecolor{ududff}{rgb}{0.30196078431372547,0.30196078431372547,1.}
	%
	\caption{13 3-cliques in $K_{18}$ for case \eqref{eq:c7}}
	\label{fig:k18alpha13-3}
\end{figure}

\begin{figure}[tbp]
	\centering
	\definecolor{yqyqyq}{rgb}{0.5019607843137255,0.5019607843137255,0.5019607843137255}
	\definecolor{ffwwqq}{rgb}{1.,0.4,0.}
	\definecolor{zzttqq}{rgb}{0.6,0.2,0.}
	\definecolor{ududff}{rgb}{0.30196078431372547,0.30196078431372547,1.}
	%
	\caption{13 3-cliques in $K_{18}$ for case \eqref{eq:c7}}
	\label{fig:k18alpha13-4}
\end{figure}

\begin{lemma}\label{lem:k18a13c7}
	In a $\{K_3,K_4\}$-decomposition of $K_{18}$, if $\alpha=13$, then the case \eqref{eq:c7} is impossible.
\end{lemma}
\begin{figure}[tbp]
	\centering
	\definecolor{wwqqcc}{rgb}{0.4,0.,0.8}
	\definecolor{wwccff}{rgb}{0.4,0.8,1.}
	\definecolor{qqqqcc}{rgb}{0.,0.,0.8}
	\definecolor{zzttqq}{rgb}{0.6,0.2,0.}
	\definecolor{ududff}{rgb}{0.30196078431372547,0.30196078431372547,1.}
	%
	\caption{13 3-cliques in $K_{18}$ for case \eqref{eq:c7}}
	\label{fig:k18alpha13-2}
\end{figure}
\begin{proof}
For contradiction suppose the case \eqref{eq:c7} is possible. Then there are $t_3=2$ 3-cliques among $w_1,\cdots,w_7$. The two 3-cliques can either be independent or share a vertex. We then need to pick 11 edges among $w_1,\cdots,w_7$ to determine all the 3-cliques. 
	
	First assume that the two 3-cliques among $w_1,\cdots,w_7$ are independent, and without loss of generality let them be $T_1=\{w_1,w_2,w_3\}$ and $T_2=\{w_4, w_5, w_6\}$. Among the 11 edges to be picked, 4 are incident to $w_7$ and the remaining 7 go between vertices in~$T_1$ and $T_2$. We therefore deduce that the 4 edges incident to $w_7$ must be picked in a balanced way between $T_1$ and $T_2$. Without loss of generality, let the 3-cliques containing $w_7$ be $T_3=\{w_1,w_7, y_1\}$, $T_4=\{w_2,w_7,y_2\}$, $T_5=\{w_4,w_7,y_3\}$ and $T_6=\{w_5,w_7,y_4\}$. Now that $w_3$ and $w_6$ each needs 3 more 3-cliques, we can take five more 3-cliques: $T_7=\{w_3,w_4,y_5\}$, $T_8=\{w_3,w_5,y_6\}$, $T_9=\{w_3,w_6,y_7\}$, $T_{10}=\{w_1,w_6,y_8\}$ and $T_{11}=\{w_2,w_6,y_9\}$. Since $w_1,w_2,w_4,w_5$ each needs 1 more 3-clique, without loss of generality, the last two 3-cliques can be taken to be $T_{12}=\{w_1,w_4,y_{10}\}$ and $T_{13}=\{w_2,w_5,y_{11}\}$. We then exclude this configuration (See Figure \ref{fig:k18alpha13-2}) by an exhaustive computer-assisted search on the remaining edges. 
	
	Then assume that the two 3-cliques among $w_1,\cdots,w_7$ share a vertex, and without loss of generality let them be $T_1=\{w_1,w_2,w_3\}$ and $T_2=\{w_1,w_4,w_5\}$. Clearly, the other two 3-cliques containing $w_1$ can be taken to be $T_3=\{w_1,w_6,y_1\}$ and $T_4=\{w_1,w_7,y_2\}$. We distinguish two cases depending on whether $w_6$ and $w_7$ share a 3-clique. 
	\begin{itemize}
		\item If $w_6$ and $w_7$ share a 3-clique, we suppose, without loss of generality, that there is a 3-clique $T_5=\{w_6,w_7,y_3\}$. Vertices $w_6$ and $w_7$ are each in 2 more 3-cliques, and these 4 3-cliques must be distributed in a balanced way between $T_1$ and $T_2$. There are thus two possibilities up to isomorphism.
		\begin{enumerate}
			\item Each of $w_6$ and $w_7$ shares one 3-clique with a vertex of $T_1$ and another with a vertex of $T_2$. Without loss of generality, we assume the 3-cliques to be $\{w_2,w_6,y_4\}$, $\{w_4,w_6,y_5\}$, $\{w_3,w_7,y_6\}$ and $\{w_5,w_7,y_7\}$. Then the remaining 4 3-cliques are uniquely determined: $\{w_2,w_4,y_8\}$, $\{w_3,w_4,y_9\}$, $\{w_2,w_5,y_{10}\}$, $\{w_3, w_5,y_{11}\}$. We then exclude this configuration (See Figure \ref{fig:k18alpha13-3} left) by an exhaustive computer-assisted search.  
			\item Each of $w_6$ and $w_7$ only shares 2 3-cliques with a vertex from one of $T_1$ and $T_2$. Without loss of generality, we assume the 3-cliques to be $\{w_2,w_6,y_4\}$, $\{w_3,w_6,y_5\}$, $\{w_4,w_7,y_6\}$, $\{w_5,w_7,y_7\}$. Then the remaining 4 3-cliques are uniquely determined: $\{w_2,w_4,y_8\}$, $\{w_3,w_4,y_9\}$, $\{w_2,w_5,y_{10}\}$, $\{w_3, w_5,y_{11}\}$. We then exclude this configuration (See Figure \ref{fig:k18alpha13-3} right) by an exhaustive computer-assisted search.  
		\end{enumerate}
		\item If $w_6$ and $w_7$ do not share a 3-clique, then each of them needs 3 more 3-cliques, and without loss of generality we may take $\{w_2,w_6,y_3\}$, $\{w_3,w_6,y_4\}$, $\{w_4,w_6,y_5\}$ and $\{w_2,w_7,y_6\}$, $\{w_4,w_7,y_7\}$, $\{w_5,w_7,y_8\}$. The remaining 3 3-cliques are uniquely determined at this point: $\{w_2,w_5,y_9\}$, $\{w_3,w_4,y_{10}\}$ and $\{w_3,w_5,y_{11}\}$. We then exclude this configuration (See Figure \ref{fig:k18alpha13-4}) by an exhaustive computer-assisted search.  
	\end{itemize}
	In conclusion, the case \eqref{eq:c7} is impossible. 
\end{proof}

\begin{lemma}
	$D(18,\{3,4\})=33$, or equivalently, $\alpha=15$.
\end{lemma}
\begin{proof}
	According to Lemma \ref{lem:k18sol} and Lemma \ref{lem:k18alphage13}, it suffices to show that $D(18,\{3,4\})\neq 32$, or equivalently $\alpha\neq 13$. We have listed all possible cases \eqref{eq:c1}--\eqref{eq:c7} in the proof of Lemma \ref{lem:k18k3k4alpha13cases}, and have excluded all of them in Lemmas \ref{lem:k18a13c1}, \ref{lem:k18alpha13c2}, \ref{lem:k18a13c3}, \ref{lem:k18alpha13c4}, \ref{lem:k18a13c5-2}, \ref{lem:k18a13c6} and \ref{lem:k18a13c7}. Therefore, $\alpha\neq 13$ has been proven. 
\end{proof}

\subsection{Optimal decomposition of \texorpdfstring{$K_{19}$}{K19}}

\begin{lemma}\label{lem:alphax-k19}
	In a $\{K_3,K_4\}$-decomposition of $K_{19}$, for any vertex $x$ in $K_{19}$, $\alpha_x\equiv 0\pmod{3}$.
\end{lemma}
\begin{proof}
	Counting the degree of $x$ gives 
	$$2\alpha_x+3\beta_x=19-1=18\equiv 0\pmod{3},$$
	so $\alpha_x\equiv 0\pmod{3}$.
\end{proof}

First, we provide an upper bound. 

\begin{lemma}\label{lem:v19alphale13}
	$D(19,\{3,4\})\le 35$, or equivalently, $\alpha\le 13$.
\end{lemma}
	\begin{table}
			\centering
			\footnotesize
			\begin{tabular}{llllllll}
					1 8 9 10 & 2 9 11 16 & 4 5 9 18 & 7 8 14 18 &  & 1 2 3 & 3 7 12 & 7 11 17 \\
					1 11 12 13 & 2 10 12 18 & 4 12 16 19 & 7 9 13 19 &  & 1 4 6 & 4 8 11 &  \\
					1 14 15 16 & 3 4 10 17 & 5 8 16 17 & 9 12 15 17 &  & 1 5 7 & 4 13 14 &  \\
					1 17 18 19 & 3 6 9 14 & 5 10 13 15 & 10 11 14 19 &  & 2 8 13 & 5 12 14 &  \\
					2 4 7 15 & 3 8 15 19 & 6 7 10 16 &  &  & 2 14 17 & 6 8 12 &  \\
					2 5 6 19 & 3 13 16 18 & 6 11 15 18 &  &  & 3 5 11 & 6 13 17 & 
				\end{tabular}
			\caption{$K_{19}$ decomposed into 13 3-cliques and 22 4-cliques.}
			\label{tab:k19byk3k4}
	\end{table}
\begin{proof}
	The block list in Table~\ref{tab:k19byk3k4} is a $\{K_3,K_4\}$-decomposition of $K_{19}$ with 13 3-cliques and 22 4-cliques. 
\end{proof}

We now raise the lower bound in the same way.

\begin{lemma}\label{lem:k19alphaneq7}
	$D(19,\{3,4\})\ge 33$, or equivalently, $\alpha\ge 9$.
\end{lemma}
\begin{proof}
	By Lemma \ref{lem:app-k19-alpha7}, it suffices to prove $D(19,\{3,4\})\neq 32$, which is equivalent to showing $\alpha\neq 7$. 
	
	For contradiction, assume there is a $\{K_3,K_4\}$-decomposition of $K_{19}$ with seven 3-cliques (and thus 25 4-cliques). For any vertex $x$ of $K_{19}$, we know $\alpha_x\equiv 0\pmod{3}$ by Lemma \ref{lem:alphax-k19}. If there were some vertex $u$ with $\alpha_u\ge 6$, then there would be at least twelve neighbours of $u$, each lying in at least three 3-cliques, hence the total number of 3-cliques would satisfy
	\begin{align*}
		\alpha\ge \frac{1\cdot 6 + 12\cdot 3}{3}=14,
	\end{align*}
	a contradiction. Therefore, $\alpha_x\le 3$ for every vertex $x$ of $K_{19}$, so there are seven vertices, each lying in three 3-cliques. Counting degrees shows that the seven 3-cliques decompose a $K_7$ on these vertices. Thus the decomposition yields a $\operatorname{PBD}(19,\{4,7^*\},1)$, which is nonexistent by Theorem \ref{thm:pbd4}. Therefore, $D(19,\{3,4\})\neq 32$. 
\end{proof}

\begin{lemma}\label{lem:v19alphage11}
	$D(19,\{3,4\})\ge 34$, or equivalently, $\alpha\ge 11$.
\end{lemma}
\begin{proof}
	According to Lemma \ref{lem:k19alphaneq7}, it suffices to show that $D(19, \{3,4\})\neq 33$, which is equivalent to $\alpha\neq 9$. 
	
	For contradiction, assume there is a $\{K_3,K_4\}$-decomposition of $K_{19}$ with nine 3-cliques (and thus 24 4-cliques). For any vertex $x$ of $K_{19}$, we know $\alpha_x\equiv 0\pmod{3}$ by Lemma \ref{lem:alphax-k19}, and we can deduce that $\alpha_x\le 3$ for any vertex $x$ just as in the proof of Lemma \ref{lem:k19alphaneq7}. Therefore, there are nine vertices, each lying in three 3-cliques, while the remaining ten vertices of $K_{19}$ lie in no 3-clique. Let $G$ be the $K_9$-subgraph of $K_{19}$ on the vertices that are contained in 3-cliques, and let $H$ be the subgraph of $G$ formed by taking the union of the nine 3-cliques. Clearly, $H$ is a 6-regular subgraph of order 9 in $G$, which implies that its complement $\overline{H}$ is a 2-regular subgraph of order 9 in $G$. Therefore, we can classify the shape of $\overline{H}$ up to isomorphism: 
	\begin{align*}
		\overline{H}\in \{C_9, C_6+C_3, C_5+C_4, 3C_3\}.
	\end{align*}
	All four possibilities for $\overline{H}$ are excluded by exhaustive computer-assisted searches. 
\end{proof}

\begin{lemma}
	$D(19,\{3,4\})= 35$, or equivalently, $\alpha= 13$.
\end{lemma}
\begin{proof}
	According to Lemma \ref{lem:v19alphale13} and Lemma \ref{lem:v19alphage11}, it suffices to show that $D(19, \{3,4\})\neq 34$, which is equivalent to $\alpha\neq 11$.
	
	For contradiction, assume there is a $\{K_3, K_4\}$-decomposition of $K_{19}$ with eleven 3-cliques (and thus 23 4-cliques). For any vertex $x$ of $K_{19}$, we know $\alpha_x\equiv 0\pmod{3}$ by Lemma \ref{lem:alphax-k19}, and we can deduce that $\alpha_x\le 3$ for any vertex $x$ just as in the proof of Lemma \ref{lem:k19alphaneq7}. Therefore, there are eleven vertices, each lying in three 3-cliques, while the remaining eight vertices of $K_{19}$ lie in no 3-clique. Let $G$ denote the $K_{11}$-subgraph of $K_{19}$ on the vertices contained in 3-cliques, and let $H$ denote the subgraph of $G$ formed by taking the union of the eleven 3-cliques. Clearly, $H$ is a $6$-regular graph of $G$. Based on our assumption, the complement of $H$ in $K_{19}$, denoted by $\overline{H}$, should be decomposable into 4-cliques. 

	An isomorphism-free computer enumeration gives 266 mutually nonisomorphic 6-regular graphs on eleven vertices, one of which must be $H$. An exhaustive computer-assisted search checks the complement of each graph in $K_{19}$ and finds no decomposition into 4-cliques.  Therefore, $\alpha\neq 11$. 
\end{proof}

\section{Explicit decompositions for \texorpdfstring{$K_{18}$ and $K_{19}$}{K18 and K19}}\label{app:k18k19}

The following block lists use the vertex sets $\{1,\ldots,18\}$ and
$\{1,\ldots,19\}$, respectively.  In each case every edge occurs exactly
once.

\subsection*{\texorpdfstring{$K_{18}$: $15K_3+18K_4$}{K18: 15 K3 + 18 K4}}

\noindent\textbf{4-cliques:}
\begin{quote}\small\sloppy
$\{1,2,10,17\}$, $\{2,4,5,7\}$, $\{5,9,13,18\}$,
$\{1,3,5,15\}$, $\{2,6,14,16\}$, $\{5,11,16,17\}$,
$\{1,4,9,16\}$, $\{2,9,12,15\}$, $\{6,8,9,17\}$,
$\{1,6,7,13\}$, $\{3,12,16,18\}$, $\{6,11,15,18\}$,
$\{1,8,11,12\}$, $\{4,12,13,14\}$, $\{7,9,11,14\}$,
$\{2,3,11,13\}$, $\{5,6,10,12\}$, $\{8,10,13,16\}$.
\end{quote}
\noindent\textbf{3-cliques:}
\begin{quote}\small\sloppy
$\{1,14,18\}$, $\{3,14,17\}$, $\{7,10,18\}$,
$\{2,8,18\}$, $\{4,8,15\}$, $\{7,12,17\}$,
$\{3,4,6\}$, $\{4,10,11\}$, $\{7,15,16\}$,
$\{3,7,8\}$, $\{4,17,18\}$, $\{10,14,15\}$,
$\{3,9,10\}$, $\{5,8,14\}$, $\{13,15,17\}$.
\end{quote}

\subsection*{\texorpdfstring{$K_{19}$: $13K_3+22K_4$}{K19: 13 K3 + 22 K4}}

\noindent\textbf{4-cliques:}
\begin{quote}\small\sloppy
$\{1,8,9,10\}$, $\{2,9,11,16\}$, $\{4,5,9,18\}$,
$\{7,8,14,18\}$, $\{1,11,12,13\}$, $\{2,10,12,18\}$,
$\{4,12,16,19\}$, $\{7,9,13,19\}$, $\{1,14,15,16\}$,
$\{3,4,10,17\}$, $\{5,8,16,17\}$, $\{9,12,15,17\}$,
$\{1,17,18,19\}$, $\{3,6,9,14\}$, $\{5,10,13,15\}$,
$\{10,11,14,19\}$, $\{2,4,7,15\}$, $\{3,8,15,19\}$,
$\{6,7,10,16\}$, $\{2,5,6,19\}$, $\{3,13,16,18\}$,
$\{6,11,15,18\}$.
\end{quote}
\noindent\textbf{3-cliques:}
\begin{quote}\small\sloppy
$\{1,2,3\}$, $\{3,7,12\}$, $\{7,11,17\}$,
$\{1,4,6\}$, $\{4,8,11\}$, $\{1,5,7\}$,
$\{4,13,14\}$, $\{2,8,13\}$, $\{5,12,14\}$,
$\{2,14,17\}$, $\{6,8,12\}$, $\{3,5,11\}$,
$\{6,13,17\}$.
\end{quote}

\end{document}